\documentclass[12pt]{article}
\usepackage{amsfonts}
\usepackage{epsfig}
\title{Outer Billiards, the Arithmetic Graph, and the Octagon}
\author{Richard Evan Schwartz \thanks{\hskip 5 pt Supported by 
N.S.F. Research Grant DMS-0072607}}

\newtheorem{theorem}{Theorem}[section]

\newtheorem{lemma}[theorem]{Lemma}

\newtheorem{corollary}[theorem]{Corollary}

\def\startproof{{\bf {\medskip}{\noindent}Proof: }}

\def\endproof{$\spadesuit$  \newline}

\def\C{\mbox{\boldmath{$C$}}}% 
\def\Q{\mbox{\boldmath{$Q$}}}% 
\def\R{\mbox{\boldmath{$R$}}}% 
\def\Z{\mbox{\boldmath{$Z$}}}% 

\begin{document}
\maketitle

\section{Introduction}

B. H. Neumann [{\bf N\/}] introduced outer billiards in
the late 1950s and J. Moser [{\bf M1\/}]
popularized the system in the 1970s as a toy model
for celestial mechanics.
Outer billiards is a discrete self-map of $\R^2-P$,
where $P$ is a bounded convex planar set as
in Figure 1.1 below.
Given $p_1 \in \R^2-P$, one
defines $p_2$ so that
the segment $\overline{p_1p_2}$ is tangent to $P$ at its
midpoint and $P$ lies to the right of the ray
$\overrightarrow{p_1p_2}$.  The map $p_1 \to p_2$
is called {\it the outer billiards map\/}.
The map is almost everywhere defined and
invertible.   See [{\bf T1\/}] for a survey of 
outer billiards.

\begin{center}
\resizebox{!}{2.2in}{\includegraphics{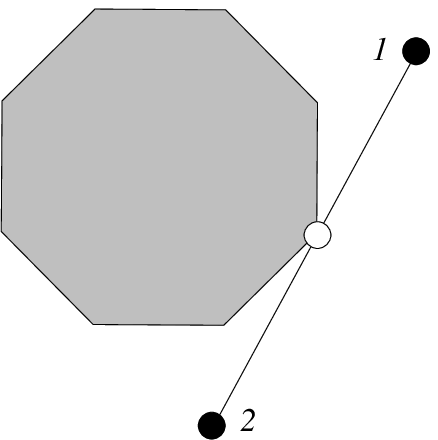}}
\newline
{\bf Figure 1.1:\/} Outer billiards relative to $P$.
\end{center}

Usually, the main object of study is the orbit $\{p_n\}$ of the
point $p_1$.   These orbits are both complicated
and interesting.   When $P$ is a convex polygon,
the orbits frequently have a fractal structure.
Aside from a few cases, the structure of the
orbits is not yet well understood. 

It turns out that it is sometimes 
productive to study not the orbit itself
but rather a certain {\it acceleration\/} of the
orbit.   Roughly speaking,  the acceleration we have
in mind amounts to considering
the first return map to a certain infinite
strip in the plane.   We call this acceleration
{\it the pinwheel dynamics\/}, because it is
based on the pinwheel map that we studied
in [{\bf S3\/}].   We give a full account
of the pinwheel map in \S \ref{pm} and
recall some of the results of [{\bf S3\/}] in
\S \ref{discuss}.   We point out, however,
that this paper does not rely on any of
the results in [{\bf S3\/}].  The work
here is self-contained.

The first purpose of this paper is to introduce
an object called the {\it arithmetic graph\/}.  The arithmetic
graph is a polygonal path $\Gamma(P,p) \subset \R^{n}$ that one associates
to the pinwheel dynamics of $p=p_1$.
One can view this graph as a geometric incarnation of the
symbolic coding of the pinwheel dynamics.
The arithmetic graph is quite similar to the lattice paths studied
by Vivaldi et. al. [{\bf V\/}] in connection with interval
exchange transformations.    We will define
the arithmetic graph in \S \ref{ag0}, right after
defining the pinwheel map.

In case $P$ is a kite -- i.e. a convex quadrilateral having a diagonal
that is a line of symmetry -- we have studied the arithmetic graph
in great detail.
In [{\bf S1\/}] we analyzed the graph
$\Gamma(P,p)$, where $P$ is the
Penrose kite and $p$ is a specially chosen point.
For this pair, we showed that $\Gamma(P,p)$ is
{\it coarsely self-similar\/}, in the sense that a
certain rescaled limit of $\Gamma(P,p)$ is a
self-similar fractal curve.  See \S \ref{lim} for
a formal definition of what we mean by a
{\it rescaled limit\/}.   

We used the coarse self-similarity of the
arithmetic graph in this case
to conclude that the orbit of $p$ is unbounded
relative to $P$.  This provided the first example
of an outer billiards system with unbounded orbits.
In [{\bf S2\/}] we studied the arithmetic graphs
relative to kites in general, and 
showed that outer billiards has unbounded orbits
when defined relative to any irrational kite.
In \S \ref{penrose} we will show some of
the nicest pictures of arithmetic graphs
associated to kites.

The irrational kites are
the only known examples \footnote{It is worth
mentioning that Dolgopyat and Fayad show in
[{\bf DF\/}] that outer billiards also has unbounded
orbits when defined relative to a semi-disk.}
 of polygons for which
outer billiards has unbounded orbits.
Since the arithmetic graph turns out to be a decisive
tool for establishing the unboundedness 
results for kites, we
think that it will also be a useful tool for
polygonal outer billiards more generally.
Compare Theorem \ref{unbound0} (a
result quoted from [{\bf S3\/}])
and the discussion following it.

Even in the case when all the orbits are known
to be bounded, as they are for so-called
{\it quasi-rational\/} polygons -- see
[{\bf VS\/}], [{\bf K\/}], [{\bf GS\/}] -- we think
that the arithmetic graph should shed light on the
dynamics.  
(See the end of \S \ref{discuss} for a definition
of {\it quasi-rational\/}.)
 In particular, the arithmetic graph
sheds light on the dynamics relative to the
regular polygons.  The regular polygons are
nice examples of quasi-rational polygons.

The vertices of a regular polygon have
coordinates that lie in a cyclotomic field.
When $p$ lies in this same field, one can view
certain projections of $\Gamma(p)$ as
``Galois conjugates'' of the orbits.
In somewhat the same way that the
full list of Galois conjugates of an algebraic number
sheds light on the number, the arithmetic graph
sheds light on the orbits associated to a polygon
with algebraic vertices.

This brings us to the main purpose of our paper.
We will concentrate on outer billiards for the
regular octagon, a polygon that we scale so
its vertices are the $8$th roots of unity.
The outer billiards dynamics for the regular
octagon are well understood from certain
points of view, but we will show that 
the study of the arithmetic graphs in
the regular octagon case leads to a big
surprise:  There is a
fractal loop $\Gamma$ in $\R^4$ which has the following
three properties.

\begin{itemize}

\item One planar projection $\pi_3$ of $\Gamma$ is an
embedded fractal loop reminiscent of the Koch
snowflake.  See Figure 1.2 below.  We call
this fractal {\it the snowflake\/}.  We give
a precise definition in \S \ref{snowflake}.

\item Another planar projection $\pi_2$ of $\Gamma$
is a fractal set reminiscent of the Sierpinski
carpet.  See Figure 1.3 below for an approximate
picture.  We call this fractal {\it the carpet\/}.
We give a precise definition in \S \ref{carpet}.

\item $\Gamma$ is the rescaled limit of the arithmetic
graphs associated to any sequence of {\it odd\/} periodic
orbits.
We give a precise definition of {\it rescaled limit\/} in
\S \ref{lim}.
\end{itemize}

It turns out that all the periodic orbits intersect a certain
strip in $3^k$ points for some integer $k$.   See
\S \ref{renorm}. We call the
orbit {\it odd\/} if $k$ is odd.
In Figure 1.2, the odd orbits are 
arithmetically closed orbits are lightly colored.
What we actually prove is slightly weaker than what we
have just said.  For technical reasons, we only consider
odd orbits that lie outside the first layer of big grey
octagons in Figure 1.3.

\begin{center}
\resizebox{!}{2.5in}{\includegraphics{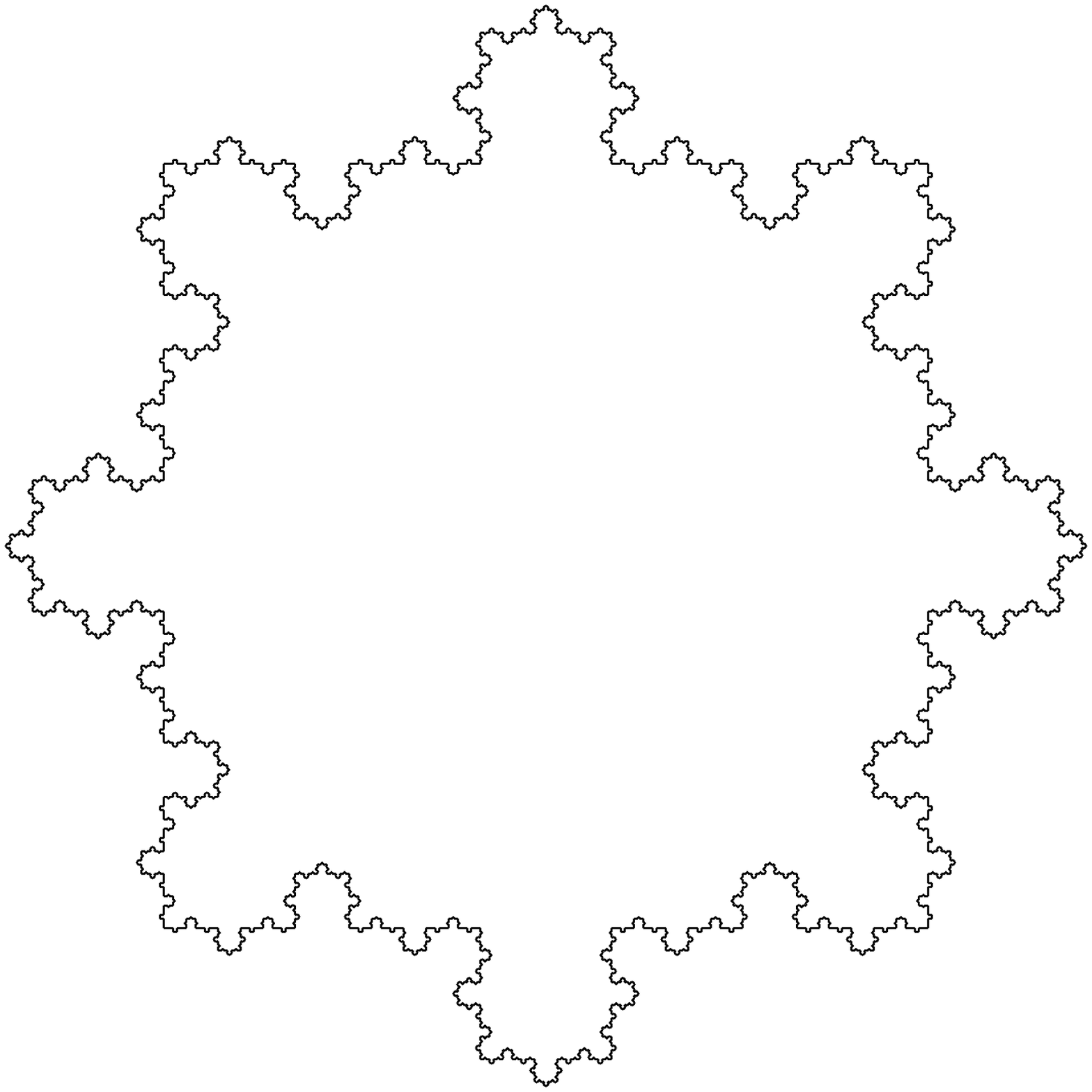}}
$\hskip 10 pt$
\resizebox{!}{2.5in}{\includegraphics{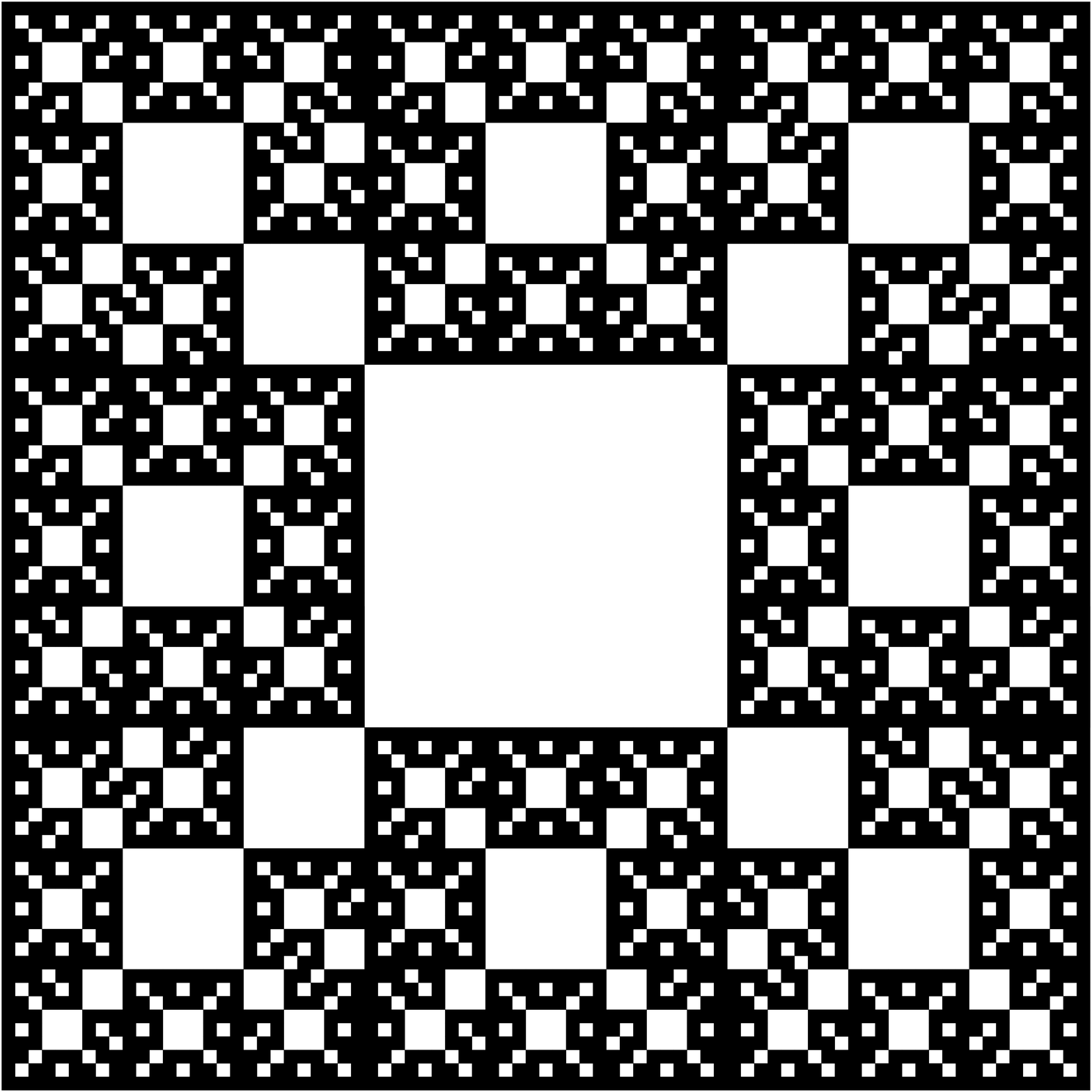}}
\newline
{\bf Figure 1.2:\/} The snowflake and the carpet
\end{center}

\begin{center}
\resizebox{!}{3.3in}{\includegraphics{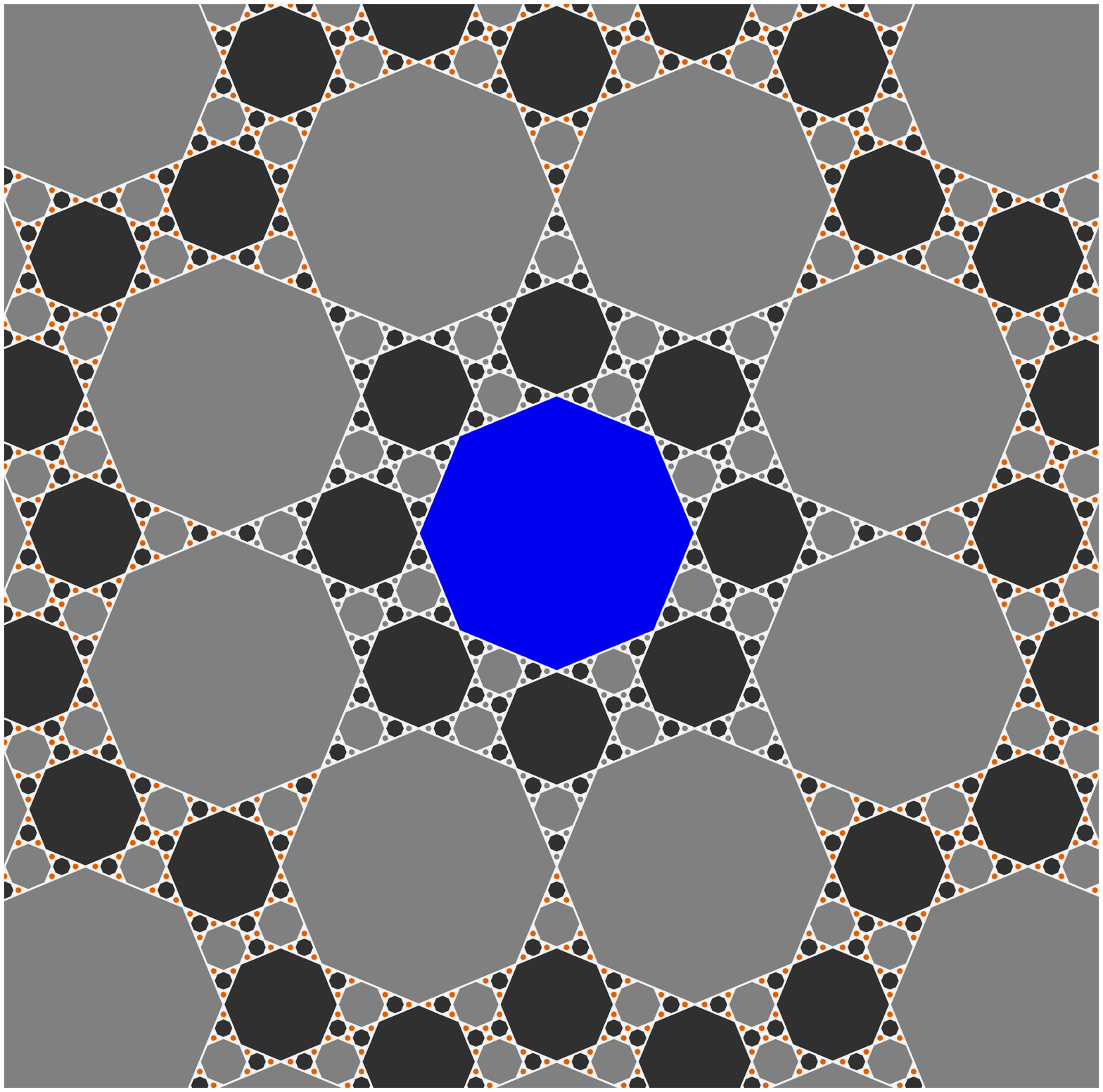}}
\newline
{\bf Figure 1.3:\/} Periodic orbits associated to the octagon
\end{center}

By considering the correspondence
between our two fractals, we produce a surjective
continuous map from the snowflake
to the carpet.
This map is ``symmetry respecting'' in a
certain formal sense that we discuss
in \S \ref{fractal}.  See \S \ref{map}
for a description of the map from
the snowflake to the carpet.
What is surprising is that this canonical
map between two seemingly unrelated fractals
actually arises ``naturally'', in connection
with outer billiards on the octagon.

For the sake of comparison, we now discuss the
situation for other regular polygons.For $n=3,4,6$,
 outer billiards on the regular $n$-gon
is rather trivial and easy to describe.  In all these
cases, there is a familiar periodic tiling of the plane,
such that the dynamics permutes the tiles.  In the
case $n=4$ we get the usual square tiling.  In
the cases $n=3,6$ we get a tiling by hexagons
and triangles.

For the case $n=5$,
Tabachnikov [{\bf T2\/}] analyzes the situation in detail.
The recent paper [{\bf BC\/}] builds on the work
in [{\bf T2\/}] and describes the
symbolic dynamics for the periodic orbits in great detail.
We are interested in something different, namely
the pinwheel dynamics, but the two are closely related.
In the case of the regular pentagon, our construction
produces a pair of embedded fractal curves, both akin to
the Koch snowflake, that arise as geometric
limits of arithmetic graphs associated to periodic orbits.
See \S \ref{penta}.
For the sake of brevity, we do not present proofs
of the statements we make for the regular pentagon.
We are mainly interested in the octagon, and
the analysis we give in the case of the octagon
could be fairly easily redone for the
pentagon.

The four cases $n=5,10,8,12$ are all quite similar
to each other, and possible to understand completely.
In these cases,
outer billiards on the regular $n$-gon has an efficient
renormalization scheme that allows one to give 
(at least in principle) a complete description of what is going on.
The special properties of these values of $n$ is that the
corresponding $n$th roots of unity lie in quadratic
number fields.   In the case $n=10$,  we get the
same two snowflakes as in the case $n=5$.
We leave the case $n=12$, which is mildly
more complicated than any of the other
quadratic cases, to the experimentally minded reader.

For any other positive integer $n$, outer billiards
on the regular $n$-gon is not well understood at all.
For the sake of comparison, we will draw some
pictures of arithmetic graphs associated to the
regular $7$-gon.  See
\S \ref{hept}.  In this case, the pictures reveal
an explosion of complexity in a geometrically striking
way.   Our $4$ figures in \S \ref{hept}
give only the faintest hint of the complexity.
\newline

Here is the plan of the paper.

\begin{itemize}
\item {\bf Basic Definitions:\/}
In \S 2 we will define the arithmetic graph,
and give explicit formulas in the case of
the regular octagon.
\item
{\bf A Picture Gallery:\/} In \S 3 we draw the arithmetic graphs
for many examples, including the regular
octagon.    After reading \S \ref{octopix},
the reader should have a good intuitive
idea of what our main theorem will say,
even though we defer the formal statement
until \S \ref{main result}.
\item
{\bf Carpet and Snowflake:\/}
In \S 4 we define the snowflake and
carpet fractals precisely, and discuss
the relationship between them.
\item
{\bf Reduction to a Small Region:\/}
We state our main result,
Theorem \ref{main}, at the beginning of \S 5.
In \S 5 we reduce Theorem \ref{main} to
the study of periodic points in a very small
region $R_1$.
\item
{\bf Toy Example of Renormalization:\/}
In \S 6 we study a simple and well-known
polygon exchange map that arises in connection
with outer billiards on the regular octagon.
The material in this chapter is not strictly necessary
for our formal proof, but the system we study here is
closely related to the one that comes from the 
pinwheel map.
\item
{\bf Pinwheel Dynamics:\/}
In \S 7, we describe the dynamics of the
pinwheel map on the region $R_1$.  Actually,
we reduce an even smaller region, $R$, which
is the top half of $R_1$.   The dynamical system on $R$
is the main dynamical system of interest to us.
Both the pinwheel map and the system 
from \S 6 have a period
$3$ renormalization that is the key to their
analysis.
\item
{\bf Substitution Scheme:\/}
In \S 8 we study the arithmetic graphs
produced by the dynamics on $R$.
Using the renormalization scheme, we describe
an efficient $2$-part method for generating the 
projections of the arithmetic graph of interest to us.
The first part of the method is combinatorial,
and involves a substitution scheme for
numerical sequences.  The second part is
geometrical, where we replace each number
in the created string by a certain vector.
We call this second half the {\it vector assignment}.
\item
{\bf Fixed Point of Renormalization:\/}
In \S 9 we modify the construction presented in
\S 8, keeping the combinatorial part the same
but changing the vector assignment.  We see that
the substitution scheme produces a kind of
{\it renormalization operator\/} defined on
the vector space of all possible vector
assignments.   This fixed point corresponds to
the Perron-Frobenius eigenvector.
 By taking the fixed point of
our operator as the new vector assignment,
we produce a nicer family of curves that
has the same rescaled limit.    
\item
{\bf Self-Similar Pattern Matching:\/}
In \S 10 we check by a combination of direct
calculation and induction that our improved
curves from \S 9 have the snowflake and
the carpet as their limit.   The key idea is
that the improved curves and the fractals
from \S 4 are self-similar in compatible ways.
 The final analysis completes the proof of the
Main Result.
\end{itemize}

The reader might be interested to know that we have
made two java programs, {\it OctoMap 1\/} and
{\it OctoMap 2\/}, which illustrate the mathematics
in this paper.  OctoMap 1 illustrates the snowflake
and carpet fractals in great detail, and gives an
interactive demonstration of the material in
\S 4.   OctoMap 2, a much more extensive
program, starts with the dynamical system produced
in \S 7 and illustrates all the key ideas that
go into the proof of the Main Theorem.
In particular, one can use OctoMap 2 to survey all
the computations we describe in \S 7-10.
One can find our applets, respectively, at the following address:
\begin{itemize}
\item http://www.math.brown.edu/$\sim$res/Java/OctoMap/Main.html
\item http://www.math.brown.edu/$\sim$res/Java/OctoMap2/Main.html
\end{itemize}
We strongly encourage the reader to look at these applets while reading
this paper.  They relate to this paper the same way that a cooked
meal relates to a recipe.
\newline

I'd like to thank Gordon Hughes and Sergei Tabachnikov about interesting
and inspiring conversations about outer billiards on regular polygons.
\newpage
\section{The Arithmetic Graph}
\label{pinwheel}

\subsection{Strip Maps}
\label{smap}

Let $\Sigma \subset \R^2$ be an infinite strip.
Let $V$ be a vector
whose tail lies on one edge of $\Sigma$ and whose head lies on the
other.   All this is shown in Figure 2.1.

\begin{center}
\resizebox{!}{3in}{\includegraphics{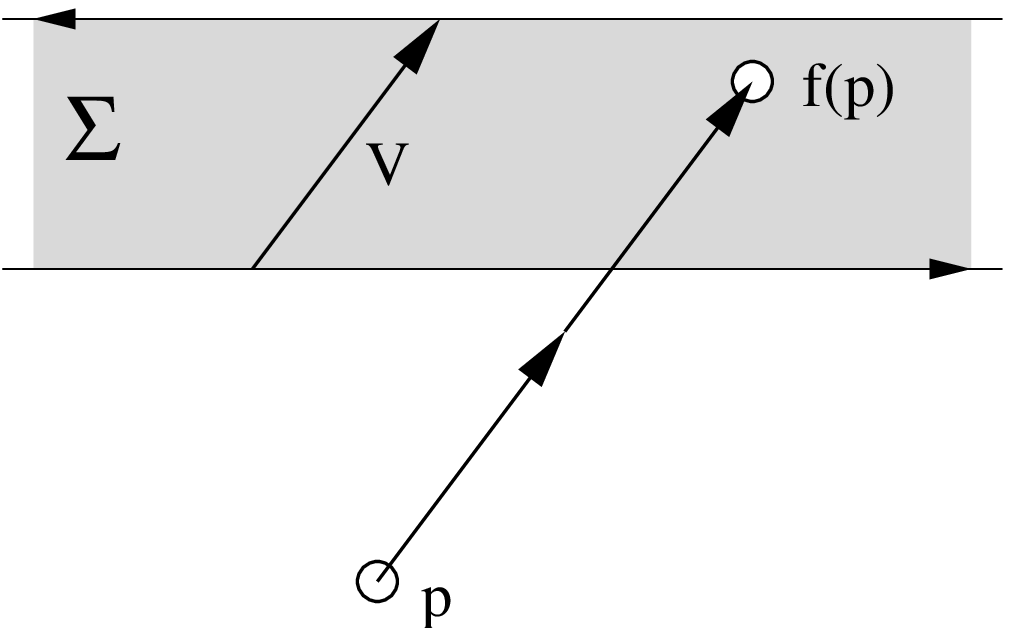}}
\newline
{\bf Figure 2.1:\/} A strip map
\end{center}  

To $(\Sigma,V)$ we associate 
a {\it strip map\/}
$f: \R^2 \to \Sigma$, defined by
$f(p)=p+nV$, 
where $n \in \Z$ is chosen so that $p+nV \in \Sigma$.
In the example shown, we have $n=2$.
\newline
\newline
{\bf Remarks:\/} \newline
(i) $f$ is not everywhere defined.  It is not defined
on a certain countable family of lines that are parallel
to $\Sigma$.  In particular $f$ is not defined on
the boundary of $\Sigma$.
\newline
(ii)  The pair $(\Sigma,-V)$ gives rise to the same
map as the pair $(\Sigma,V)$.  \newline
(iii) In case $\Sigma$ is the horizontain strip
bounded by the lines $y=0$ and $y=1$, and
$V=(0,1)$, the map
map $f$ has the formula $f(x,y)=(x,[y])$, where
$[y]$ is the fractional part of $y$.  In this case,
$f$ is not defined on the horizontal lines of
integer height.
Any other strip map is conjugate to this one by some affine
transformation.

\subsection{The Pinwheel Map}
\label{pm}

Let $P$ be a convex $n$-gon.  To $P$ we associate
$n$ special {\it pinwheel strip\/}.  Each pinwheel strip $\Sigma$ is
such that one component of $\partial \Sigma$ contains
an edge $e$ of $P$, and vertices of $P$ farthest
from this boundary component lie on the centerline
of $\Sigma$.    See Figure 2.2

\begin{center}
\resizebox{!}{3in}{\includegraphics{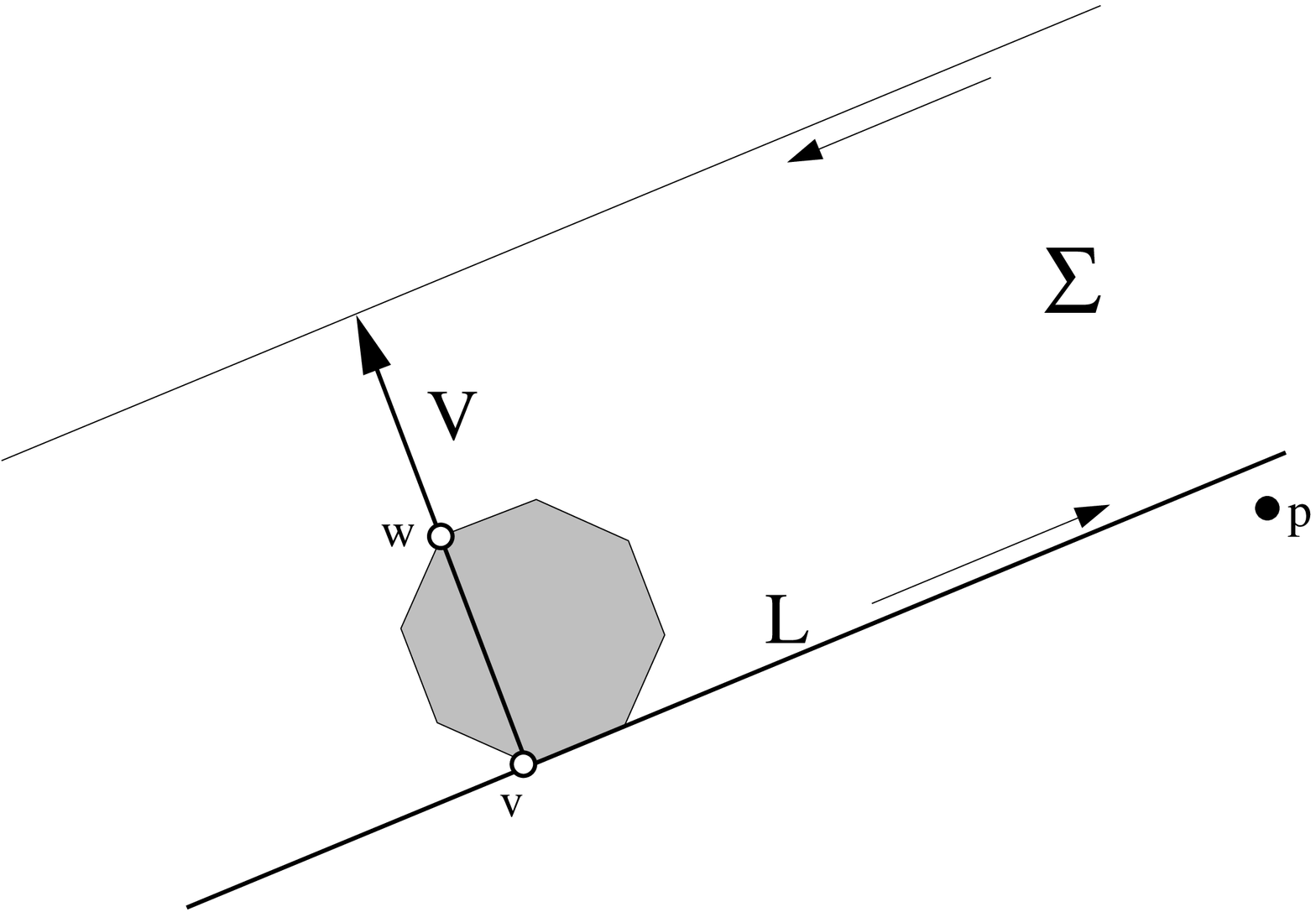}}
\newline
{\bf Figure 2.2:\/} A pinwheel strip associated to a regular octagon
\end{center}  

We orient the boundaries of a pinwheel strip $\Sigma$
so that a person walking along a boundary component
would see $\Sigma$ on the left.  Thus, the two boundary
components are oriented in opposite directions.
Say that a {\it pointed strip\/} is a strip, together with a choice
of boundary component.   We think of a pointed strip as
pointing in the direction of the orientation on its
preferred boundary component.   To each $n$-gon we
associate $2n$ pointed strips, $2$ per pinwheel strip.
We denote pointed strips by a pair $(\Sigma,L)$, where
$L$ is a boundary component of the strip $\Sigma$.

To the pair $(\Sigma,L)$ we associate a vector $V$ as
follows.  We choose $\epsilon>0$ and let
$p$ be a point which is $\epsilon$ units from $L$,
outside of $\Sigma$, and $1/\epsilon$ units away
from the origin.  Of the two possible locations
for $p$ that the above conditions determine,
we choose the one toward which $L$ points.
Figure 2.2 shows what we mean.   We then
let $V=\phi^2(p)-p$.  Here $\phi^2$ is the
square of the outer billiards map.  Our definition is independent
of $\epsilon$, provided that $\epsilon$ is sufficiently small.

In general, 
$V=2(w-v)$, 
where $v$ is the vertex bisecting the segment
joining $p$ to $\phi(p)$ and $w$ is the
vertex bisecting the segment joining $\phi(p)$
to $\phi^2(p)$. 
It is not hard to check that $\Sigma$ and $V$ are
related as in \S \ref{smap}.

\begin{center}
\resizebox{!}{4.5in}{\includegraphics{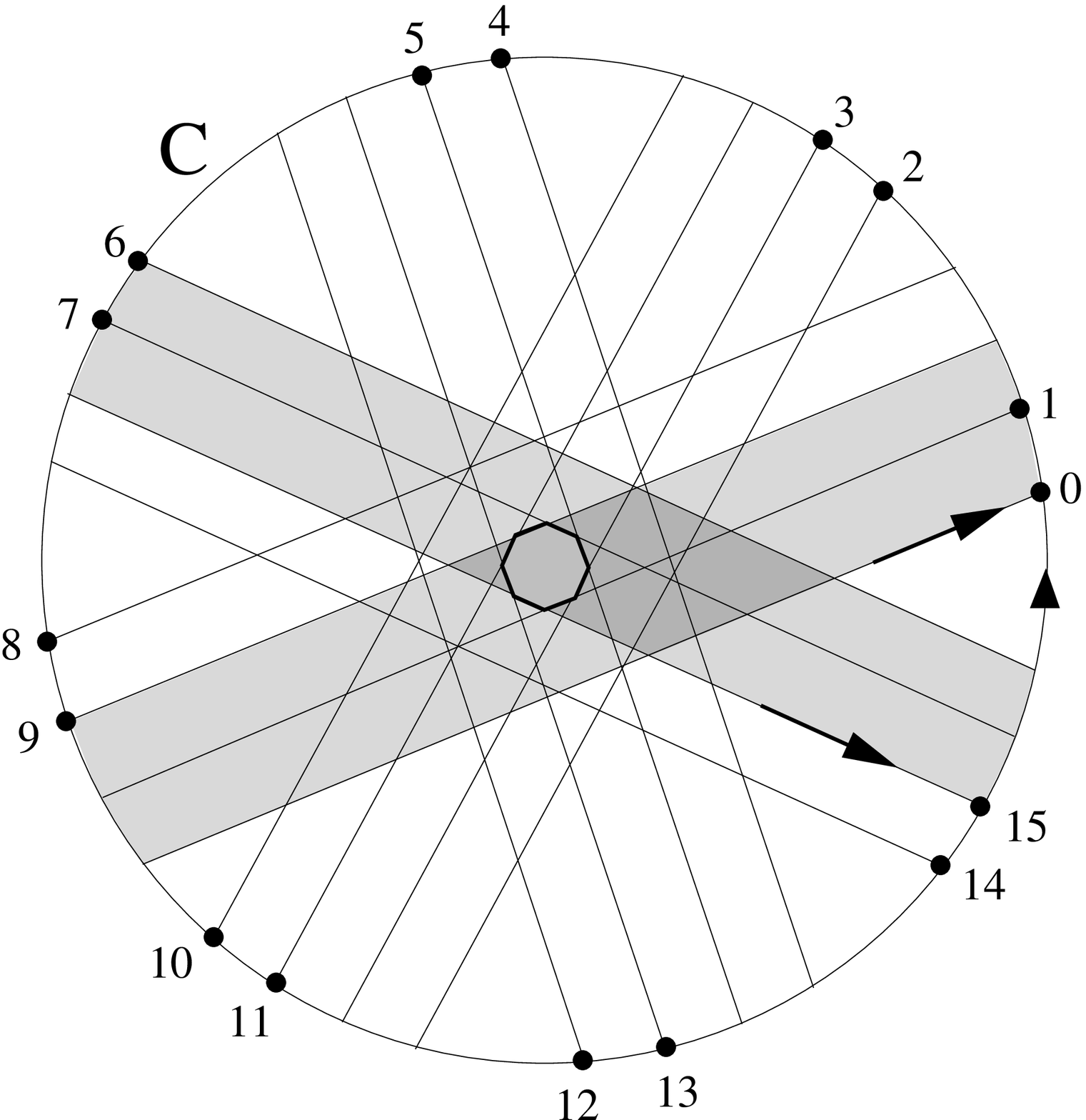}}
\newline
{\bf Figure 2.3:\/} The pointed strips associated to a regular octagon
\end{center}  

Thus, we have associated to $P$ a total of $2n$ triples
$(\Sigma_k,L_k,V_k)$.   Here $(\Sigma_k,L_k)$ is a pointed
pinwheel strip and $V_k$ is the associated vector.
We cyclically order these $2n$ triples in the following 
way.  We choose a large circle $C$ centered at the origin,
and orient $C$ counterclockwise.   Each pair
$(\Sigma_k,L_k)$ determines a point $c_k \in C$,
namely the intersection point of $L_k \cap P$ towards
which $L_k$ points.   We order our
triples so that the points $c_1,...,c_{2n}$
occur in counterclockwise order along $C$.
Figure 2.3 shows these points, in case $P$ is
the regular octagon.
In Figure 2.3 we highlight the pointed
strips $\Sigma_{15}$ and $\Sigma_0=\Sigma_{16}$.   

Let $f_k$ be the strip map associated to the
pair $(\Sigma_k,V_k)$.   We say that the {\it pinwheel map\/} is the
composition
\begin{equation}
\Phi: f_{2n} \circ \ldots \circ f_0: \Sigma_0 \to \Sigma_0
\end{equation}

We say  that $p \in \Sigma_0$ is {\it preferred\/} if,
for all $x \in P$, the vector $p-x$ has positive
dot product with vectors parallel to $L_0$.   
We call the subset of preferred points the
{\it preferred half\/} of $\Sigma_0$ and denote it by $H_0$.
We usually restrict $\Phi$ to $H_0$.   Note that $\Phi$ need
not carry $H_0$ to $H_0$.   However, if $p \in H_0$ is sufficiently far from
$\Sigma_0-H_0$ then $\Phi(p) \in H_0$.  
\newline
\newline
{\bf Remark:\/}
 For our purposes here,
the boundary of $H_0$ (that separates $H_0$ from $\Sigma_0-H_0)$
is not important.   What is important is the {\it end\/} of
$\Sigma_0$ that $H_0$ determines.   Below, and just for
the case of the regular octagon, we will replace
$H_0$ by a similar set, namely $\Sigma_0^1$ below, that
is invariant under $\Phi$.
\newline

Let $\phi$ denote the outer billiards map and let
$\phi^2$ be the square of the outer billiards map.
Let $O_+(\phi^2,p)$ denote the forwards $\phi^2$-orbit of $p$.

\begin{lemma}
\label{far2}
The following is true for all points $p \in H_0$
that lie sufficiently far from $P$.
Let $q$ be the first point in $O_+(\phi^2,p)$ that lies in $H_0$.  Then
$q=\Phi(p)$.
\end{lemma}

\startproof
For each point $p \in \R^2-P$ for which $\phi^2$ is well-defined,
there is a vector $V(p)$ such that $\phi^2(p)-p=V(p)$.  The function
$p \to V(p)$ is locally constant.
Let $\Delta$ be the disk bounded by the large circle $C$, as in Figure 2.2.
Let $W_k$ be the open acute wedge shaped region of $\R^2-\Delta$ bounded by the
lines $L_{k-1}$ and $L_k$.    The vector $V_k$ has the property that
$V(p)=V_k$ for at least some points $p \in W_k$.
But one can check easily that $\phi^2$ is defined on all of $W_k$
provided that the circle $C$ is taken large enough. Hence $V(p)=V_k$
for all $p \in W_k$.

Starting at $p=p_0 \in W_1 \cap H_0$, the successive points in $O_+(\phi^2,p_0)$
have the form $p_0+mV_1$ for $m=1,2,3...$.  This
continues until we reach a point $p_1=p_0+m_1V_1 \in \Sigma_1$.  But then
$p_1=f_1(p_0)$.    
Starting at $p_1\in W_2 \cap \Sigma_1$, the successive points in $O_+(\phi^2,p_1)$
have the form $p_1+mV_2$ for $m=1,2,3...$.  This
continues until we reach a point $p_2=p_1+m_2V_2 \in \Sigma_2$.  But then
$p_2=f_2(p_1)$.     Continuing in this way, we get $p_{n} \in \Sigma_0$.  But
$p_n$ and $p_0$ lie in opposite components of $\Sigma-\Delta$.  We have gone
halfway around.    Continuing the process, we finally arrive at
$p_{2n}=q=\Phi(p_0) \in H_0$, as claimed.
\endproof

\subsection{The Arithmetic Graph}
\label{ag0}

Let $p=p_0 \in \Sigma_0$ be a point.   As in
the proof of Lemma \ref{far2}, 
there are integers $m_1,...,m_{2n}$ such that
\begin{equation}
\label{gr2}
p_k:=f_k(p_{k-1})=p_{k-1}+m_k V_k \in \Sigma_k.
\end{equation}

Let $v_1,...,v_n$ be the vertices of $P$, ordered
clockwise around $P$.
Let $e_1,...,e_n$ be the standard basis of $\R^n$.
For each vector $V=V_k$ there are indices
$i$ and $j$ such that
$V_k=2(v_i-v_j)$.   Of course,
$i$ and $j$ depend on $k$, but we are suppressing this
from our notation.  We define
\begin{equation}
\label{lift0}
\widetilde V=2(e_i-e_j); \hskip 30 pt k=1,...,2n.
\end{equation}

We define
\begin{equation}
\label{spectrum}
\widetilde \Phi(p)=
\sum_{k=1}^{2n} m_k \widetilde V_k; \hskip 30 pt
\widetilde \Phi^{k+1}(p)=\widetilde \Phi(\widetilde \Phi^{k}(p)).
\end{equation}
We define $\Gamma(P,p)$ to be the polygonal path in $\R^n$ that
starts at $0$ and has consecutive vertices
$\widetilde \Phi(p)$, $\widetilde \Phi^2(p)$, etc.
\newline
\newline
{\bf Remarks:\/} \newline
(i) It might seem at first that $\widetilde \Phi$ is somehow a ``lift''
of $\Phi$ to $\R^n$.   This is only partially true.
What is true is that 
$\widetilde \Phi^k$ is a map from $\Sigma_0$ to $\R^n$.
This map is only defined on points where $\Phi^k$ is well defined.
\newline
(ii)
By construction the arithmetic graph is a {\it lattice path\/}
in that its vertices lie in $\Z^n$.    As remarked in the
introduction, one should compare the lattice paths of
[{\bf V\/}]. 
\newline
(iii)
There is a constant $C'$, depending only on $P$, 
such that $|m_{k+n}-m_k|<C'$.   The reason is that
the successive points $p_1,...,p_{2n}$ are the
vertices of a convex polygon that is centrally symmetric
to within a bounded error.  Hence,
the edges of $\Gamma(P,p)$ have length at most $C$,
for some constant $C$ depending only on $P$.
This property is important when we take
rescaled limits of the graph. 
\newline
(iv) 
There is an obvious projection $\pi: \R^n \to \R^2$ given
by $\pi(e_k)=v_k$.    By construction,
$p+\pi(\Gamma(P,p))$, meaning the translation
of $\pi(\Gamma(P,p))$ by $p$, is
 exactly the forwards $\Phi$ orbit of $p$.
What makes $\Gamma(P,p)$ interesting is that sometimes
other projections are much more revealing.   We will illustrate
this in the next chapter with many examples.

\subsection{Formulas for the Octagon}
\label{octo}

Since we are concentrating on the regular octagon,
it seems worthwhile giving explicit formulas
in this case.    Recall that the vertices of
the regular octagon are $\omega^k$ for
$k=0,...,7$.   We let
$(ij)$ denote the line through $\omega^i$ and $\omega^k$,
oriented from $i$ to $j$.   We let $k(ij)$ denote the
image of $(ij)$ after we apply reflection in $\omega^k$.
Finally, we let $[ij]$ denote the vector that points from
$\omega^i$ to $2\omega^j-\omega_i$.   
The vector $-[ij]$ is literally the negative of $[ij]$. 

\begin{center}
\resizebox{!}{3.8in}{\includegraphics{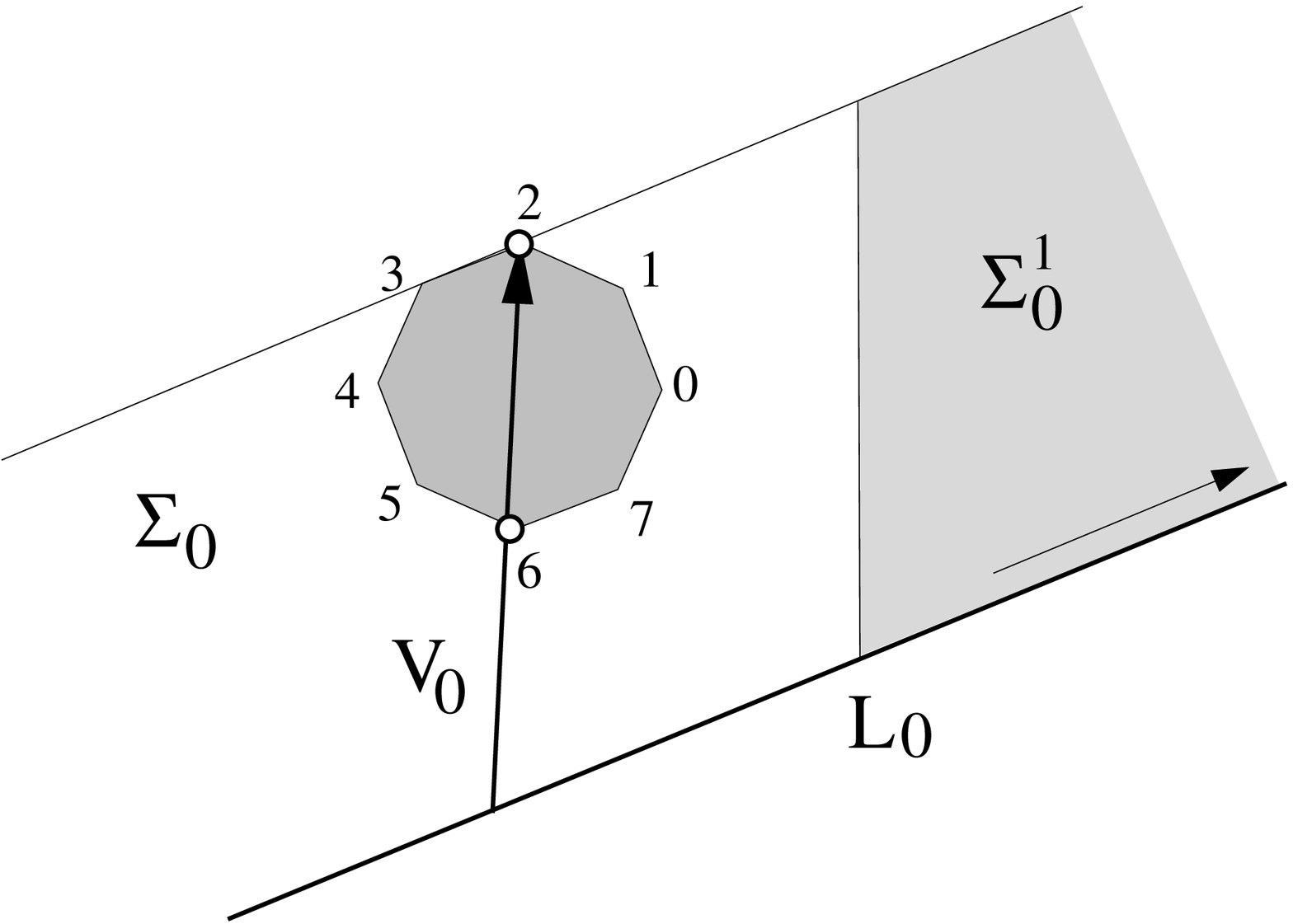}}
\newline
{\bf Figure 2.4:\/} The triple $(\Sigma_0,L_0,V_0)$.
\end{center}  

\noindent
{\bf Definition:\/}
The shaded region in Figure 2.4 is the subset of $\Sigma_0$
to the right of the line $y=2+\sqrt 2$.   We let $\Sigma_0^1$
denote this set.
\newline

To illustrate our notation by way of example,
The triple $(\Sigma_0,L_0,V_0)$ is specified by
\begin{equation}
{\bf 0:\/} \hskip 50 pt 6(23) \hskip 30 pt (23) \hskip 30 pt -[26].
\end{equation}
We first list $L_0$, then the other component of $\Sigma_0$,
and then $V_0$.   See Figure 2.4.

The listing for $(\Sigma_2,L_2,V_2)$ is obtained from
the one for $(\Sigma_0,L_0,V_0)$ simply by incrementing
all the indices by $1$.   That is
$$
{\bf 2:\/} \hskip 50 pt 7(34) \hskip 30 pt (34) \hskip 30 pt -[37].
$$
And so on.

\begin{center}
\resizebox{!}{3.8in}{\includegraphics{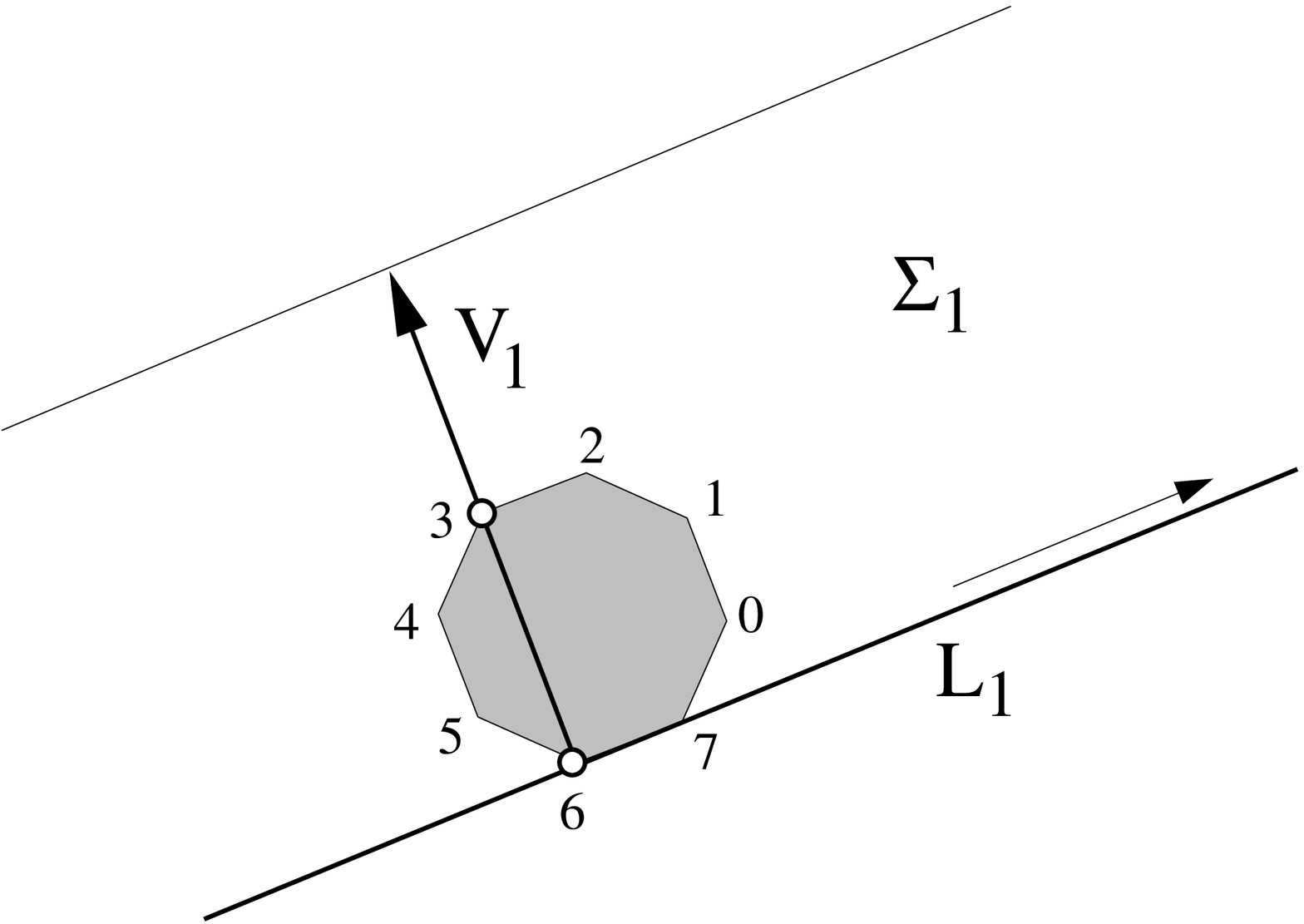}}
\newline
{\bf Figure 2.5:\/} The triple $(\Sigma_1,L_1,V_1)$.
\end{center}  

The triple $(\Sigma_1,L_1,V_1)$ is specified by
\begin{equation}
{\bf 1:\/} \hskip 50 pt  (67) \hskip 30 pt 3(67) \hskip 30 pt [63]
\end{equation}
See Figure 3.4. 

The listing for $(\Sigma_3,L_3,V_3)$ is obtained from
the one for $(\Sigma_1,L_1,V_1)$ simply by incrementing
all the indices by $1$.   That is
$$
{\bf 3:\/} \hskip 50 pt  (70) \hskip 30 pt 4(70) \hskip 30 pt [74]
$$
And so on.

\subsection{Discussion}
\label{discuss}

Lemma \ref{far2} gives a precise relationship between the
pinwheel map and the outer billiards map for points
in $H_0$ that are far from $P$.    The relationship
between $\phi^2$ and $\Phi$ for points near $P$ is
rather subtle, and it was the purpose of the paper
[{\bf S3\/}] to explain the relationship.   
Here we summarize some of the results from
[{\bf S3\/}].  These results are not needed for
this paper, but they are worth knowing.

Let $O_+(\Phi,p)$ denote the forward $\Phi$-orbit of $p \in H_0$.
Here is a consequence of our work in [{\bf S3\/}]

\begin{theorem}
\label{unbound}
Suppose that all our constructions are made
with respect to a convex polygon with no
parallel sides.  Then, 
the following holds for any sufficiently large
disk $\Delta$ centered at the origin.
Let $p \in \Sigma_0-\Delta$.    Then 
$O_+(\phi^2,p)$ returns to $H_0-\Delta$
if and only if $O_+(\Phi,p)$ returns to
$H_0-\Delta$, and the two points of return
are the same.
\end{theorem}

There are two subtle points to Theorem \ref{unbound}.
First, the forward orbits of $p$, in either case, might wind
many times around $P$ before returning to $H_0-\Delta_0$.  
Second, Lemma \ref{far2} only applies to points
of $H_0$ that are sufficiently far from $P$.   After
our orbits exit $H_0-\Delta$, then might wander quite
close to $P$ and only ``come back out'' much later on.
However, magically, they ``come back out''
in exactly the same way.   
One consequence of Theorem \ref{unbound} is 

\begin{theorem}
\label{unbound0}
Suppose that all our constructions are made
with respect to a convex polygon with no
parallel sides.
The pinwheel map $\Phi$ has unbounded orbits relative to
$P$ if and only if outer billiards has unbounded orbits
relative to $P$.
\end{theorem}

We only proved these results for polygons without
parallel sides because we wanted to avoid certain
technical complications.  We fully expect that the same
results hold for all convex polygons.
However, we have not yet worked out the details.
For the present paper, which deals (rigorously) just with the
regular octagon, we will prove something stronger than
Theorem \ref{unbound}.  Namely, in \S \ref{qr} we prove

\begin{lemma}[Invariance]
For $k=1,2,3...$ let
$\Sigma_0^k \subset \Sigma_0$ be the subset
consisting of points lying to
the right of the line $x=k(1+\sqrt 2)$ .
Then $\Sigma_0^k$.
$\Phi$-invariant for all $k \geq 1$ 
\end{lemma}

The argument in Lemma \ref{far2} works directly
for all orbits in $\Sigma_0^1$.   Accordindly,
for points in $\Sigma_0^1$, we can be sure that
the pinwheel dynamics and the first return
dynamics of the square outer billiards map
coincide.   

The Invariance Lemma implies that
all forward orbits of both $\phi$ and $\Phi$ are bounded.
By symmetry, the same goes for the backwards orbits.
Hence, the Invariance Lemma implies Theorem
\ref{unbound0} for the regular octagon.
Compare \S \ref{id}.

The Invariance Lemma is a special case of a general result
concerning quasi-rational polygons.
The polygon $P$ is called {\it quasi-rational\/} if it may be scaled
so that
$${\rm area\/}(\Sigma_k \cap \Sigma_{k+1}) \in \Z \cup \infty$$
for all $k=1,...,2n$.   Here indices are taken mod $2n$, as
usual.   In case no sides of $P$ are parallel, the above areas
are all finite.  For regular polygons, one can scale so that
all the finite areas are $1$.   Hence, regular polygons
are quasi-rational.   Likewise, polygons with rational
vertices are quasi-rational.

 As we mentioned in the introduction,  
it is proved in [{\bf VS\/}], [{\bf K\/}], and [{\bf GS\/}]
that all outer billiards orbits are bounded for a quasi-rational
polygon.  
  In [{\bf S3\/}] we give a self-contained proof
of this result, in case $P$ has no parallel sides.

\newpage

\section{Examples of Arithmetic Graphs}
\label{graph}

\subsection{Rational Kites}
\label{penrose}

As in [{\bf S2\/}] we let $K(A)$ be the kite with vertices
\begin{equation}
(-1,0); \hskip 30 pt
(0,1); \hskip 30 pt
(0,-1); \hskip 30 pt
(A,0); \hskip 30 pt A \in (0,1).
\end{equation}
The case $A=1$ corresponds to the square, a trivial case we
ignore.   $K(A)$ is called (ir)rational iff $A$ is (ir)rational.

When $A=p/q$, we call the point $(1/q,1)$ the {\it fundamental
point\/},  for reasons we explain in great detail in [{\bf S2\/}]
and briefly as follows.
The orbit of any point $(x,1)$ has the same combinatorics
as the point $(1/q,1)$ for any $0<x<2/q$.  For this
reason, the point $(1/q,1)$ is a representative of the
dynamics of points ``arbitrarily close'' to the top
vertex of $K(p/q)$ and on the same horizontal line
as the vertex.

We will consider the arithmetic graphs
\begin{equation}
\Gamma_{p/q}=\Gamma\big(K(p/q),(1/q,1)\big)
\end{equation}
for various choices of $p/q$.
We will draw a certain projection of this graph into the plane.
We will not specify the projection we use explicitly because
it is tedious to do so.  The fact is that 
$\Gamma_{p/q}$ lies in a thin neighborhood of a
$2$-plane in $\R^4$ and a random projection
will produce pictures very similar to ours.  Indeed,
any projection looks about the same, up to an
affine transformation of $\R^2$.

Consider the rational sequence $\{a_n\}$ which starts
$0,1$ and obeys the rule  $a_{n}=4a_{n-1}+a_{n-2}$.
This sequence starts out
$$0,1,4,17,72,305, \ldots$$
We have
$$
\lim_{n \to \infty} a_{n-1}/a_n=\sqrt 5-2=\phi^{-3}.
$$
Here $\phi$ is the golden ratio.  The kite
$K(\phi^{-3})$ is affinely equivalent to the 
Penrose kite.    The first four
quotients are
$$1/4 \hskip 30 pt 4/17 \hskip 30 pt 17/72 \hskip 30 pt
72/305 \ldots$$

Figure 3.1 shows $\Gamma_A$
for the first four of these quotients.
\begin{center}
\resizebox{!}{2.5in}{\includegraphics{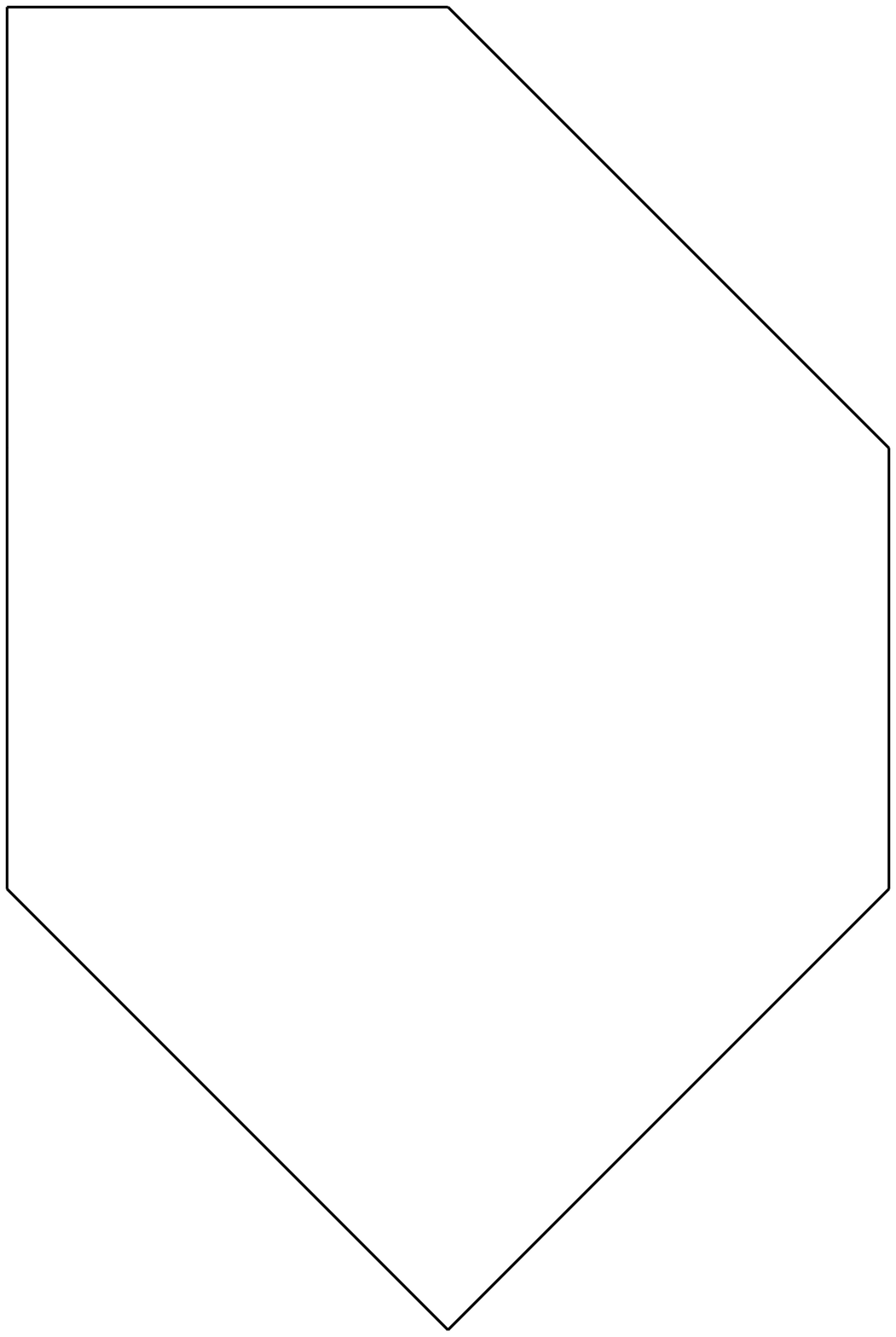}}
\resizebox{!}{2.5in}{\includegraphics{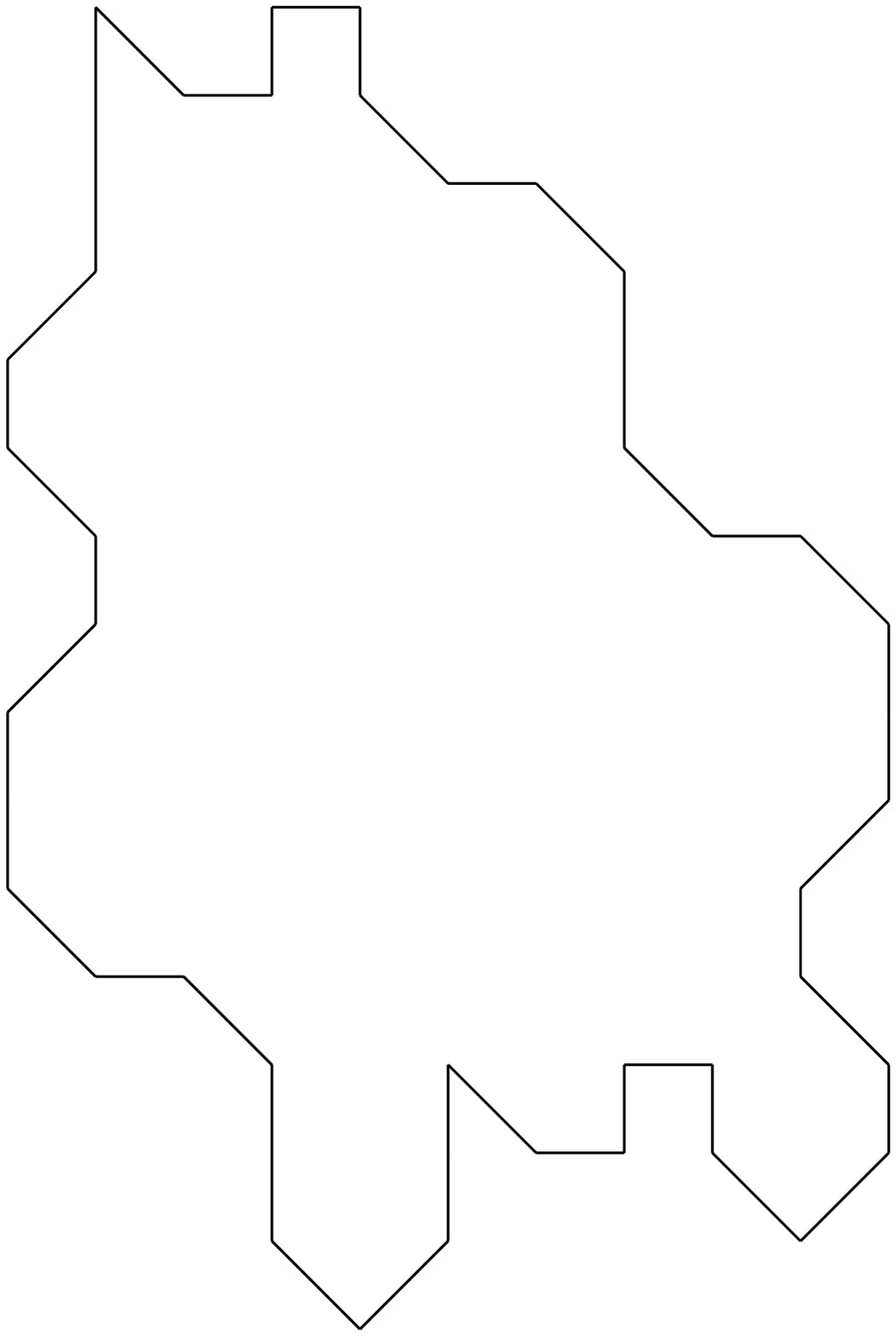}}
\newline
\resizebox{!}{2.5in}{\includegraphics{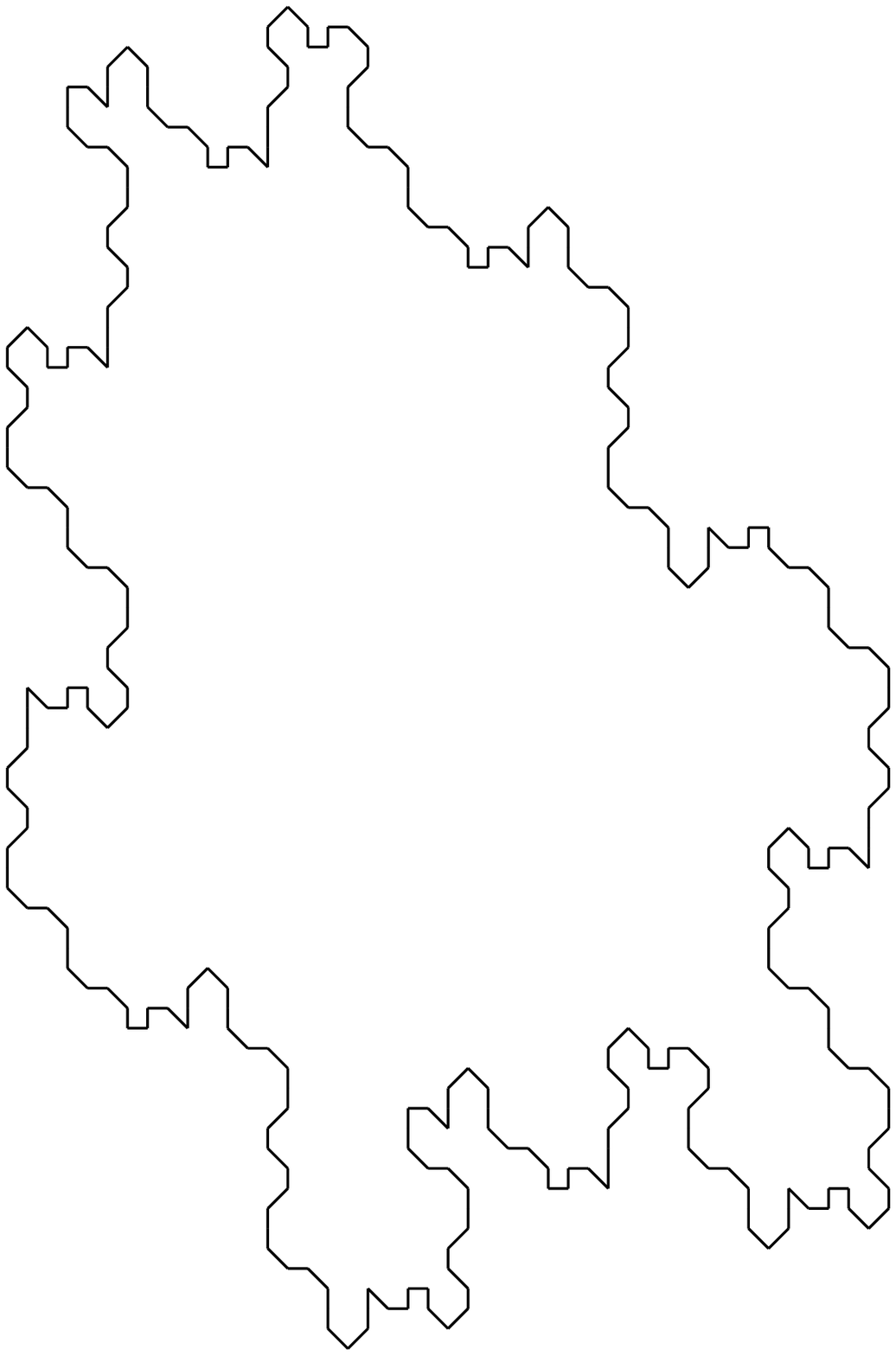}}
\resizebox{!}{2.5in}{\includegraphics{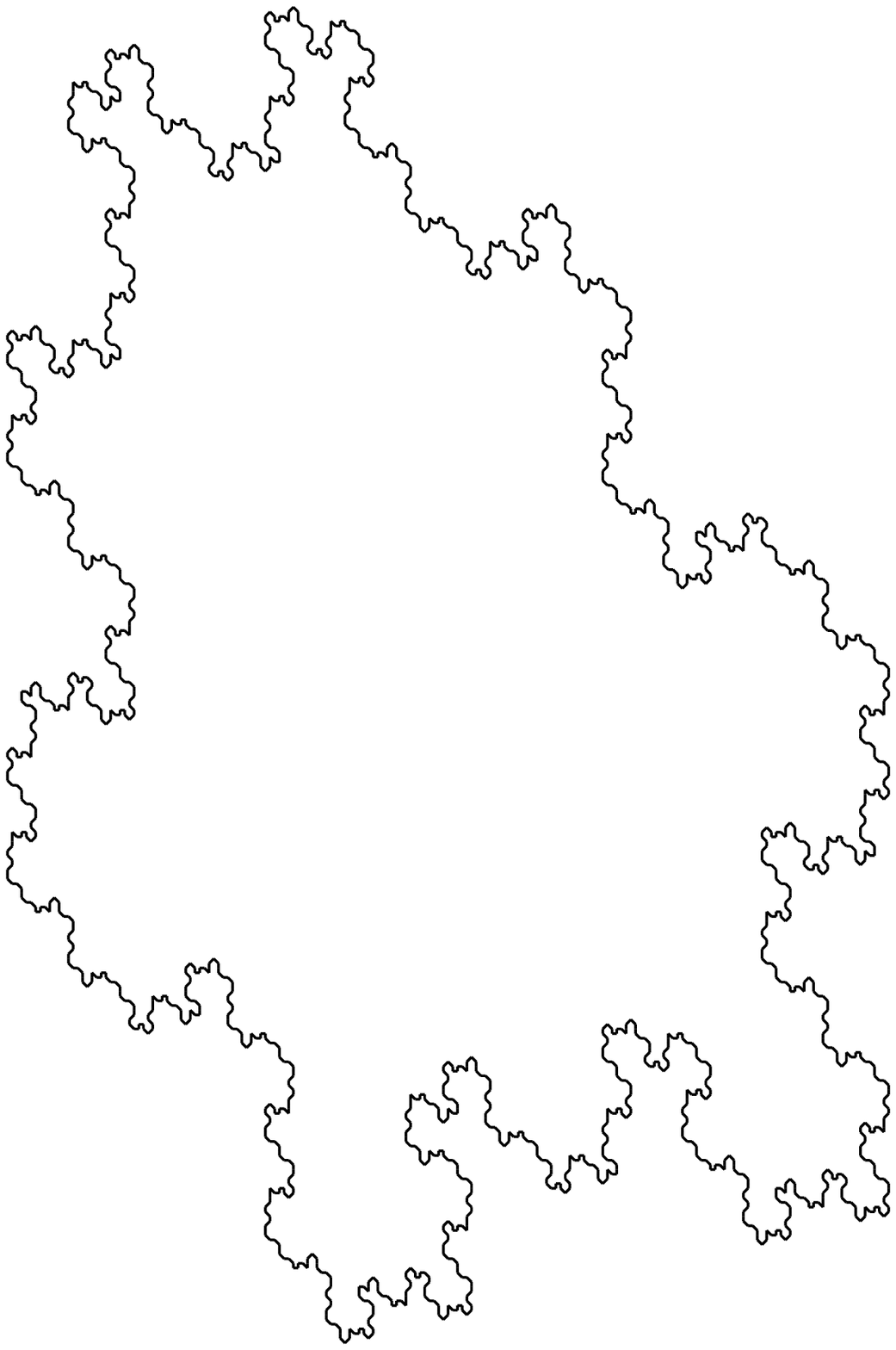}}
\newline
{\bf Figure 3.1:\/} $\Gamma_A$ for $A=1/4, 4/17, 17/72, 72/305$.
\end{center}  

We want to emphasize that our projection is such that the
vertices of these polygons all lie in $\Z^2$.   Were the
polygons drawn to scale, each one would be about
$\phi^3 \approx 4.23$ times as large as the
predecessor.  However, when we rescale them so that
they are about the same size, a fractal structure emerges.
Later in the paper, we will formalize what we mean to
taking the {\it rescaled limit\/} of a sequence of
arithmetic graphs, but we hope that the above informal
discussion makes the general idea clear.

Another nice sequence of pictures is given by the
quotients of the rational sequence $\{a_n\}$ where
$a_0=1$ and $a_1=2$ and $a_n=2a_{n-1}+a_{n-2}$.
In this case $\lim a_{n-1}/a_n=\sqrt 2-1$.  One
of the quotients is $169/408$.   Figure 3.2 shows
the picture of the corresponding arithmetic graph.
One can see the 
fractal structure emerging just from this single polygon.

\begin{center}
\resizebox{!}{4.5in}{\includegraphics{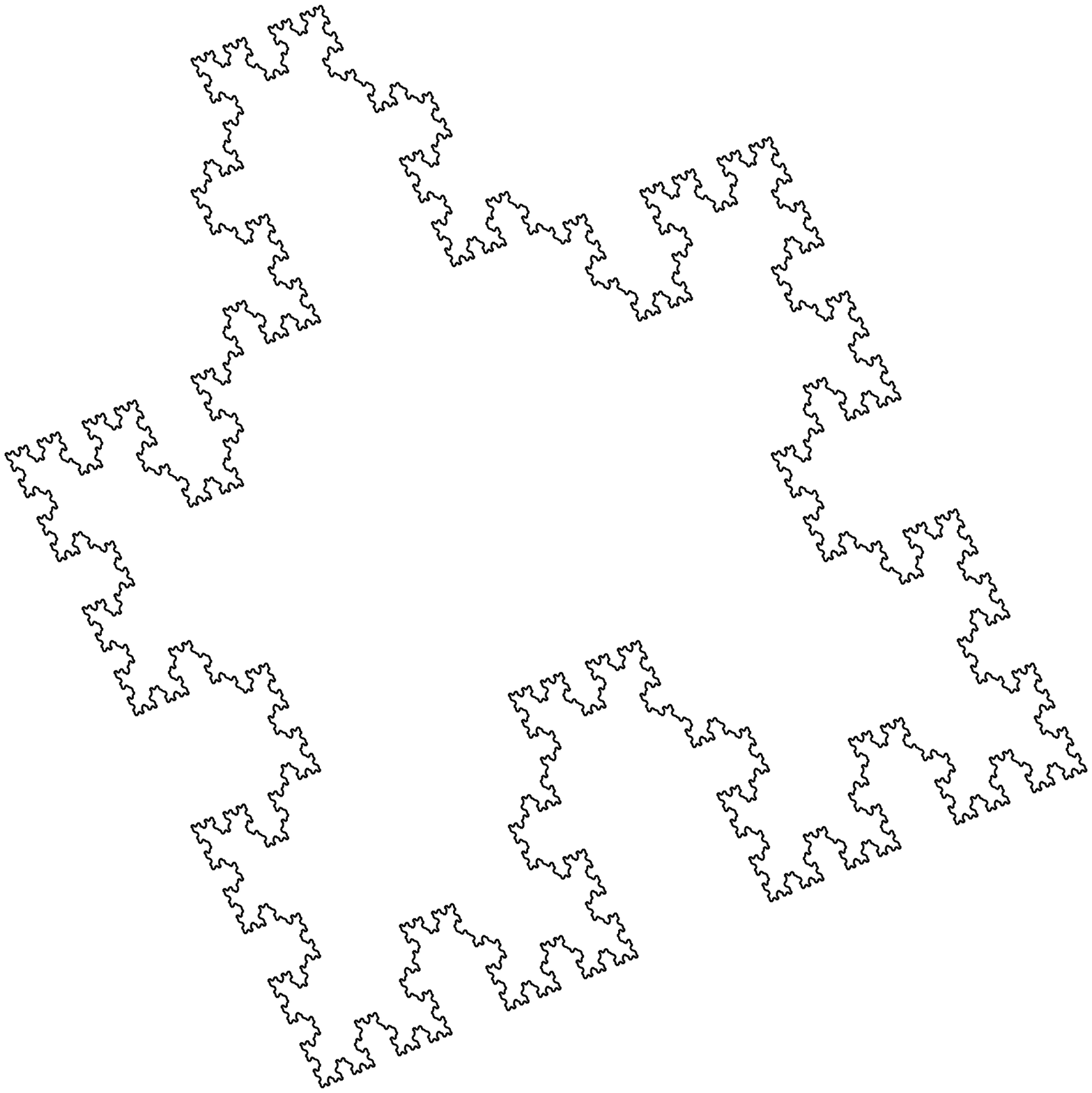}}
\newline
{\bf Figure 3.2:\/}  $\Gamma_A$ for $A=169/408$.
\end{center}  

There is one interesting feature of $K_{p/q}$ that
one sees almost immediately.  If $pq$ is even, then
$K_{p/q}$ is a closed embedded loop.  If $pq$ is odd, then
$K_{p/q}$ is an open embedded polygonal curve.   The difference
derives from the fact that the orbit of $(1/q,1)$ is
stable when $pq$ is even and unstable when $pq$ is
odd. By this we mean that a small change in the
kite parameter destroys the orbit of $(1/q,1)$ when
$pq$ is odd, but not when $pq$ is even.  We prove
the above statements in [{\bf S2\/}].

My program Billiard King \footnote{You can download this program
from my website http://www.math.brown.edu/$\sim$res}
allows the reader to draw these pictures 
for any (smallish) rational parameter.
Our monograph [{\bf S2\/}] discusses these graphs in
great detail.  In fact, one could consider [{\bf S2\/}]
as an exploration of the arithmetic graphs associated
to kites.

\subsection{The Regular Pentagon}
\label{penta}

Let $P$ be the regular pentagon.  We scale $P$ so that its vertices
have the form $\omega^k$ for $k=0,1,2,3,4$.  Here
$\omega=\exp(2 \pi i/5)$ is the usual primitive
$5$th root of unity.     In this case, we have
maps
$\pi_k: \R^5 \to \C$ given by
\begin{equation}
\pi_k(X)=\sum_{j=0}^4 x_j \omega^{kj}.
\end{equation}
Here $X=(x_0,x_1,x_2,x_3,x_4)$.
The map $\pi_1$ is just the projection mentioned above.
Note that $\pi_k$ and $\pi_{5-k}$ agree up to complex
conjugation.  Thus, the only remaining interesting projection
is $\pi_2$.

Relative to any convex polygon, any periodic point $x$ lies in
a convex polygon $P_x$ consisting of points that have exactly
the same dynamics.  The outer billiards map moves the tile
$P_x$ around isometrically.   The pictures we show of
outer billiards will draw the orbits of these tiles, rather than
the orbits of individual points, because the pictures are
more revealing.   Figure 3.3 shows the picture for the
regular pentagon.

\begin{center}
\resizebox{!}{3.5in}{\includegraphics{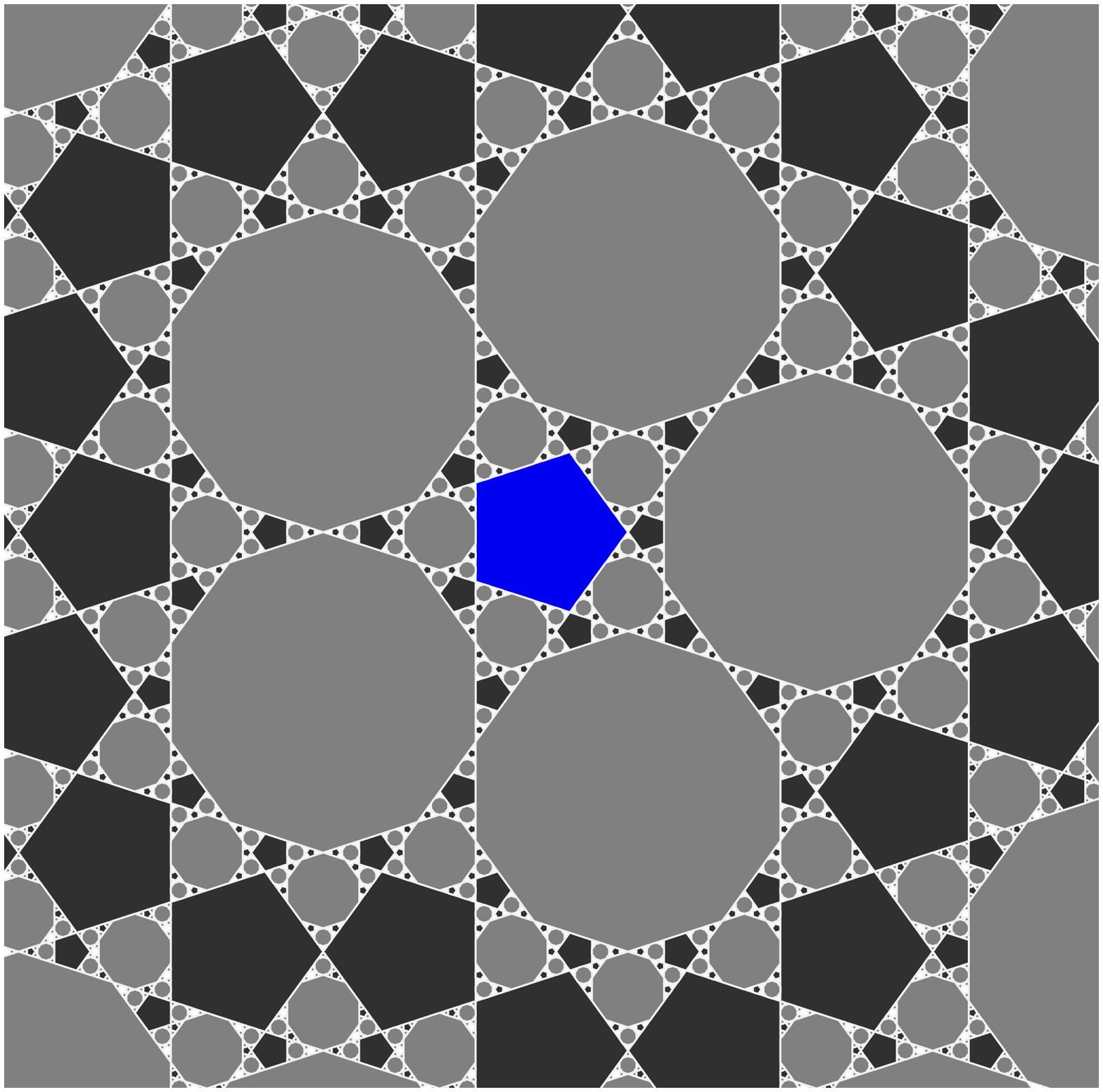}}
\newline
{\bf Figure 3.3:\/}  Some of the periodic orbits for the regular pentagon
\end{center}  

In the case of the regular pentagon, there are two kinds of
orbits, those lying in regular pentagon tiles and those lying
in regular $10$-gon tiles.  This is worked out in
[{\bf T2\/}], and also discussed in [{\bf BC\/}].
The smaller the tile, the longer the periodic orbit.

The left half of
Figure 3.4 shows a plot of $\pi_2(\Gamma(P,p))$, where $p$ is a point
of $\Sigma_0$ that lies inside a pentagonal tile.  We choose
the period to be fairly long, so that the actual size of the
graph is quite large.  One can see the emerging fractal
structure in the rescaled picture.  The right half of
Figure 3.4 shows the  corresponding picture for a periodic
point that lies inside a $10$-gon tile.  The reader can
probably see that these tiles somehow fit together
``hand-in-glove''.

\begin{center}
\resizebox{!}{2in}{\includegraphics{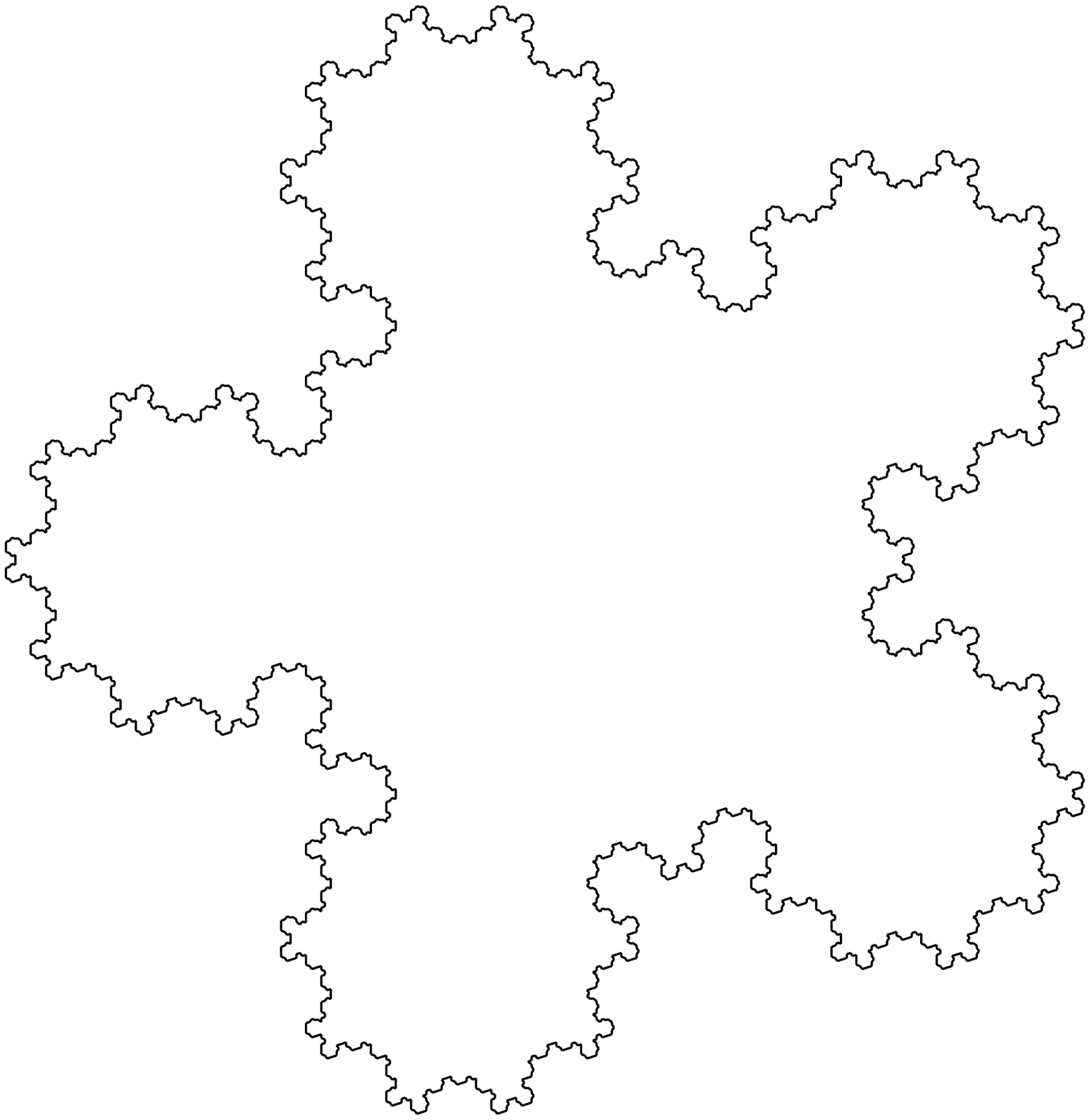}}
\resizebox{!}{3in}{\includegraphics{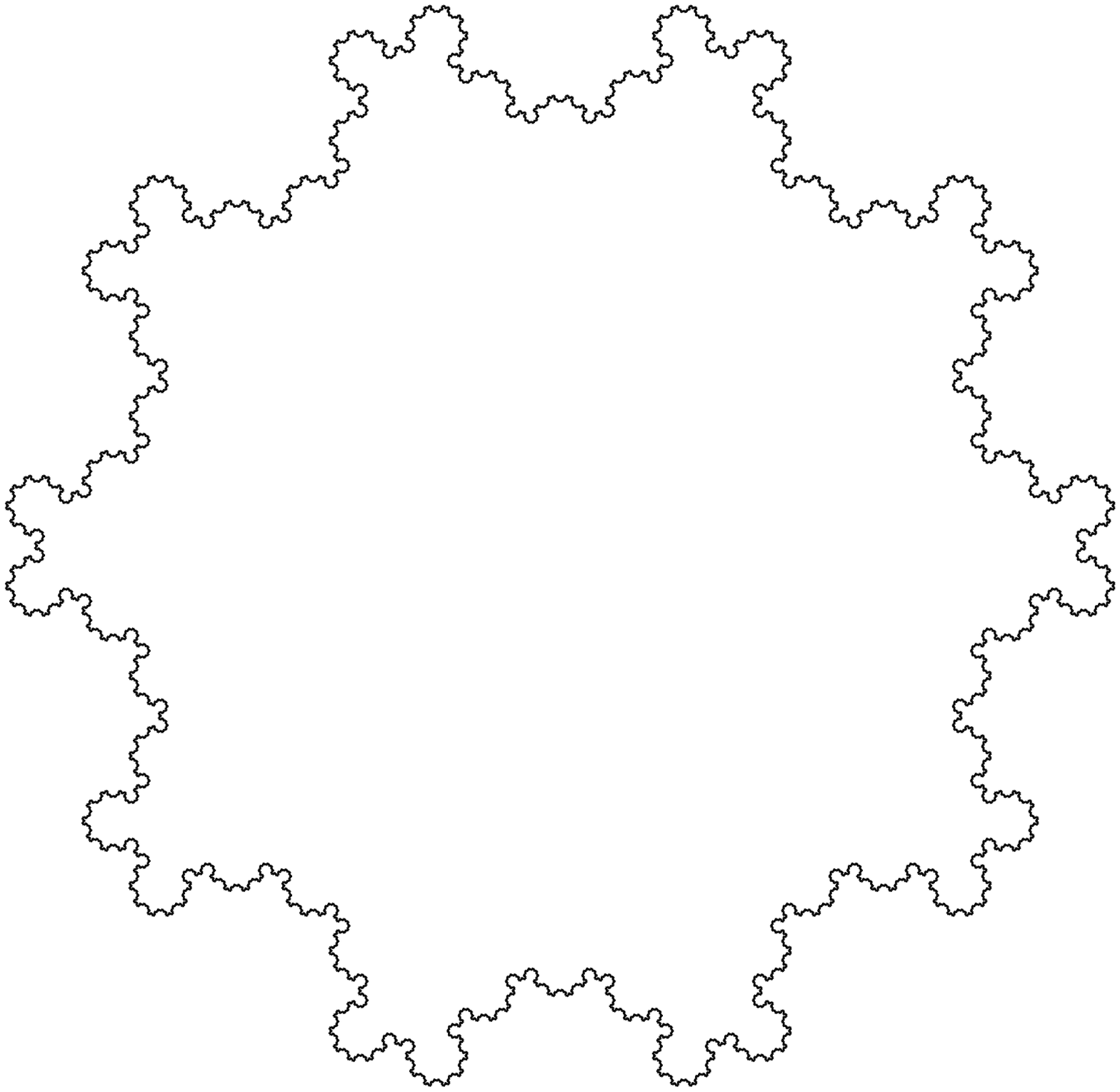}}
\newline
{\bf Figure 3.4:\/}  Some of the periodic orbits for the regular pentagon
\end{center}  

Any choice of long periodic orbit will yield a picture that is essentially
identical to one of the two above.  
In the case of the regular pentagon, the arithmetic graphs lie within a thin
tubular neighborhood of a $2$-plane in $\R^5$.  Any other projection
of the graph would look like an affine image of one of the two pictures
shown above.  This is why we say that, in the
regular pentagon case, there are two fractal curves that somehow
control the structure of the pinwheel dynamics.

As we remarked in the introduction, we do not offer proofs of these
results, because we will give complete proofs of the analogous
results for the regular octagon.

\subsection{The Regular Octagon}
\label{octopix}

Let $P$ be the regular octagon.  We scale $P$ as in
\S \ref{octo}, and we define the projection maps
as in \S \ref{penta}.     In the case, there are
two interesting projections, $\pi_2$ and $\pi_3$.
(It turns out that $\pi_4$ always maps the graph
into a line segment.)

As we discuss in the introduction, in connection with Figure 1.3,
we call a periodic orbit {\it odd\/} if it intersects $\Sigma_0^1$,
the half-strip from Figure 2.4., in
$3^k$ points, for $k$ odd.    It turnso out that ``half'' the
periodic orbits that start in $\Sigma_0^1$ have this property.
The odd periodic orbits are the lightly colored ones in Figure 3.5
that happen to intersect $\Sigma_0^1$.    This includes all
the lightly colored octagons that lie outside the first large
ring of $8$ octagons.

\begin{center}
\resizebox{!}{5in}{\includegraphics{Pix/octagon0.ps}}
\newline
{\bf Figure 3.5:\/} The tiling associated to the octagon
\end{center}  

Figure 3.6 shows some $\pi_3$ projections of the
arithmetic graphs associated to odd orbits.
For $\pi_3$, the projections corresponding to the
even orbits look similar.
Note that Figure 1.2, which was drawn using the
subdivision method discussed in \S \ref{snowflake},
emerges as we rescale the pictures.

\begin{center}
\resizebox{!}{2.6in}{\includegraphics{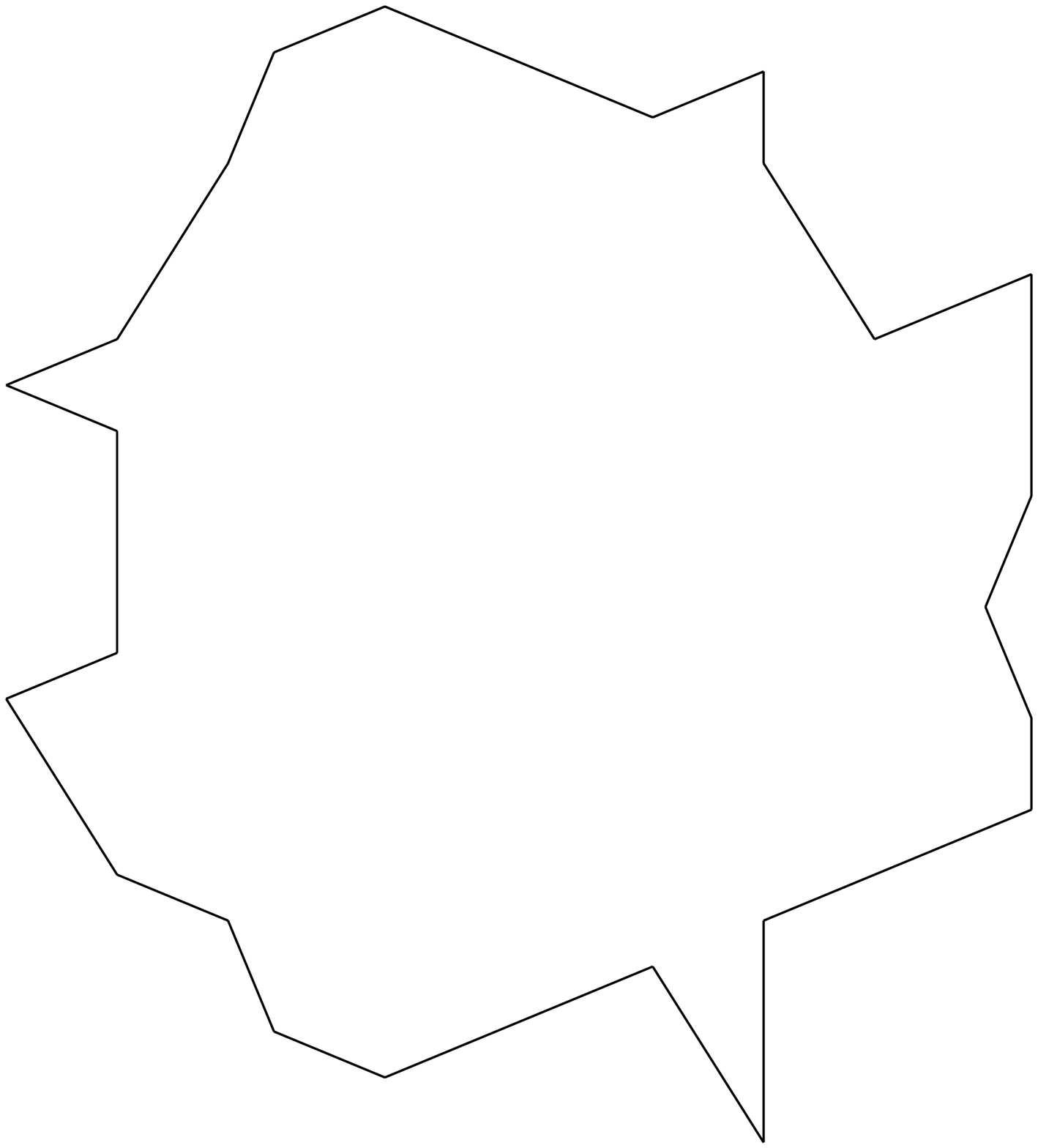}}
\resizebox{!}{2.6in}{\includegraphics{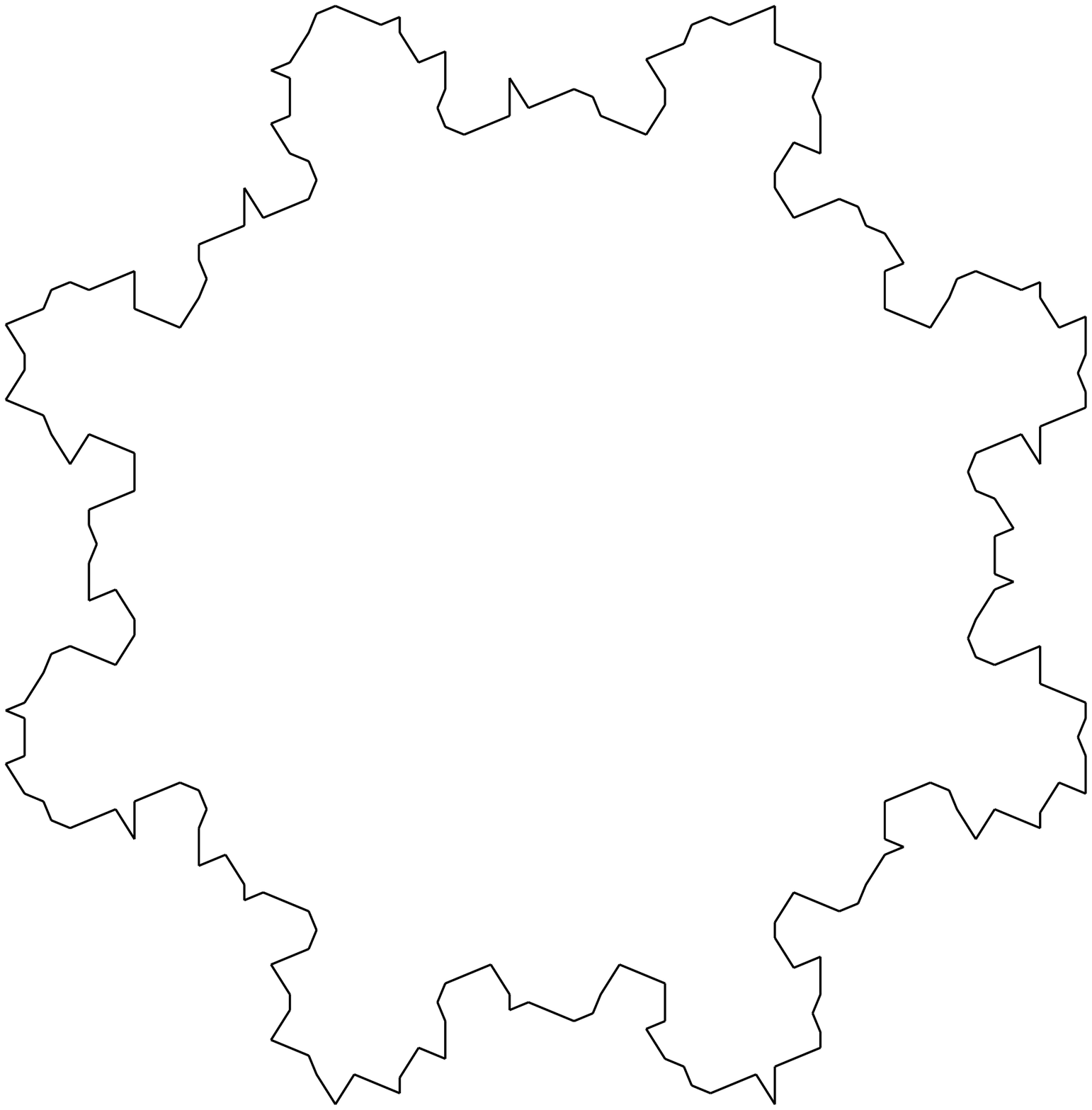}}
\newline
\resizebox{!}{2.6in}{\includegraphics{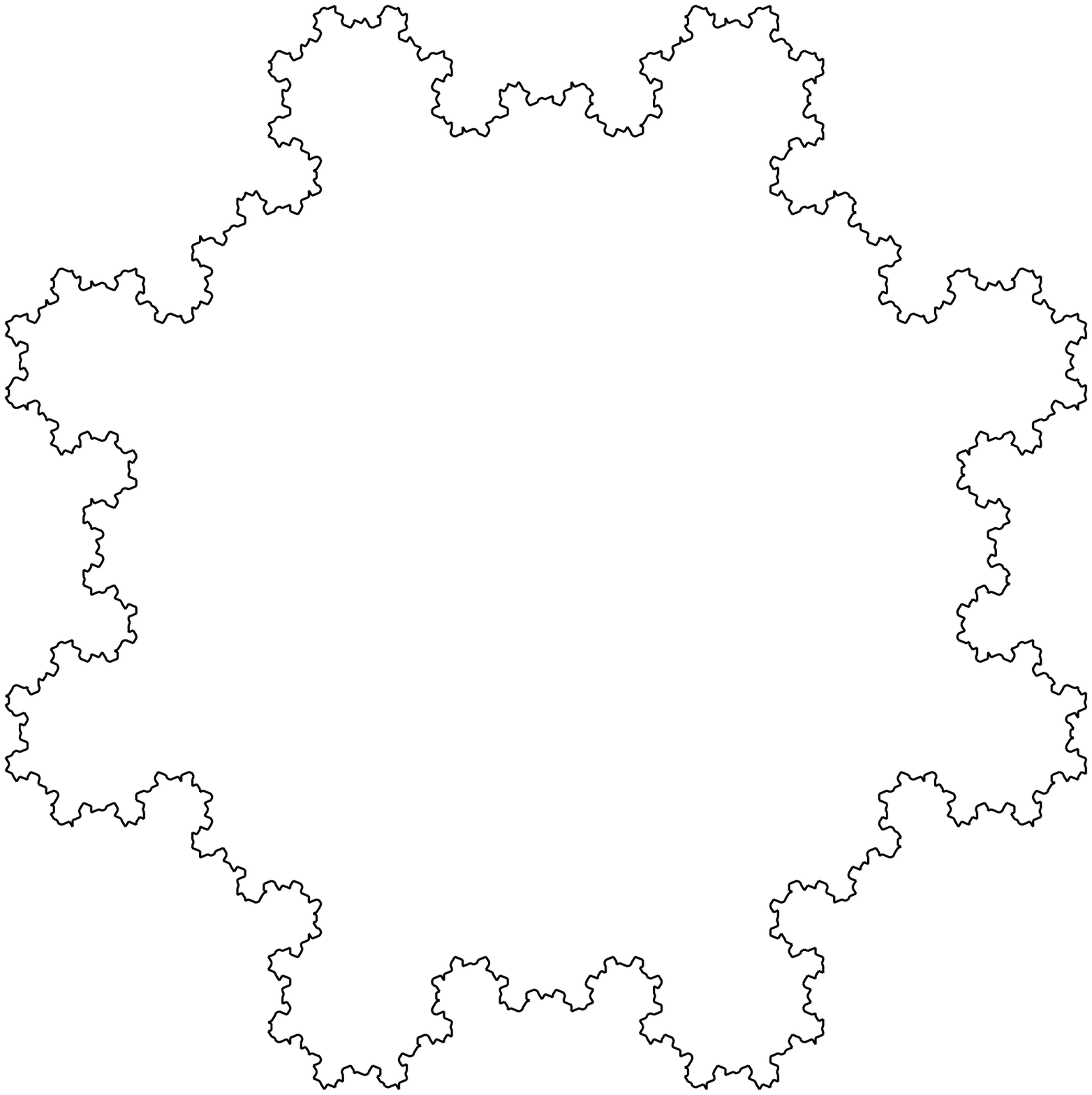}}
\resizebox{!}{2.6in}{\includegraphics{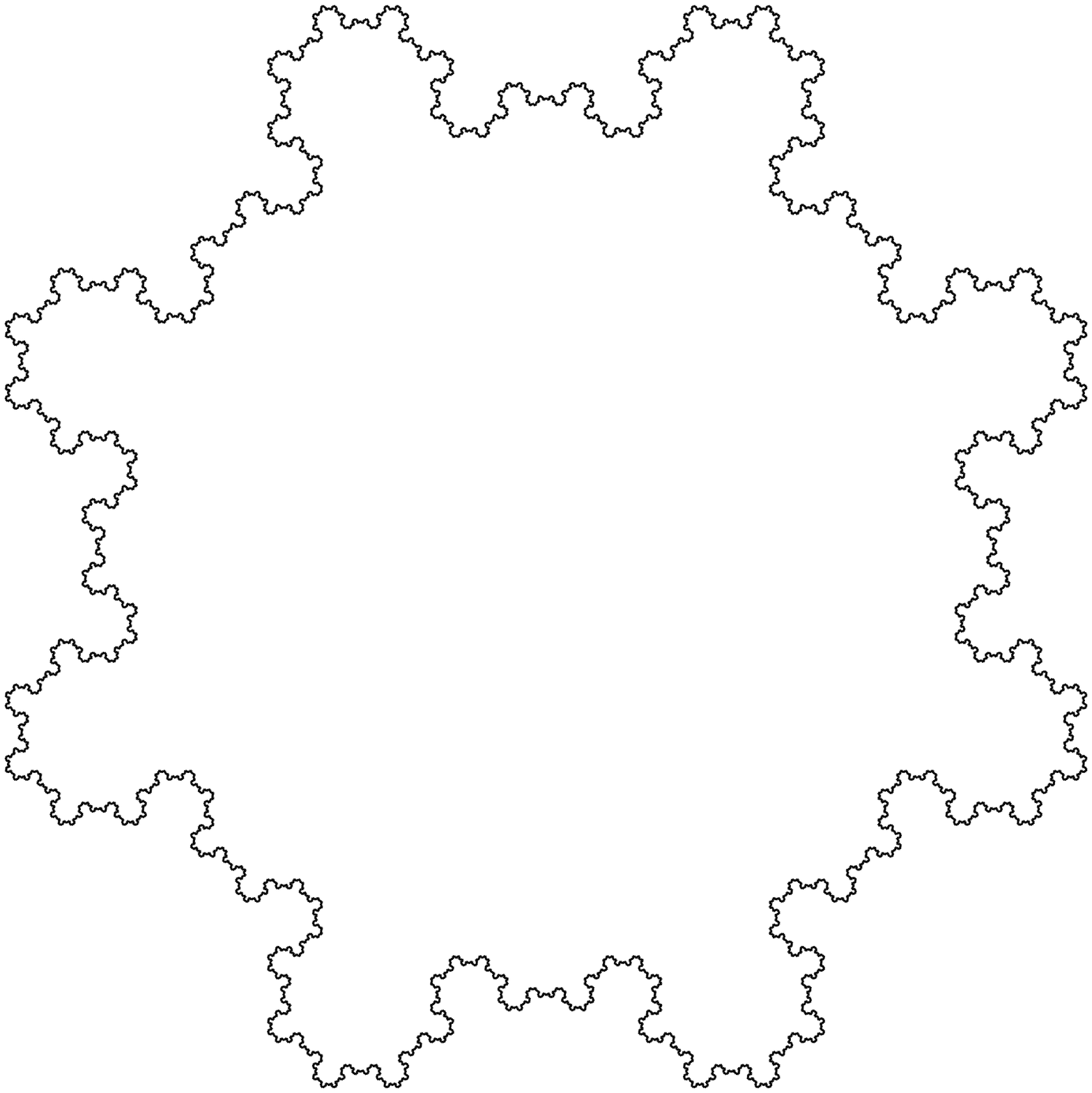}}
\newline
{\bf Figure 3.7:\/} The snowflake emerges.
\end{center}

Now we are going to show the $\pi_2$ projections
of the same orbits.  

\begin{center}
\resizebox{!}{2.6in}{\includegraphics{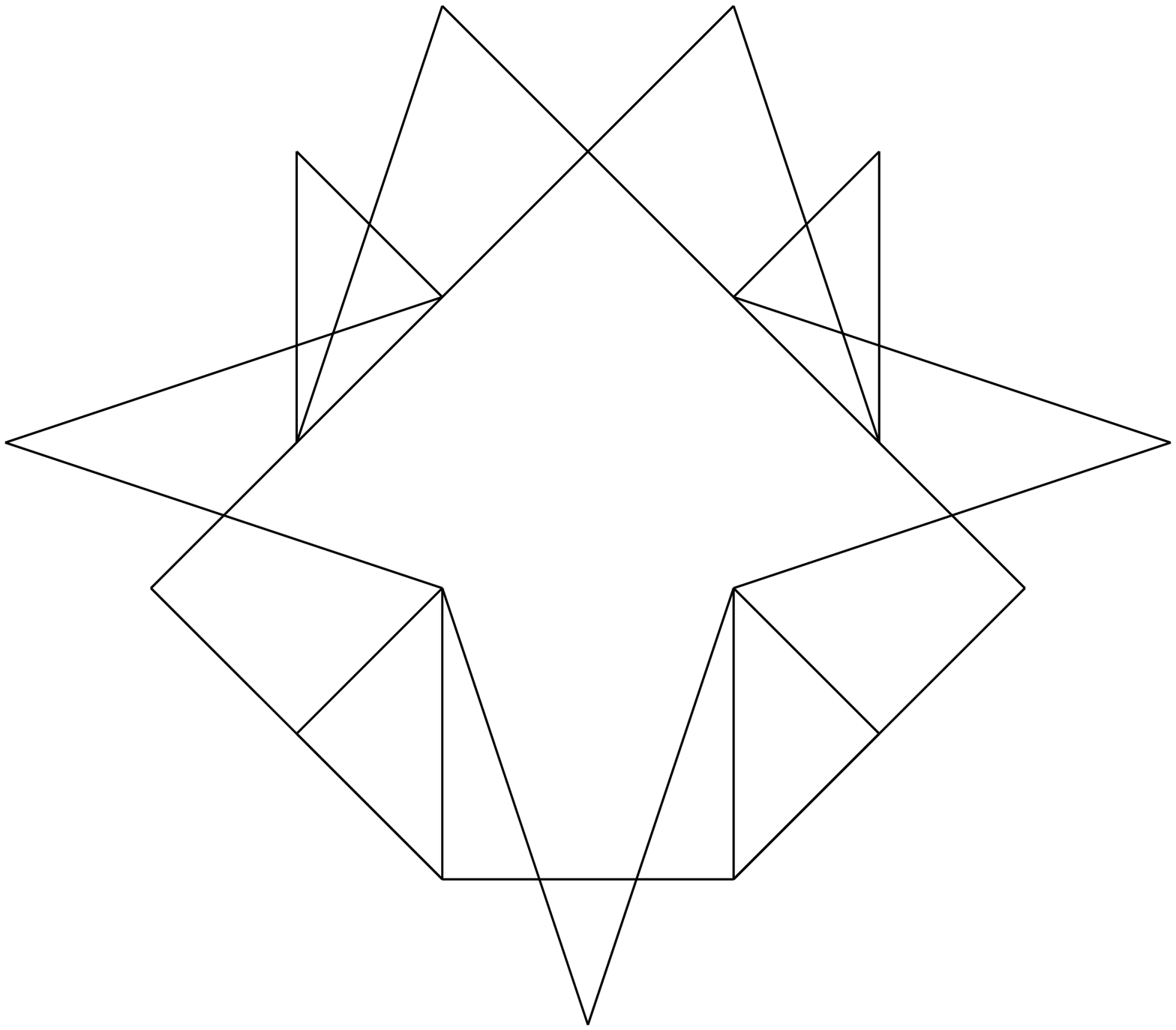}}
\resizebox{!}{2.6in}{\includegraphics{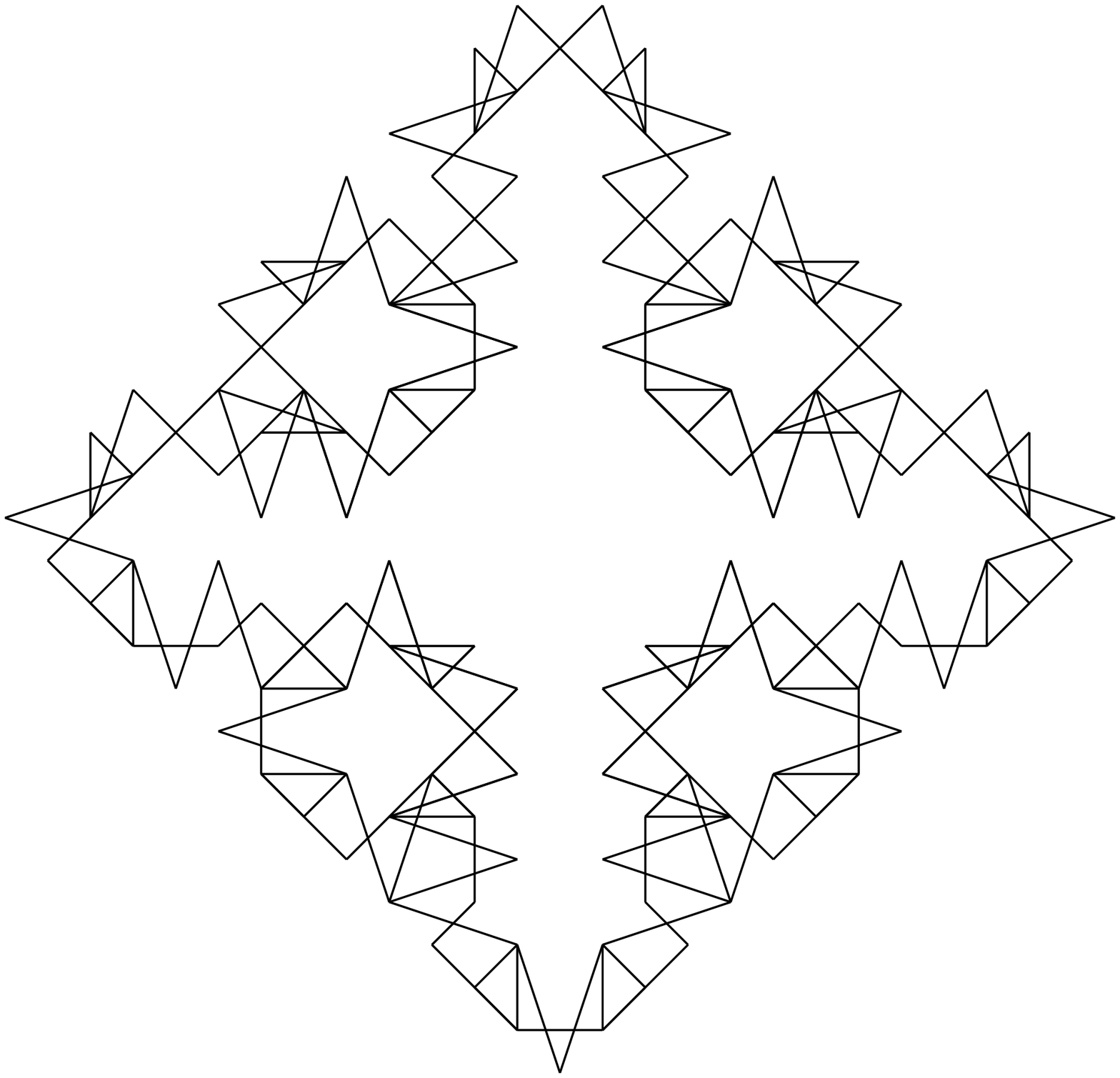}}
\newline
\resizebox{!}{2.6in}{\includegraphics{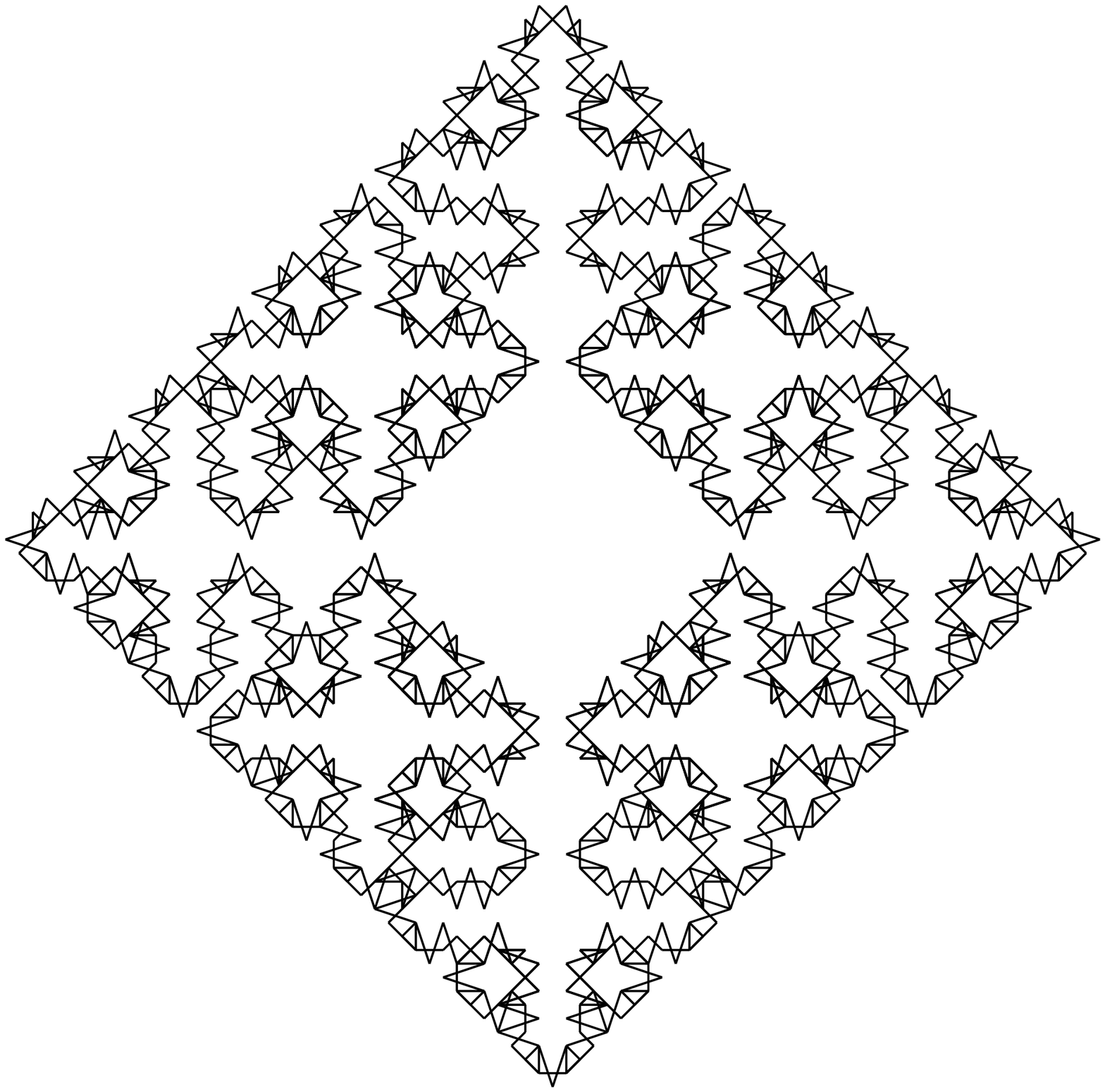}}
\resizebox{!}{2.6in}{\includegraphics{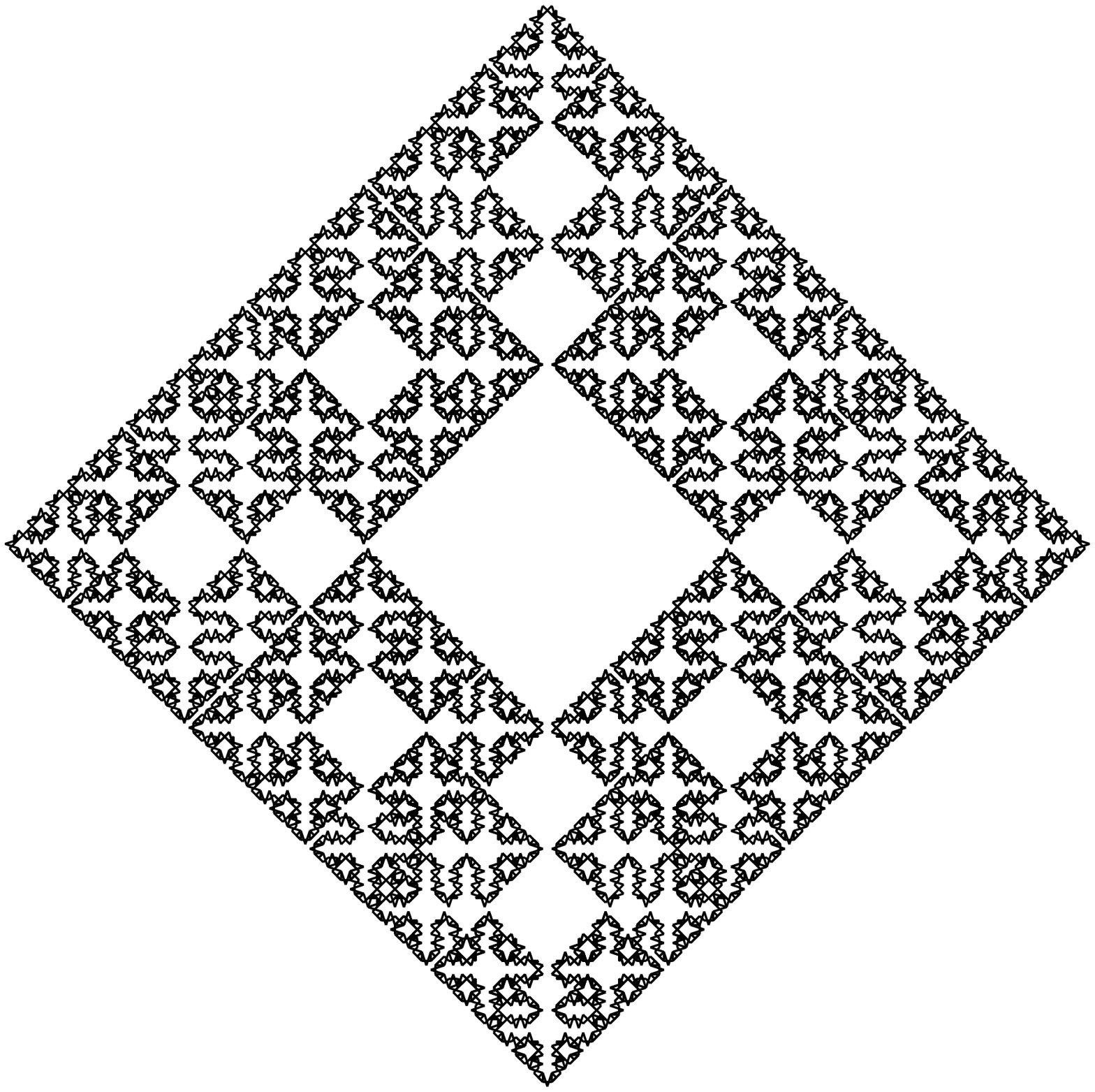}}
\newline
{\bf Figure 3.7:\/} The carpet emerges in the odd case.
\end{center}

For the even orbits, the $\pi_2$ projection is an open curve
that looks fairly similar.  We refer the reader to the program
OctoMap 2, where one can see much better pictures, in
both the odd and even cases.

\subsection{The Regular Heptagon}
\label{hept}

We can make the same constructions for the regular heptagon as we
made for the regular pentagon and octagon.  The analogue of Figure 3.3 exists,
but is not yet understood.  The two interesting projections
in this case are $\pi_2$ and $\pi_3$.   Unlike the two cases we
considered above, there is an explosion
of distinct pictures in this case.
Sometimes the $\pi_2$ projection looks nicer and sometimes
the $\pi_3$ projection looks nicer.  Sometimes both projections
look incomprehensibly complex.  Here are two simple examples,
both $\pi_2$ projections.

\begin{center}
\resizebox{!}{2.5in}{\includegraphics{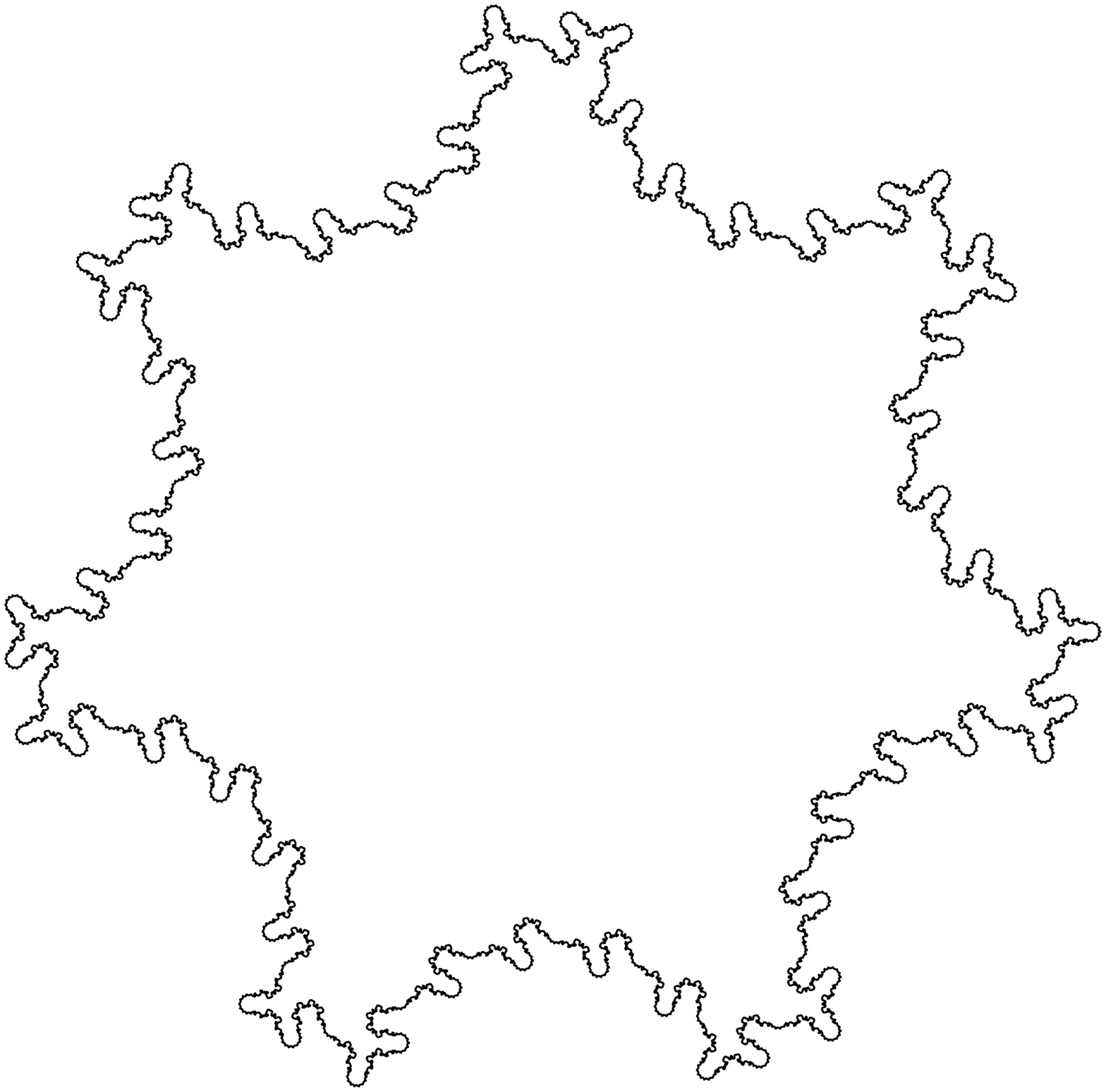}}
\resizebox{!}{2.5in}{\includegraphics{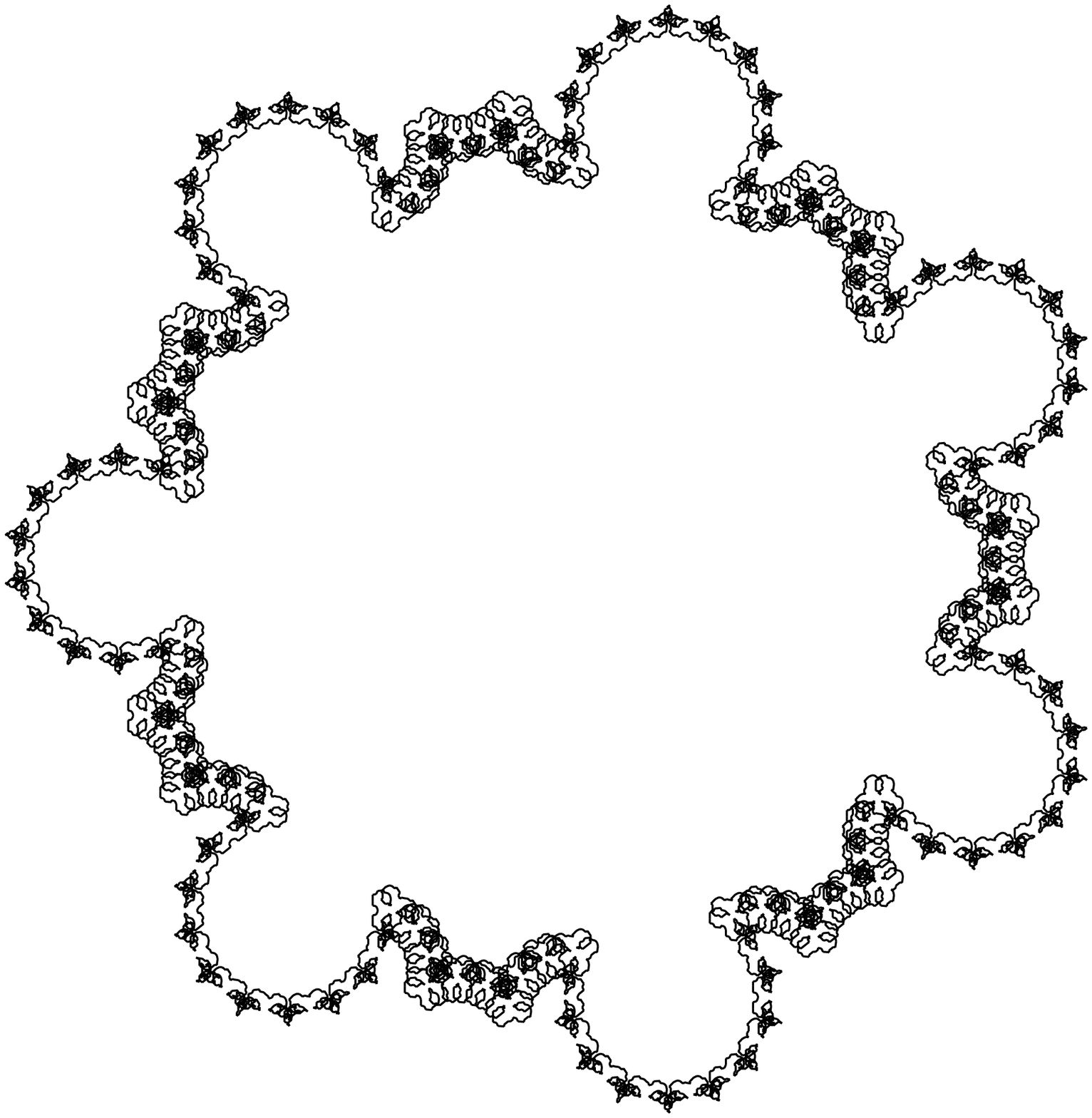}}
\newline
{\bf Figure 3.10:\/}  Some examples associated to the regular $7$-gon.
\end{center}

These examples
do not even begin to capture the vast
array of pictures one sees for the regular $7$-gon.
We encourage the experimentally-minded reader to 
explore the situation for themselves.

\newpage

\section{The Snowflake and the Carpet}
\label{fractal}

\subsection{The Hausdorff Topology}
\label{lim}

Let $K_1,K_2 \subset \R^2$ be two compact subsets.
We define
$\delta(K_1,K_2)$ to be the smallest $\epsilon$ such that
$K_1$ is contained in the $\epsilon$-tubular neighborhood
of $K_2$ and {\it vice versa\/}.   The function
$\delta$ is known as the {\it Hausdorff distance\/}
between $K_1$ and $K_2$.
Given any compact
subset $\Omega \subset \R^2$ we let
$M(\Omega)$ denote the set of compact subsets of $\Omega$.
Equipped with the function $\delta$, the space
$M(\Omega)$ is itself a compact metric space.

When we have a sequence of uniformly
bounded compact subsets $\{K_n\}$ and
we say that it converges, we mean that it
converges in $M(\Omega)$ for some large compact
$\Omega$ that contains all the individual sets.
The set $\Omega$ just serves as a kind of
container, so that we can speak about the
convergence as taking place within a
compact metric space.  The notion
of convergence we get is independent of
the choice of $\Omega$.
Indeed, we could say more simply that $K_n \to K$ iff
$d(K_n,K) \to 0$.  
We call $K$ the {\it Hausdorff limit\/}
of $\{K_n\}$.

One case of interest to is is when we have a
sequence $\{P_n\}$ of polygons whose
diameter tends to $\infty$.   If we can
a compact subset $K$ and a sequence of
similiarities $\{S_n\}$ such that
$S_n(P_n) \to K$ in the Hausdorff topology,
then we call $K$ a {\it rescaled limit\/} of
$\{P_n\}$.  In this case, the contraction
factor of $S_n$ necessarily tends to $0$.

\subsection{The Snowflake}
\label{snowflake}

We call the snowflake $I_3$.
Figure 4.1 shows the subdivision of a right-angled
isosceles triangle into $5$ smaller ones.   If the long side of
the triangle has length $1$, then the square on the right
hand side of the figure has side length $\sqrt 2-1$.
The three small triangles in the middle have the
same size.

\begin{center}
\resizebox{!}{1in}{\includegraphics{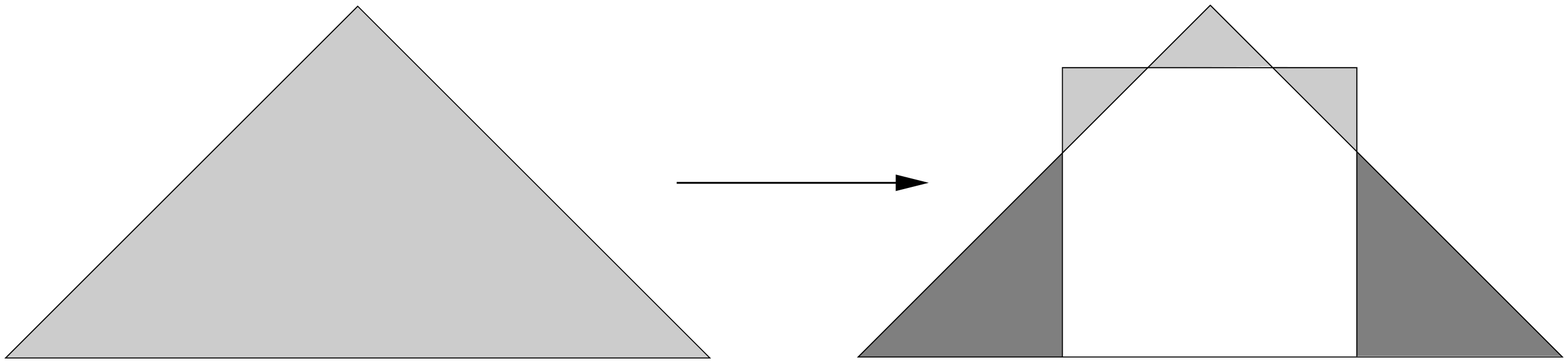}}
\newline
{\bf Figure 4.1:\/} Subdivision rule
\end{center}

To produce the snowflake, we start with the union of $8$ isosceles
triangles shown at the top right in Figure 4.2.
Then we apply the subdivision rule iteratively.
The Hausdorff limit is our snowflake. We denote this
limit by $I_3$.

\begin{center}
\resizebox{!}{2.2in}{\includegraphics{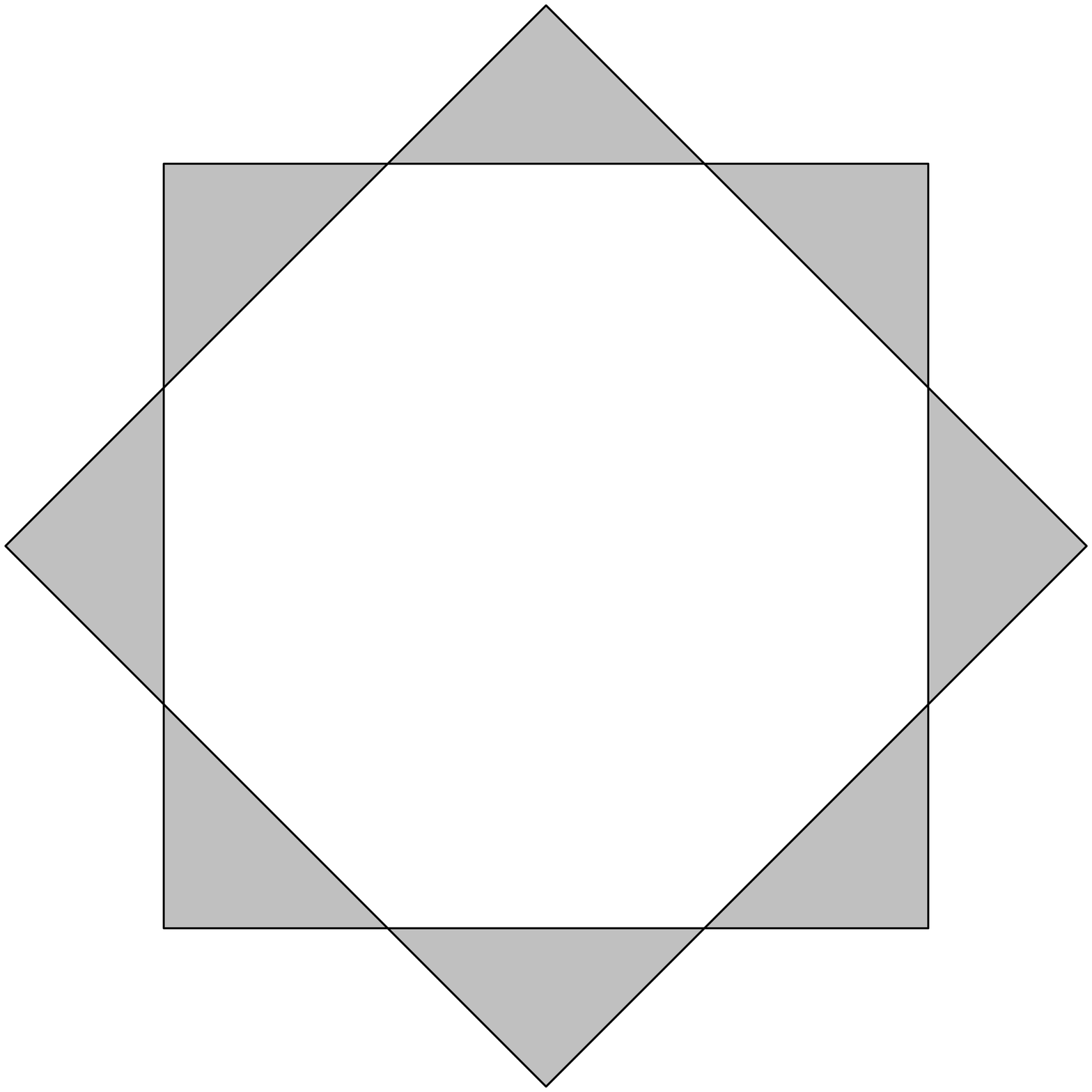}}
\resizebox{!}{2.2in}{\includegraphics{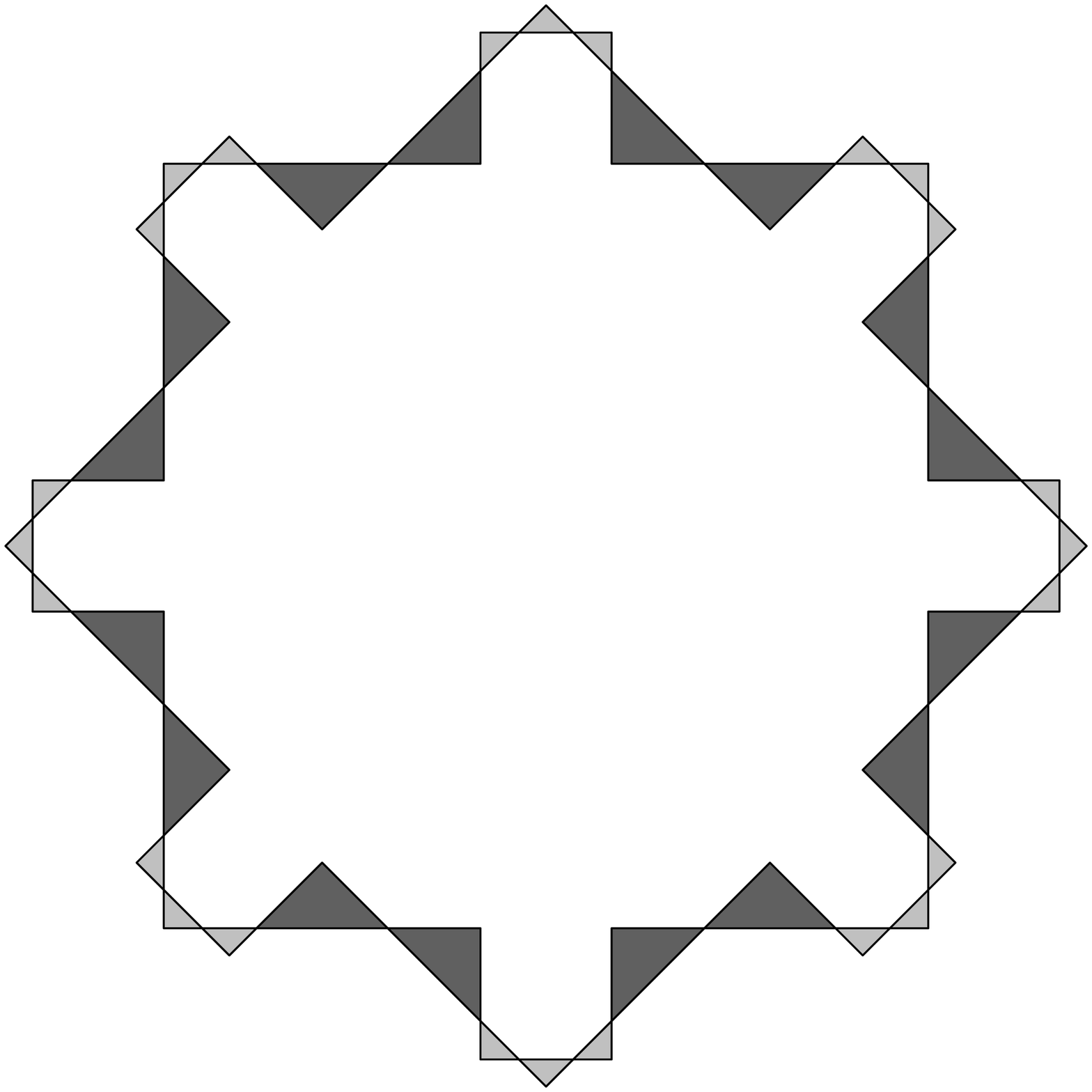}}
\newline
\resizebox{!}{2.3in}{\includegraphics{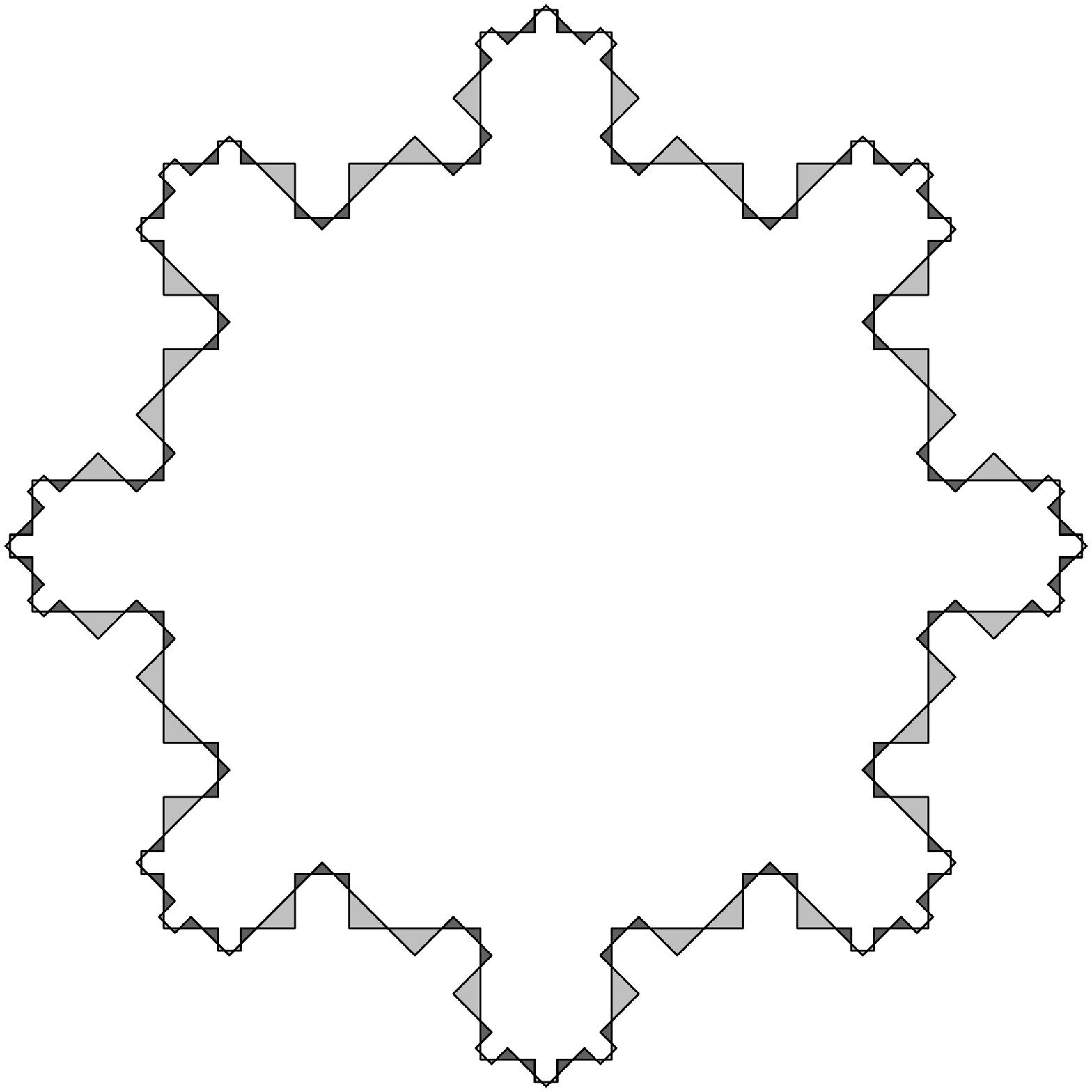}}
\resizebox{!}{2.3in}{\includegraphics{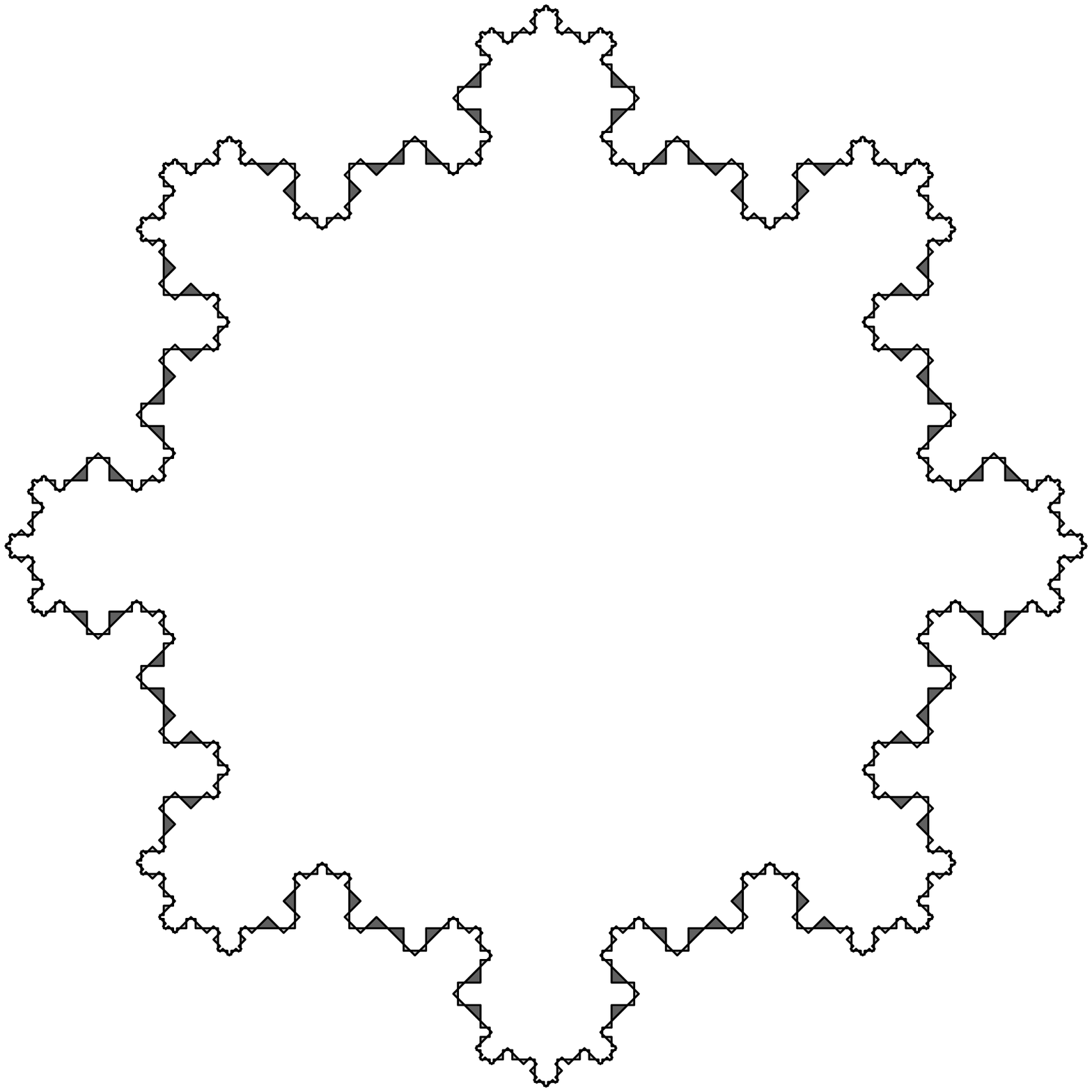}}
\newline
{\bf Figure 4.2:\/} Producing the snowflake by iteration.
\end{center}

One nice feature of our construction is that $I_3$
contains all the vertices of the triangles at every
stage of the construction, and the union of these
vertices is dense in $I_3$.   Moreover, all these
vertices lie in $\Z[\omega]$, where $\omega$ is
the usual primitive $8$th root of unity.
The triangles in our pictures are $2$ colored in a natural way,
according as to whether they point inward or outward.
Assuming that we label a triange $0$ if it points inward
and $1$ if it points outward, our subdivision rule
changes the colors according to the following rule:
\begin{equation}
0 \to 10001; \hskip 30 pt 1 \to 01110.
\end{equation}
We will see the significance of this structure below 
when we examine the carpet.

\subsection{A Hidden Symmetry of the Snowflake}
\label{hidden}

For any right isosceles triangle $T$, let
$L(T)$ denote the curve one obtains by
iterating the basic snowflake substitution rule and
taking a limit.   Assuming that $T$ has
a distinguished side $\sigma$, let
$L(T,\sigma)$ denote the portion
of $L(T)$ whose endpoints are the
endpoints of $\sigma$.

\begin{center}
\resizebox{!}{1.8in}{\includegraphics{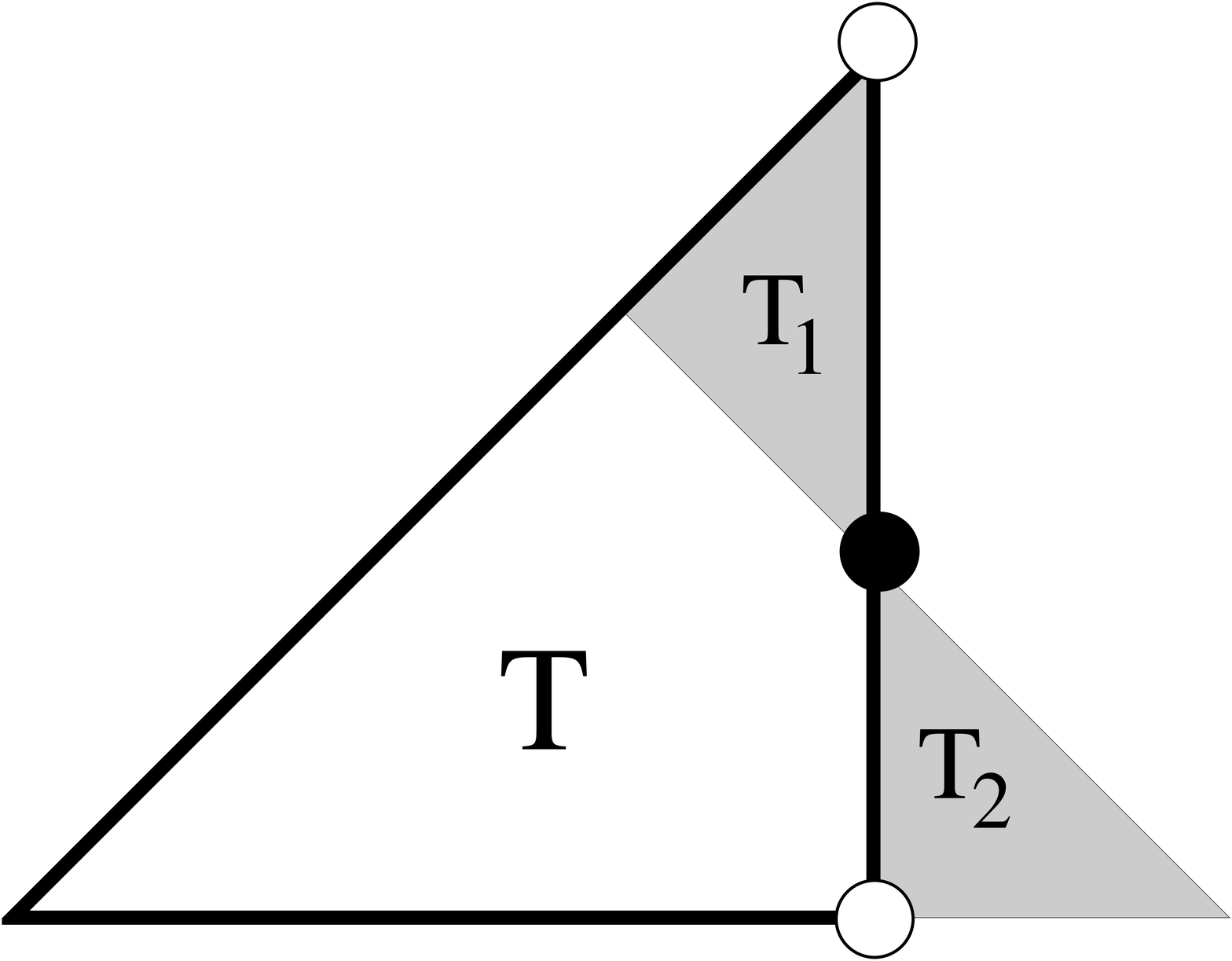}}
\newline
{\bf Figure 4.3:\/} Some special triangles
\end{center}

Let $T$, $T_1$, and $T_2$ be the three right isosceles triangles
shown in Figure 4.3.  The two triangles $T_1$ and $T_2$ have
the same size as each other, and are $\sqrt 2-1$ times as large
as $T$.   Let $A=L(T,\sigma)$, where $\sigma$ is the
edge of $T$ bounded by the two white vertices.
Let $A_1=L(T_1,\sigma_1)$, where $\sigma_1$ is
edge of $T_1$ bounded by the white and black
vertex.  Define $A_2=L(T_2,\sigma_2)$ similarly.

\begin{lemma}
\label{hiddenlemma}
$A= A_1 \cup A_2$.
\end{lemma}

\startproof
Let $d$ be the Hausdorff distance 
between $A$ qnd $A_1 \cup A_2$.
The substitution of $T$ consists of $5$ smaller triangles,
$T_1',...,T_5'$, with $T_1'=T_1$ and
$T_2' \cup T_3'$ related to $T_2$ just as
$T_1 \cup T_2$ is related to $T$.   
There are edges $\sigma_2'$ and $\sigma_3'$ of
$T_2'$ and $T_3'$ such that
$\sigma_2=\sigma_2'\cup \sigma_3'$.
The remaining triangles $T_4'$ and $T_5'$ lie
to the left of the bottom white vertex of $T$.

Since $T_1=T_1'$, we have $A_1 \subset A$.
At the same time, the pair of arcs
$(A,A_1 \cup A_2)$ is similar to the pair
of arcs $(A_2,A_2' \cup A_3')$, where
$A_2'=L(T_2',\sigma_2')$ and
$A_3'=L(T_3',\sigma_3')$.
From this we see that the Hausdorff distance
between $A$ and $A_1 \cup A_2$ is the
same as the Hausdorff distance between
$A_2$ and $A_2' \cup A_3'$.  But, by
scaling, the latter distance is $(\sqrt 2 -1)d$.
So, we have the equation
$$d=(\sqrt 2-1)d,$$
which of course forces $d=0$.
\endproof

\subsection{The Carpet}
\label{carpet}

We call the carpet $I_2$.
The carpet is produced by substitution rule that is combinatorially
identical to the one that produces the snowflake.  This time,
we have two shapes, parallelograms and trapezoids.  The initial
shapes have their vertices in $\Z^2$.  Our figures below are
accurate.

Figure 4.4 illustrates how the parallelogram $P$ is replaced by
the sequence 
$$S_P(P)=T_1 \cup P_2 \cup P_3 \cup P_4 \cup T_5$$ 
of $5$ smaller parallelograms and trapezoids.

\begin{center}
\resizebox{!}{1.8in}{\includegraphics{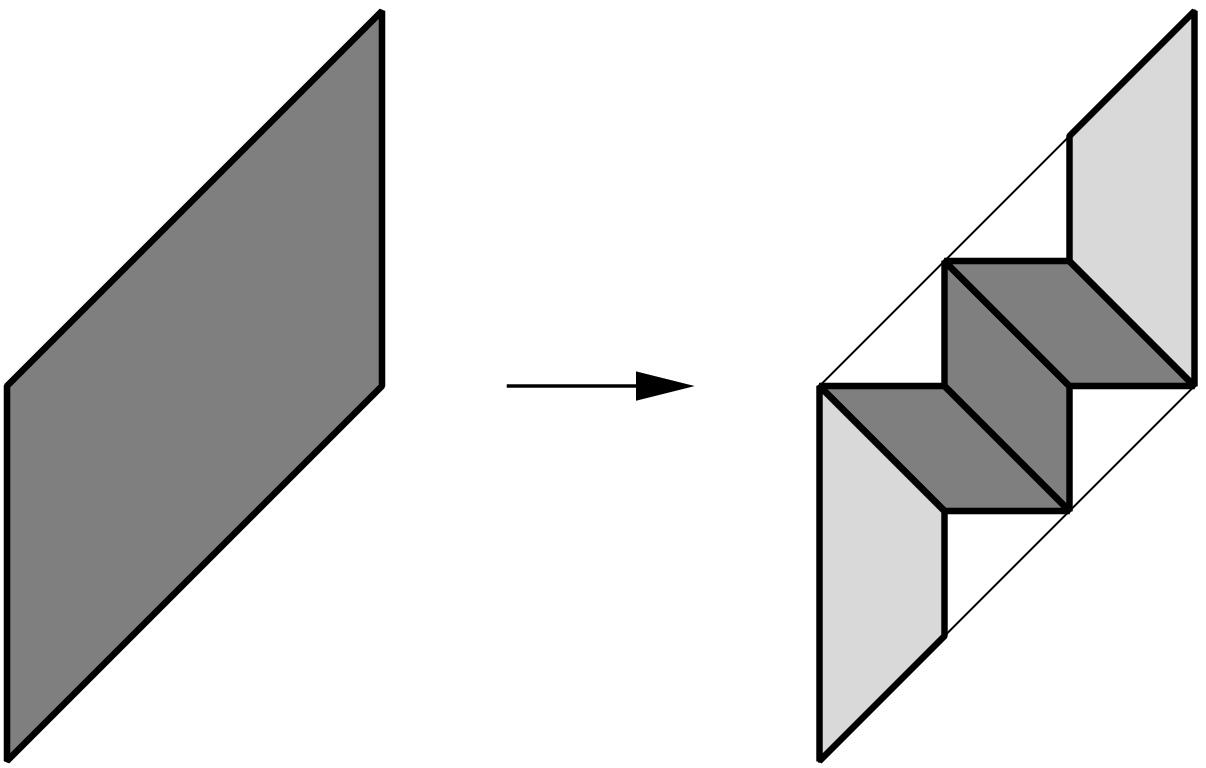}}
\newline
{\bf Figure 4.4:\/} Substitution rule for the parallelogram
\end{center}

Figure 4.5 shows how the trapezoid $T$ is replaced by the sequence
$$S_T(T)=P_1 \cup T_2 \cup T_3 \cup T_4 \cup P_5$$ 
of smaller parallelograms and trapezoids.
We say {\it sequence\/} here rather than {\it union\/} because
the pieces in the sequence $S_T(T)$ overlap each other.
To make the substitution clear, we first
draw $P_1 \cup T_2 \cup T_3$.
Then we add $P_4$ and $P_5$.   The partition is
invariant under reflection in the vertical line of symmetry.

\begin{center}
\resizebox{!}{1.8in}{\includegraphics{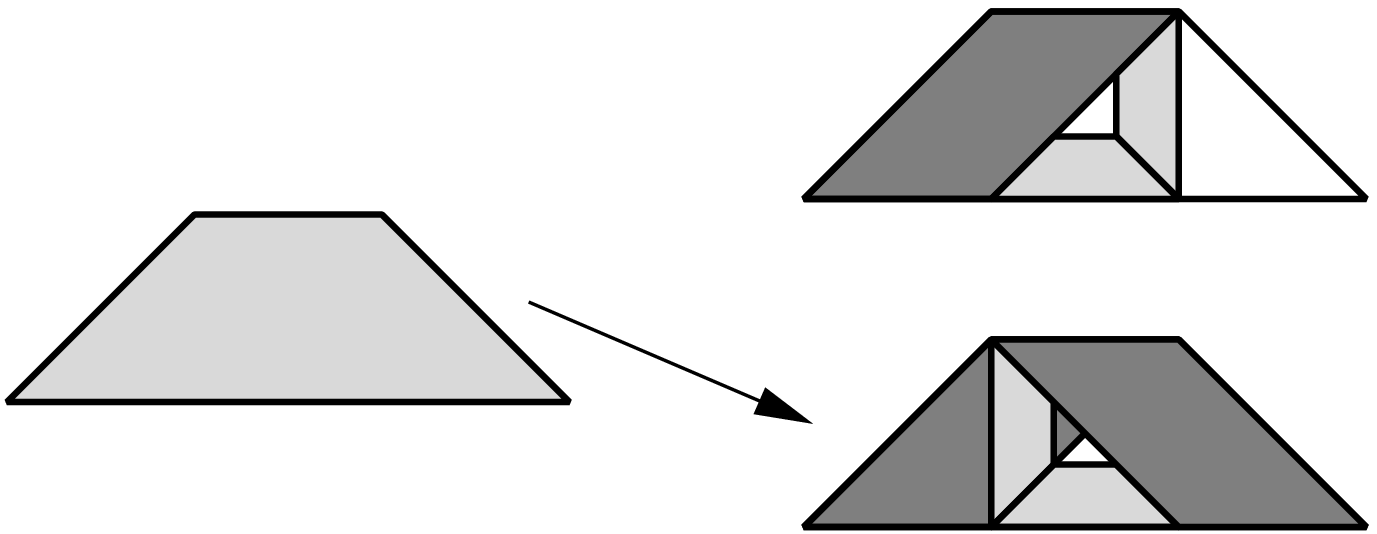}}
\newline
{\bf Figure 4.5:\/} Substitution rule for the trapezoid
\end{center}

Let $X$ be a cyclically ordered
finite list of parallelograms and trapezoids,
similar to the ones above.
We let $|X|$ be the union of shapes in $X$.
We let $S_P(X)$ denote the result of substituting
$S_P(P)$ for all parallelograms $P \in X$.
Likewise we define $S_T(X)$.  Note that
$S_P(X)$ and $S_T(X)$ both have natural cyclic
orders, inherited from the ordering on $X$.
Our carpet is
\begin{equation}
\label{2lim}
I_2=\lim_{n \to \infty} |I_2(n)| \hskip 30 pt
I_2(n)=S_P(S_T(I_2(0)).
\end{equation}
Here $I_2(0)$ is a union of $8$ copies of $T$ arranged in a square
pattern that winds around twice.  (We use $8$ rather than $4$
so that the seeds for $I_2$ and $I_3$ have the same cardinality.)
Figure 4.6 shows the sets $I_2(n)$ produced by iterated
substitution, for $n=0,1,2,3$.   Since the pieces
in overlap, there is more than one partition we could
draw.  We have chosen the partition in which every
trapezoid is drawn above every parallelogram. 

\begin{center}
\resizebox{!}{2.3in}{\includegraphics{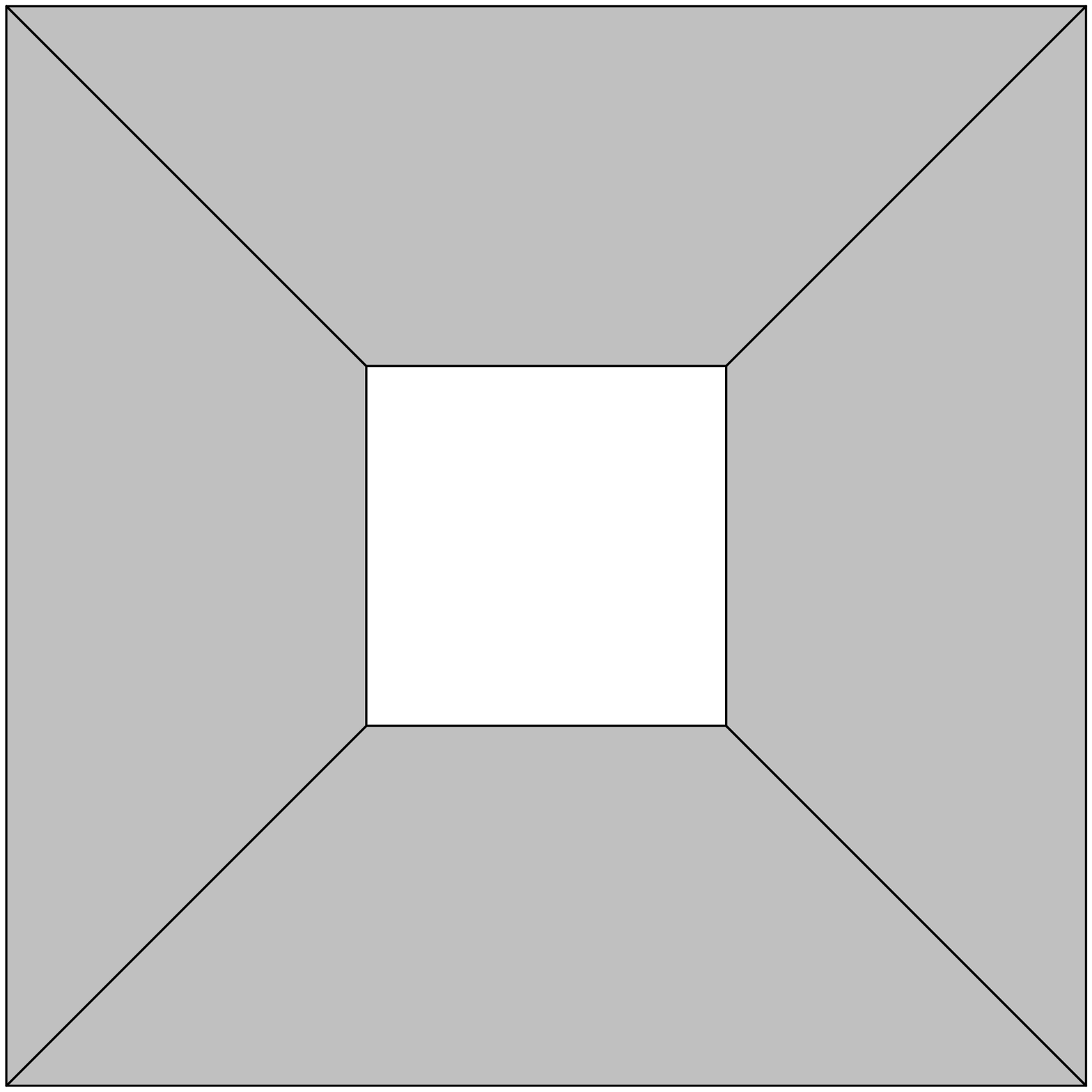}}
\resizebox{!}{2.3in}{\includegraphics{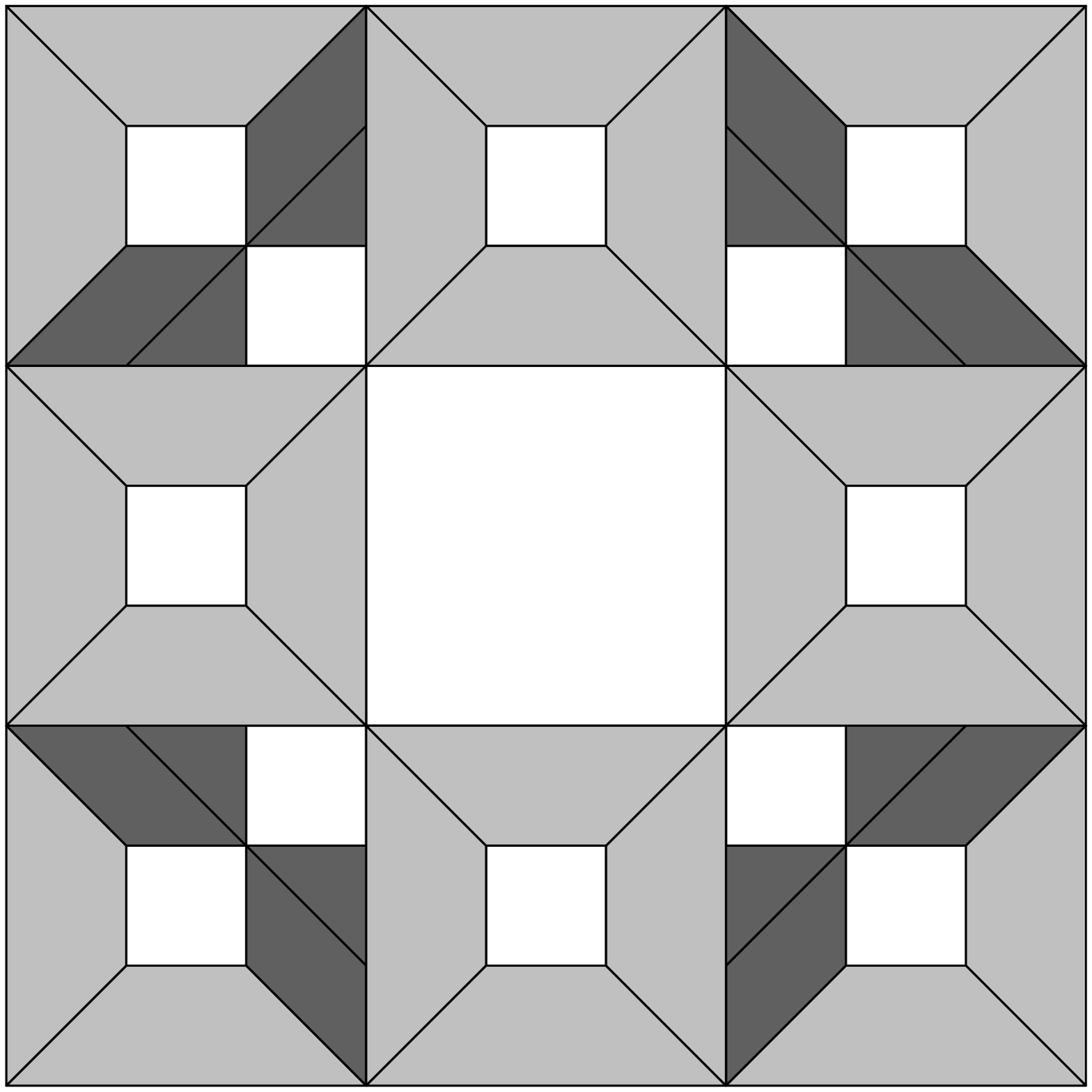}}
\newline
\resizebox{!}{2.3in}{\includegraphics{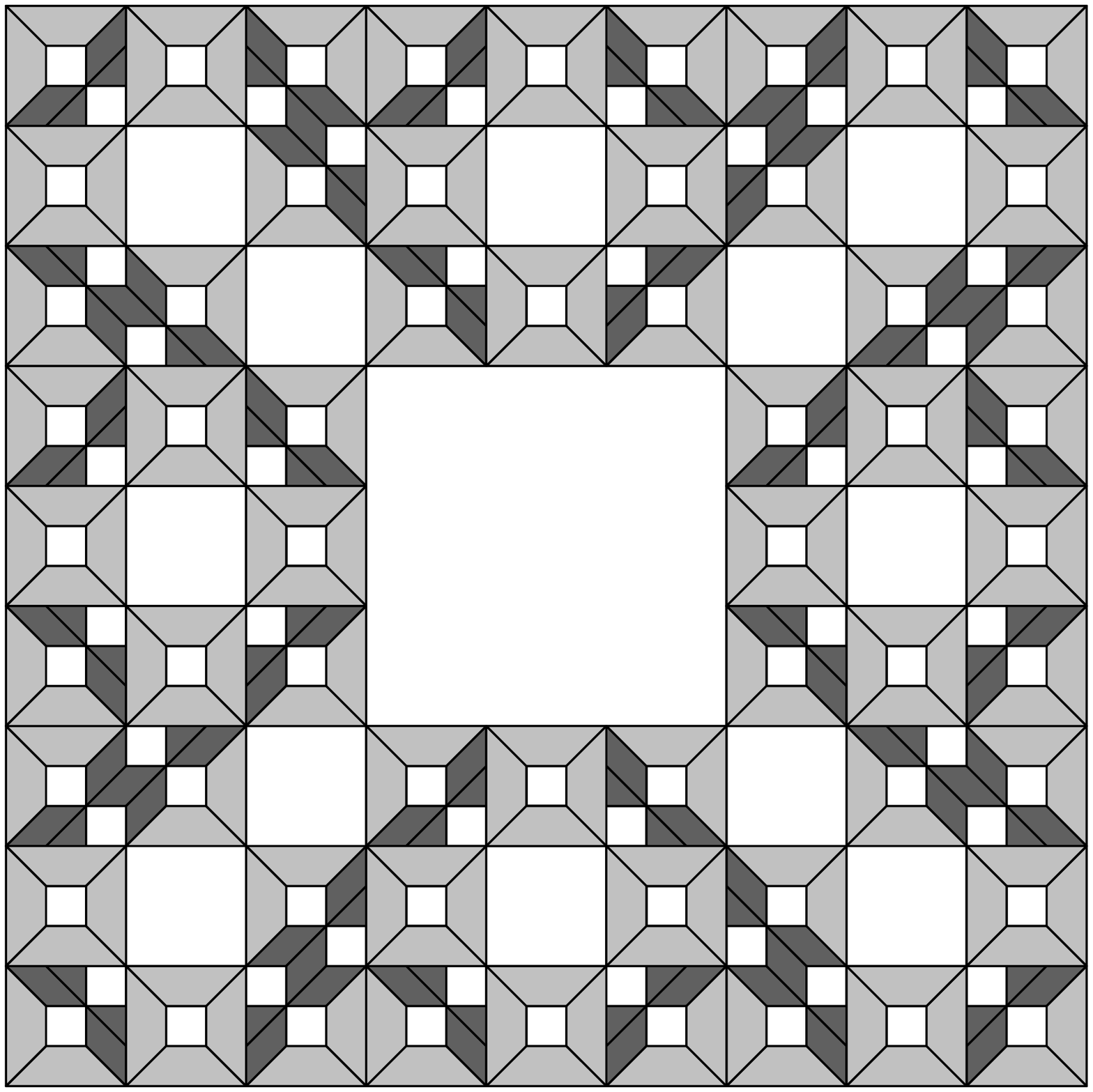}}
\resizebox{!}{2.3in}{\includegraphics{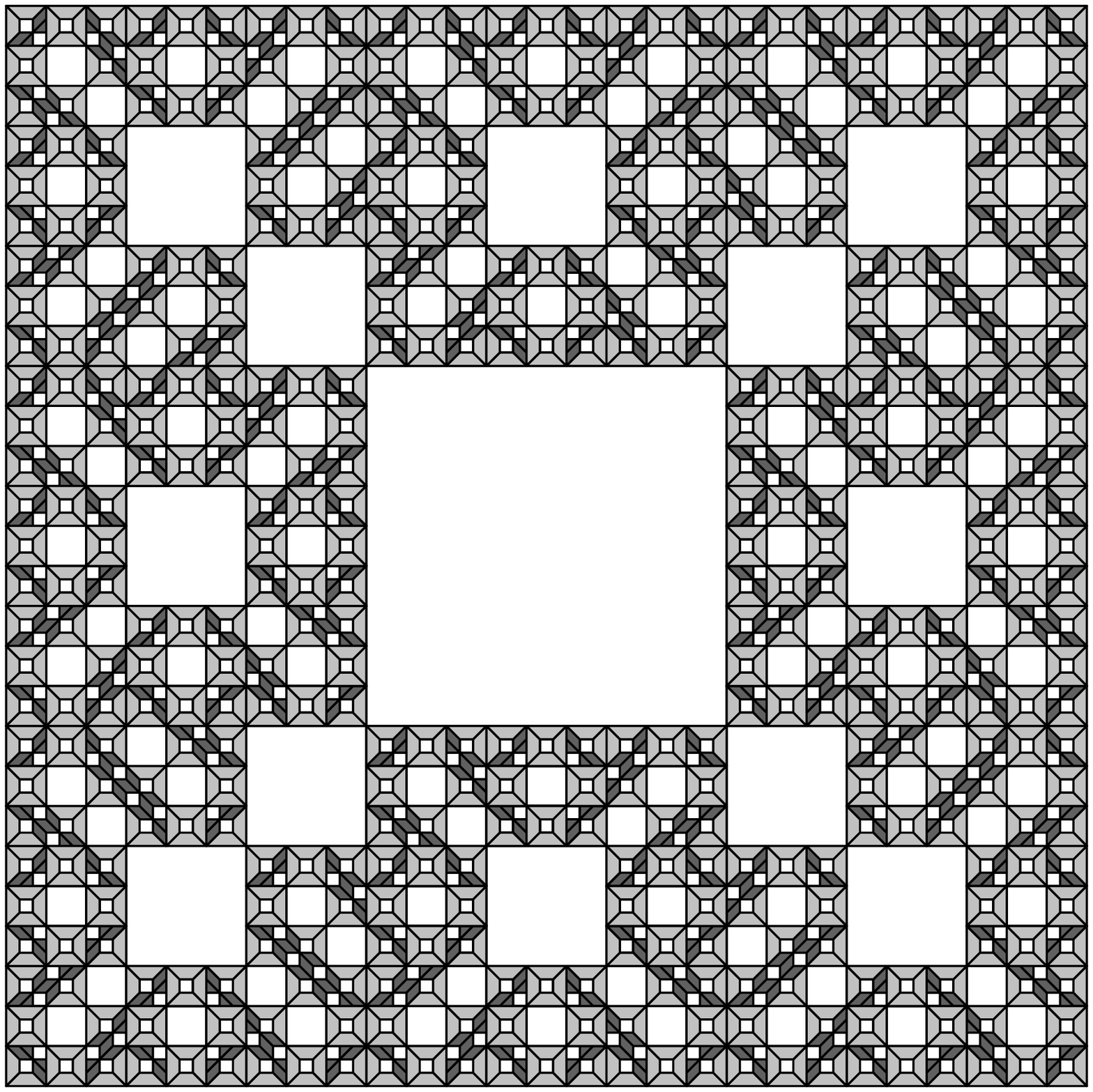}}
\newline
{\bf Figure 4.6:\/} Iterating the substitution rule
\end{center}

\subsection{Another View of the Carpet}
\label{view2}

Figure 4.7  shows another substitution rule that generates
$I_2$.   In Figure 4.7, we are showing how to replace
a square by a subset of a square.  Here we think of the
dark triangles and quadrilaterals as markings that tell
us how to orient the pieces when we iterate the basic rule.
We obtained this rule by looking at Figure 4.6.
The fractal obtained from this alternate rule is a
subset of the Sierpinski carpet.  It is obtained from
the Sierpinski carpet by systematically deleting
certain squares.  

\begin{center}
\resizebox{!}{4.5in}{\includegraphics{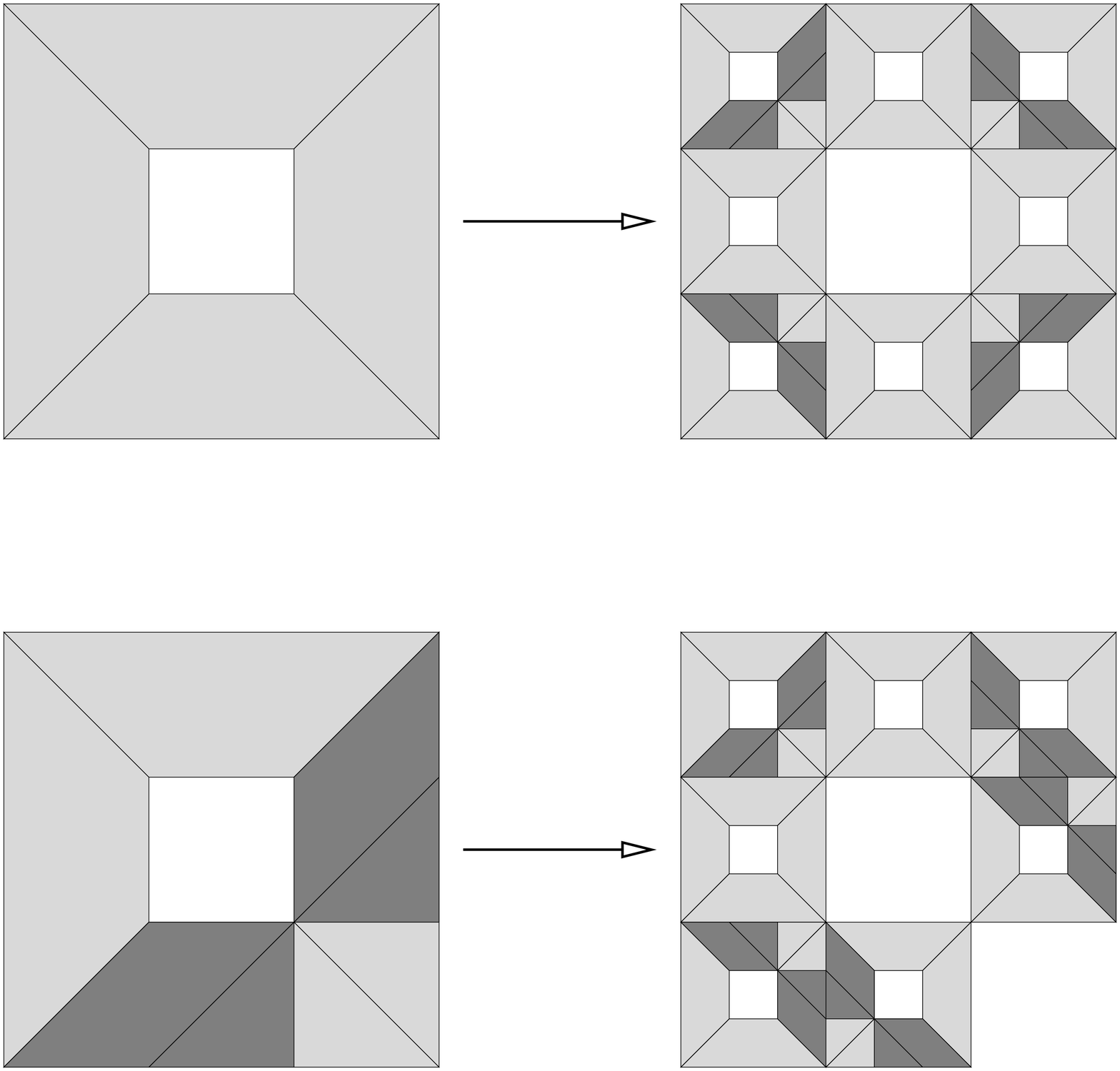}}
\newline
{\bf Figure 4.7:\/} An alternate substitution rule
\end{center}

\begin{lemma}
$I_2$ is the fractal produced by the
alternate substitution rule.
\end{lemma}

\startproof 
We will
start with the original rule  and keep modifying it until
we get to the alternate rule.  Along the way, we show that
the modifications don't change the limiting fractal.

For starters, we divide each parallelogram in half.   In this way,
we interpret the original substitution rule as a being generated
by a rule for trapezoids and a rule for isosceles triangles.
The left hand side of Figure 4.8 shows the new $S_P(P)$,
The right hand side of Figure 4.8 shows the new
$S_T(T)$.

\begin{center}
\resizebox{!}{1in}{\includegraphics{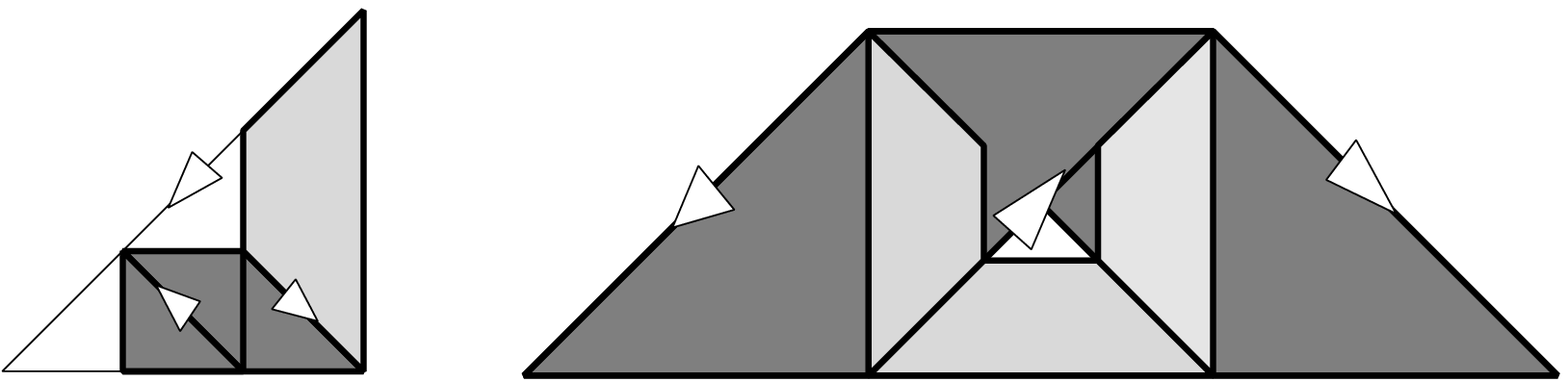}}
\newline
{\bf Figure 4.8:\/} The first modification
\end{center}

We have redrawn the second halves of Figures 4.4 and 4.5.  In
the case of Figure 4.5, we have drawn all the trapezoids on top.
In drawing this piece, we arbitrarily chose to draw the left
(divided) parallelogram on top of the right one.  The substitution
rule for the isosceles triangle does not respect the symmetry
of the triangle.  We record the asymmetry by orienting the
hypotenuse of each triangle, as shown.
 The new
substitution rule clearly generates $I_2$.    We still call our
new rules $S_R$ and $S_T$.  Here $R$ stands for a right
isosceles triangle and $T$ still stands for the trapezoid.

Figure 4.9 shows two successive modifications of the rule
for $S_T$.

\begin{center}
\resizebox{!}{.85in}{\includegraphics{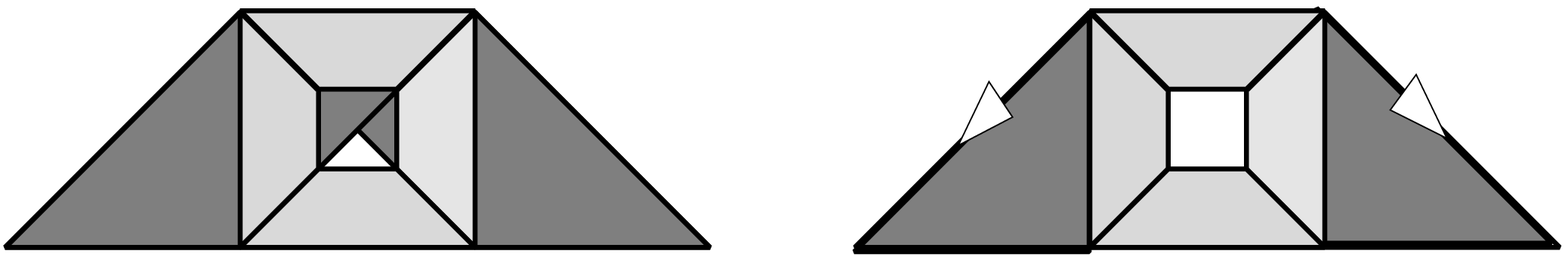}}
\newline
{\bf Figure 4.9:\/} Two modifications of $S_T$.
\end{center}

For the first modification, we add another trapezoid to $S_T(T)$.
This new trapezoid appears in $S_R(S_T(T))$, so the new rule
rule generates the same limit.    We give the same names to
our new substitution rules.    For the second modification, we
delete two of the isosceles triangles.  Call the deleted
triangles $R_1$ and $R_2$.   Observe that
$S_R(R_k) \subset S_T(T)$ for $k=1,2$.  Therefore,
the latest subdivision rule has the same limit as the
original.  We again call the generators of the
rule $S_T$ and $S_R$.  This time $S_T$ and $S_R$
replace one piece by a union of pieces that have
pairwise disjoint interiors.  

From here it is a simple exercise to check that $S_R \circ S_T$
implements the rules shown in Figure 4.7.  
\endproof

\subsection{The Map from the Snowflake to the Carpet}
\label{map}

The initial seed for $I_3$ is a cyclically
ordered list $I_3(0)$ of $8$ isosceles triangles $\tau_1,...,\tau_8$.
The initial seed for $I_2$ is a cyclically ordered list $I_2(0)$
of $8$ trapezoids $T_1,...,T_8$, which wrap around twice,
as we have already mentioned.
We have a correspondence $T_i \leftrightarrow \tau_i$ between
the two seeds.

Since the two substitution rules are combinatorially identical,
we inductively get a surjective map from the triangles
in $I_3(n)$ to the quadrilaterals in $I_2(n)$.
This map respects the cyclic ordering on both sets,
and it also respects the relation of {\it descendence\/}.
Say that a shape $\tau_n$ of $I_j(n)$ is a
{\it descendent\/} of a shape $\tau_m$ of
$I_j(m)$ if the substitution rule, applied
$n-m$ times, to $\tau_m$, produces a list of
shapes that contains $\tau_n$.  

We define $\phi: I_3 \to I_2$ by taking the
limit of our  correspondence.   Given $x \in I_3$ we
can find a sequence $\{\tau_n\}$ of triangles,
with $\tau_n \in X_3(n)$, that accumulates on $x$.
We define $\phi(x)$ to be the limit of the corresponding
quadrilaterals in $I_2$.    

\begin{lemma}
$\phi$ is well-defined.
\end{lemma}

\startproof
Suppose that $\{\tau_n'\}$ is another sequence of triangles
that converges to $x \in I_2$.  Let $S_n$ and $S_n'$ 
be the shapes of $I_2(n)$ corresponding
respectively to $\tau_n$ and $\tau_n'$.  
It suffices to prove that ${\rm diam\/}(S_n \cup S_n') \to 0$.

Call two triangles in $I_3(m)$ {\it close\/} if they are either
identical or else consecutive in the cyclic order.
Let $d_n$ be the largest integer such that
$\tau_n$ and $\tau_n'$ are descendents of
close triangles in $I_3(d_n)$.   Since the distance
between $\tau_n$ and $\tau_n'$ converges to
$0$, and $I_3$ is evidently an embedded curve,
$d_n \to \infty$.

We make the same definitions as above for $I_2$.
We can find close shapes $\Sigma_n$ and
$\Sigma_n'$ in $I_2(d_n)$ such that $S_n$ is a
descendent of $\Sigma_n$ and $S_n'$ is
a descendent of $\Sigma_n'$.
Note that $\Sigma_n$ and $\Sigma_n'$ have intersecting
boundaries in all cases.   Also, the diameters of these shapes
tends to $0$ with $n$.  Hence 
${\rm diam\/}(\Sigma_n \cup \Sigma_n') \to 0$.
Finally, we have $S_n \subset \Sigma_n$ and
$S_n' \subset \Sigma_n'$.  
Hence ${\rm diam\/}(S_n \cup S_n') \to 0$ as desired.
\endproof

Essentially the same proof shows that $\phi$ is continuous.
Since all the approximating correspondences are
surjective, $\phi$ is surjective as well.
\newline
\newline
{\bf Remark:\/}
Note that we could get an equally canonical map if we compose
$\phi$ with some isometry of $I_2$ or $I_3$.  This basic (though
trivial) ambiguity comes up in Statement 3 of our Main Theorem below.

\subsection{Geometry of the Map}

Now we give some geometrical information about the map
$\phi: I_3 \to I_2$ defined in the previous section.
Note that $I_3$ contains the vertices of triangles in
$I_3(n)$ for all $n$.  Say that a {\it special point\/}
of $I_3$ is a right-angled vertex of $I_3(n)$ for
some $n$.  The set $A_3$ of special points is
dense in $I_3$.

Say that a {\it special point\/} of $I_2$ is a midpoint of
the long edge of some trapezoid in $I_2(n)$.  We could
equally define the special points to be the centers of
the parallelograms.  The two definitions coincide.
Let $A_2$ denote the set of special points of $I_2$.
The set $A_2$ is the set of midpoints of the edges of
the countable family of squares that arises from the
substitution rule explained in \S \ref{view2}.

\begin{lemma}
$\phi(A_3)=A_2$.
\end{lemma}

\startproof
Let $v \in A_3$ be some point.  Then $v$ is the
right-angled vertex of a triangle $\tau_0$ of $I_3(n)$ for
some $n$.   Note that $v$ is also the right vertex
of a triangle $\tau_k$ of $I_3(n+k)$ for $k=1,2,3...$
The triangles $\tau_1,\tau_2,...$
all have the same type, either inward pointing
or outward pointing.

Consider the outward pointing case first.
In this case the quadrilateral $T_k$ of
$I_2(n+k)$ corresponding to
$\tau_k$ is a trapezoid.   Inspecting
our subdivision rule, we see that 
$T_{k+1}$ is the middle trapesoid
of $S_T(T_k)$.   Looking at Figure 4.5, we see that
 $T_k$ and $T_{k+1}$ have
a common distinguished point.  This holds for all
$j$.  Hence $\bigcap T_k=\phi(v)$ is the
common distinguished point of all these trapezoids.
Hence $\phi(v) \in A_2$.

In the inward pointing case, the quadrilateral
$P_k$ corresponding to $\tau_k$ is a parallelogram.
Here $P_{k+1}$ is the middle parallelogram of
$S_P(P_k)$.  But these parallelograms have a common
center.   The intersection
$\bigcap P_k=\phi(v)$ is this common center.
Looking at Figure 4.5 again, we see that the centers
of the parallelograms coincide with distinguished
points of trapezoids.   So, again
$\phi(v) \in A_2$.

Our argument so far shows that $\phi(A_3) \subset A_2$.
Any $x \in A_2$ is the distinguished point of some
trapezoid $T$ of some $I_2(n)$.  But then there is
a triangle $\tau$ of $I_3(n)$ which corresponds
to $\tau$.  By construction $x=\phi(v)$, for the
right vertex $v$ of $\tau$.  Therefore $\phi(A_3)=A_2$.
\endproof

Let $B_3 \subset I_3$ denote the set of points
that are acute vertices of triangles in $I_3(n)$
for all the different $n$.  Then $B_3$ is dense in $I_3$.
Let $B_2 \subset I_2$ denote the set of vertices
of quadrilaterals in $I_2(n)$ all the different $n$.
The set $B_2$ is the set of corners of the distinguished
squares mentioned above in connection with $A_2$.
The curve $I_3$ looks ``oscillatory'' in the
neighborhood of any $x \in B_3$.  
Any line through $x$ intersects 
$I_3$ infinitely often in every neighborhood of $x$.

\begin{lemma}
$\phi(B_3)=B_2$.
\end{lemma}

\startproof
The proof is similar to the one for the
previous result, so we will be a bit sketchy.
Any $x \in B_3(n)$ is the acute vertex of a triangle
$\tau_k$ in $I_3(n+k)$.   
The triangle $\tau_{k+1}$ is either the first or last triangle
in the list of $5$ triangles in $I_3(n+k+1)$ that replaces
$\tau_k$.   Whether $\tau_{k+1}$ is ``first'' or ``last''
is independent of $k$.   

Let $T_k$ be the quadrilateral of $I_2(n+k)$ corresponding
to $\tau_k$.  By construction, these quadrilaterals all share a
common vertex, and $\phi(x)$ must be this common vertex.
This shows that $\phi(B_3) \subset B_2$.   The proof that
$B_2 \subset \phi(B_3)$ is just as in the previous result.
\endproof

\noindent
{\bf Remark:\/}
It is worth pointing out one more feature of the map $\phi$:  It
is far from $1$ to $1$.   The set $I_2(n)$, considered as a
list of quadrilaterals, contains multiple copies of the same
trapezoid.  Many trapezoids appear $2^n$ times in $I_2(n)$.
From this we see that, for any $n$, there are $2^n$ segments
of $I_3$ that $\phi$ identifies.    For this reason, the
canonical surjection from $I_3$ to $I_2$ is a many-to-one
map.   We think if it is a kind of fractal version of the
universal covering map from the line to the circle.

\subsection{A Hidden Symmetry of the Carpet}
\label{hidden2}

The snowflake has a hidden symmetry, as discussed in
\S \ref{hidden}. Given the close correspondence between
the carpet and the snowflake, we would expect a
combinatorially identical hidden symmetry to appear
for the carpet.   This is indeed the case.  In this section
we describe the symmetry and sketch the proof.

The hidden symmetry of the carpet manifests itself
in two ways, and this at first seems different from
what happens with the snowflake.  However, in the
case of the snowflake, there are really two kinds
of right-angled isosceles triangles, so the symmetry
there actually does arise in two ways.    

\begin{center}
\resizebox{!}{3.2in}{\includegraphics{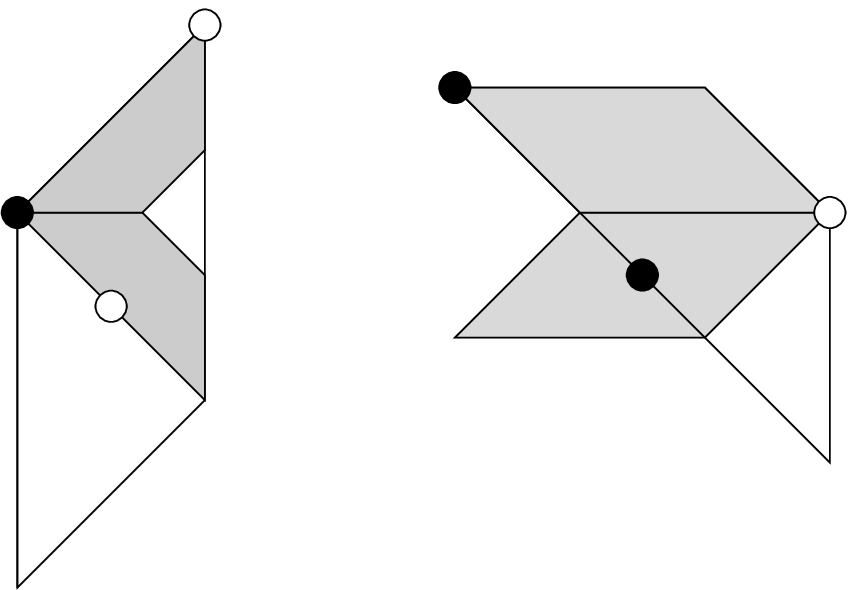}}
\newline
{\bf Figure 4.7:\/} Some special quadrilaterals.
\end{center}

The left side of Figure 4.7 shows a large parallelogram and two smaller
shaded trapezoids.   We can generate a subset of a suitably
scaled copy of the carpet starting with the large parallelogram
as a seed and then considering the portion of the resulting
limit that connects the two white dots in Figure 4.7.
Let $A$ denote the resulting set.   Likewise, we
start with each shaded trapezoid as a seed, and consider
the portion of the resulting limit that connects a black dot
to a white dot.  Let $A_1$ and $A_2$ be the resulting sets.

The right hand side of Figure 4.7 shows a trapezoid and
two shaded parallelograms.   We let $B$ and $B_1$ and $B_2$
be the portions of carpets produced by the same construction
as we did for the $A$s, but interchanging the roles of
{\it black\/} and {\it white\/}.

\begin{lemma}
\label{hiddenlemma2}
$A=A_1 \cup A_2$ and $B=B_1 \cup B_2$.
\end{lemma}

\startproof
The proof here is formally isomorphic to what we did 
for the snowflake. We omit the details.
\endproof

\newpage

\newpage

\section{The Main Result}
\label{small}

\subsection{Statement of the Result}
\label{main result}

Let $\omega=\exp(2 \pi i/8)$ be the usual primitive $8$th root of unity.
Let $P$ be the regular octagon whose vertices are powers of $\omega$.
Let $\Sigma_0$ be the pinwheel strip associated to $P$ as in \S \ref{octopix}.
Let $\Sigma_0^1 \subset \Sigma_0$ denote the half-strip that
lies to the right of the line $x=1+\sqrt 2$.   See Figure 5.1.

\begin{center}
\resizebox{!}{3.5in}{\includegraphics{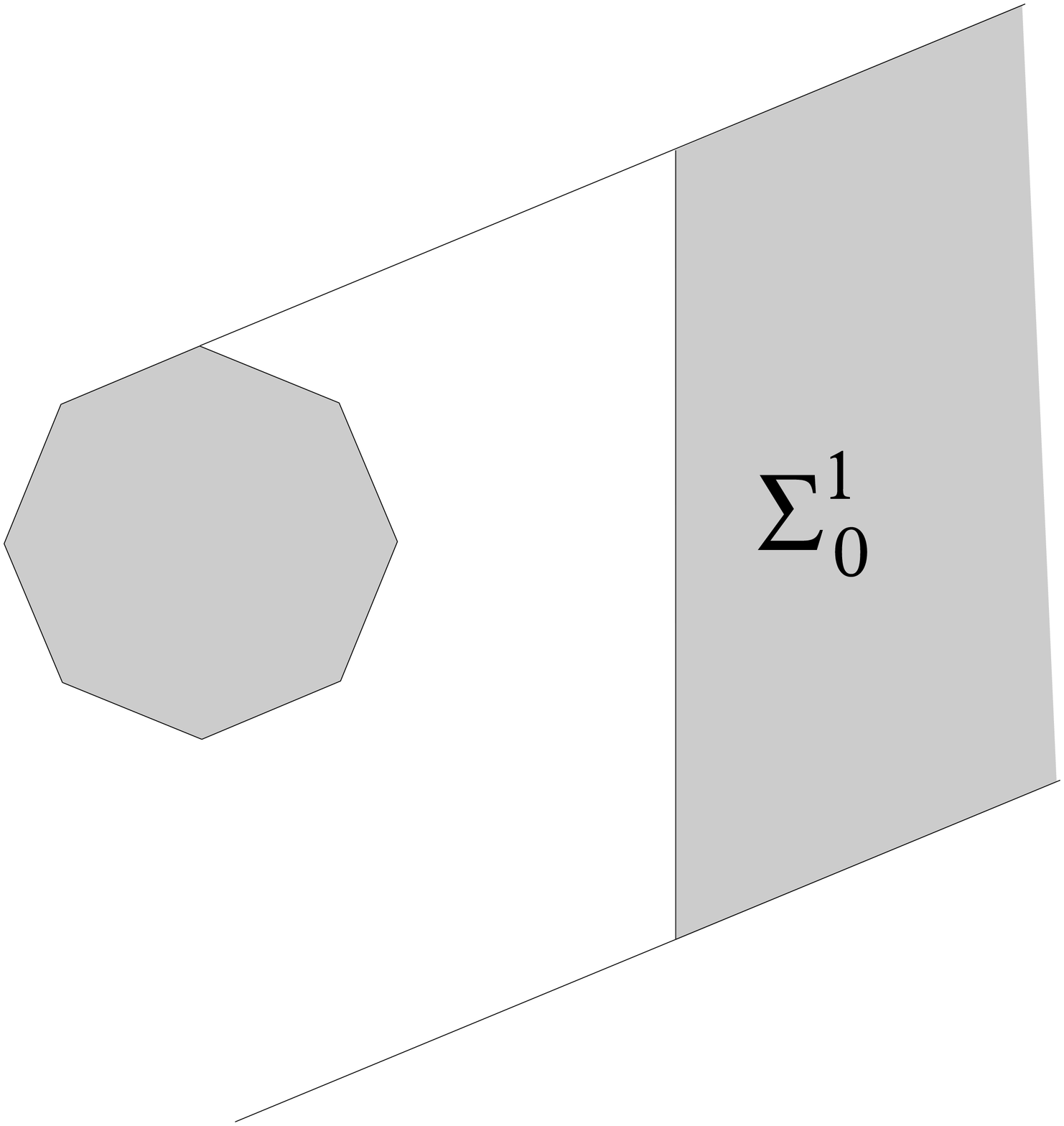}}
\newline
{\bf Figure 5.1:\/} The relevant sets
\end{center}

\S \ref{countproof} we prove

\begin{lemma}
\label{ocount}
Relative to the the pinwheel map for the regular
octagon, any periodic orbit 
that intersects $\Sigma_0^1$
nontrivially does so in exactly $3^k$ points
for some $k=0,1,2...$.
\end{lemma}

Given a periodic point $x \in \Sigma_0^1$, we define $|x|$
to be the exponent $k$ such that the pinwheel
orbit of $x$ intersects 
$\Sigma_0^1$ exactly $3^k$ times.
We call a sequence $\{x_n\} \in \Sigma_0^1$ a {\it good sequence\/} if
$|x_n| \to \infty$ as $n \to \infty$ and
he congruence of $|x_n|$ mod $4$ is odd, and independent of $n$.

Recall that $\pi_2$ and
$\pi_3$ are the projections from $\R^8$ to $\R^2$ defined
in the previous chapter.   Define
\begin{equation}
s_2=\sqrt 3; \hskip 30 pt
s_3=1+\sqrt 2.
\end{equation}
Let $\phi: I_3 \to I_2$ be the canonical surjection described above.

\begin{theorem}[Main Theorem]
\label{main}
Let $\{x_n\}$ be a good sequence. Also, let
$A_n \subset \R^8$ be the arithmetic graph of $x_n$.  
Then the following is true.
\begin{enumerate}
\item $\pi_2(A_n)$ and $\pi_3(A_n)$ are closed polygons for all $n$.
\item Let $S_{k,n}$ be dilation by
$s_k^{-|x_n|}$ about the origin.
Then the following limit exists.
$$\Gamma_k=\lim_{n \to \infty} S_{k,n} \circ  \pi_k (A_n).$$

\item $\Gamma_3=I_3$, where 
$I_3$ is scaled so that one of the $8$ isosceles triangles
in its seed has vertices 
$$
(0,0); \hskip 20 pt
(1/2,-s_3/2) \hskip 30 pt
(-s_3/2,-1/2).
$$

\item $\Gamma_2=I_2$, where $I_2$ is scaled
so that its center is 
$(3/2,3/2)$ and one of its main corners is $(-3/2,3/2)$.

\item Identifying points on $\pi_3(A_n)$ and
$\pi_2(A_n)$ if they come from the same point on $A_n$, we
get a polygonal surjectin
$$\phi_n:  S_{3,n}\pi_3(A_n) \to S_{2,n} \pi_2(A_n).$$
We have $\lim_{n \to \infty} \phi_n = F \circ \phi$.
Here $F$ is an isometry of $I_2$ and
$\phi: I_3 \to I_2$ is the canonical surjection from
the snowflake to the carpet.
\end{enumerate}
\end{theorem}

\noindent
{\bf Remarks:\/} \newline
(i)
Experimentally, we see that Theorem \ref{main} also holds when
we consider sequences of points in $\Sigma_0-\Sigma_0^1$, but
these cases lead to annoying complications.  We are interested
in establishing a nice robust version of the phenomenon, but not
necessarily the sharpest possible version. \newline
(ii) 
The limits are exactly the same when $|x_n| \equiv 1$ mod $4$ and
when $|x_n| \equiv 3$ mod $4$.    The only thing that changes
is the global isometry $F$.  In the former case, $F$ is
orientation preserving (and can be taken to be the identity)
and in the other case $F$ is orientation reversing.
\newline

The rest of the paper is devoted to proving Theorem \ref{main}.

\subsection{Restricting the Domain}
Let $P$ be the regular octagon, and let $\Sigma_0$ be the
pinwheel strip described in \S \ref{octo}.  Consider the octagons
\begin{equation}
O_k^{\pm}=P+k(\sqrt 2+1,\pm 1); \hskip 30 pt k=1,2,3 \ldots
\end{equation}

$O_k^+$ and $O_k^-$ share a common vertex and
fit inside $\Sigma_0$ as shown in Figure 5.2.   Figure 5.2
shows the cases $k=1,2$, but the remaining cases look the same.

\begin{center}
\resizebox{!}{3in}{\includegraphics{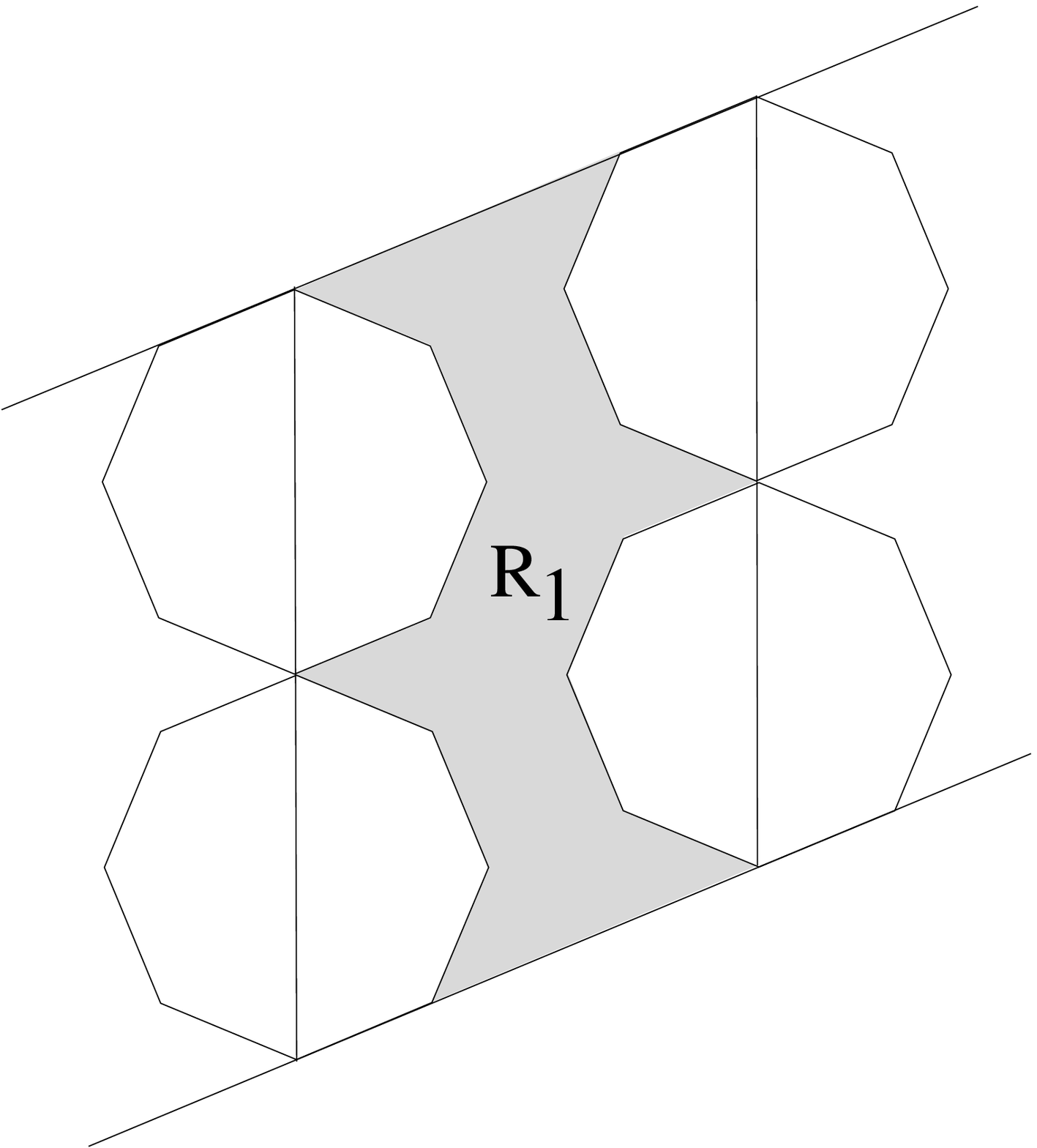}}
\newline
{\bf Figure 5.2:\/} The region $R_1$.
\end{center}

We let $R_k$ be the region of $\Sigma_0$
bounded by $O_k^+ \cup O_k^-$ and
$O_{k+1}^+ \cup O_{k+1}^-$. Figure 5.2
shows the region $R_1$, which is the main
region of interest to us.  Note that $R_k$ is translation
equivalent to $R_1$ for all $k \geq 1$,  so the
picture in Figure 5.2 is typical.

The first purpose of this chapter is to prove the
Invariance Lemma.   This result implies that each
$R_k$ is
forward invariant under the pinwheel map.  
Following this, we will prove

\begin{lemma}
\label{main2}
Each $R_k$ is an invariant set for the pinwheel map.
Furthermore, if Theorem \ref{main} is true for
sequences in $R_1$ then it is true as originally stated.
\end{lemma}

\subsection{The Pattern of Octagons}
\label{id}
\label{qr}

let $P$ be as above. The octagons we described above are
part of the infinite pattern of octagons shown in Figure 5.3.
Figure 5.3 shows a part of an infinite
set of octagons that, as it turns out, is invariant under
the outer billiards map.  Figure 5.3 also shows the
strip $\Sigma_0$.

\begin{center}
\resizebox{!}{4in}{\includegraphics{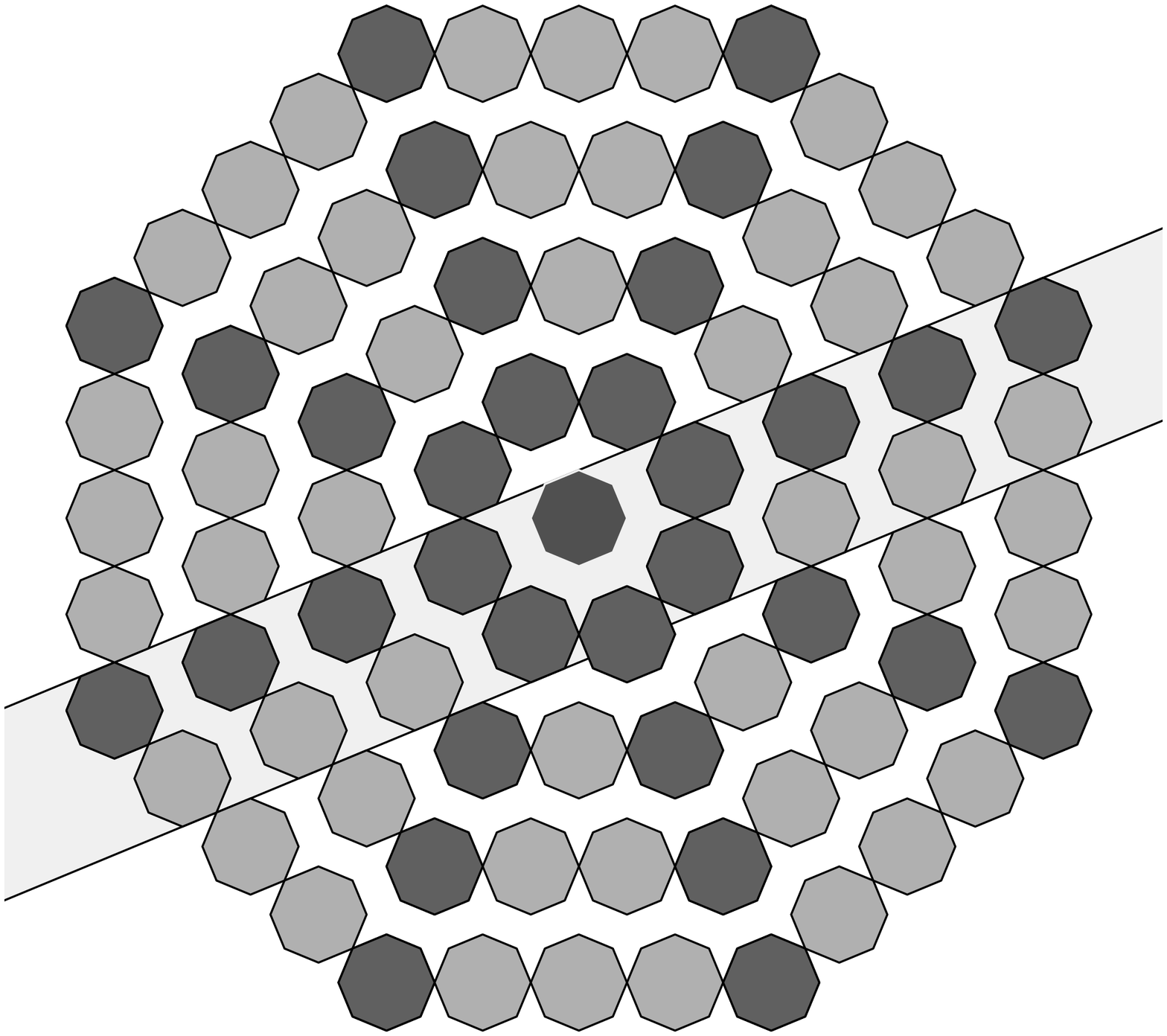}}
\newline
{\bf Figure 5.3:\/} Octagonal orbits
\end{center}

The $8$ dark rows of octagons are obtained by translating
the central one by a vector of the form
\begin{equation}
\label{translate}
k\omega^n(\sqrt 2 + 1,1); \hskip 30 pt
k,n \in \Z; \hskip 30 pt \omega=\exp(2 \pi i/8).
\end{equation}
The remaining octagons interpolate between these dark ones.
The octagons are situated in such a way that the outer billiards map
is entirely defined on the interior of each one.
The octagons come in an infinite family of rings, or
{\it necklaces\/} $N_1, N_2, N_3,...$.  The octagons
$O_k^{\pm}$ are precisely the two octagons of $\Sigma_0 \cap N_k$.
Compare  [{\bf GS\/}] and [{\bf K\/}].

\begin{lemma}
Each necklace is an outer billiards orbit.  Moreover, the
region between consecutive necklaces is invariant under the
outer billiards map.
\end{lemma}

\startproof
We change notation slightly.
Let $O(n,k;0)$ be the octagon we get by translating the central one by 
the vector in Equation \ref{translate}.
Let $O(n,k;i)$ denote the octagon that is $i$ ``clicks'' away from
$O(n,k;0)$ going counterclockwise around $N_k$. 

The outer billiard map has the same action on all the
octagons lying on the vertical line segment between
$O(-1;k;0)$ and $O(1,k;-1)$ because such octagons
all lie in the sector shown in Figure 5.4.
The outer billiards map reflects each of the octagons
under construction through the apex of the sector in
Figure 5.4.  We check easily that the action is as follows.
\begin{equation}
O(-1,k;0) \to O(4,k;1); \hskip 30 pt
O(1,k;-1) \to O(5,k;0).
\end{equation}

\begin{center}
\resizebox{!}{1.8in}{\includegraphics{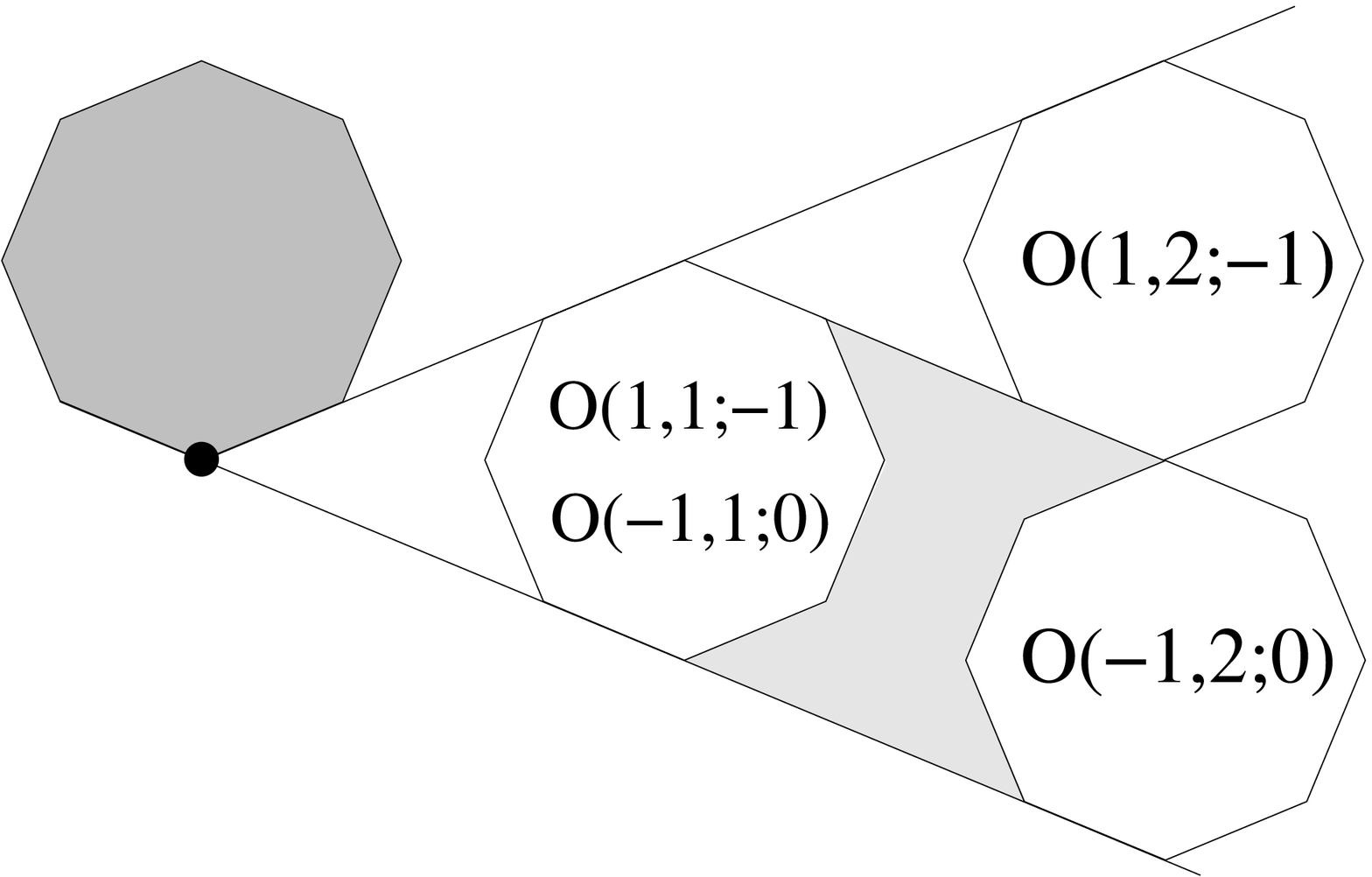}}
\newline
{\bf Figure 5.4:\/} Octagons in a sector
\end{center}

The octagons $O(4,k;1)$ and $O(5,k;0)$ have their centers
on the same vertical line segment $L$, as shown for
$k=2$ in Figure 5.4.  By symmetry, the
octagons between $O(-1,k;0)$ and $O(1,k;-1)$ map to the
octagons of the grid that have their centers on $L$.
These octagons all belong to $N_k$.
For the remaining octagons in $N_k$, the result follows
from the $8$-fold symmetry of the picture.

Call two octagons $O_1 \in N_k$ and $O_2 \in N_{k+1}$ {\it close\/}
if they are as close together as possible.   Two close octagons determine
a polygonal region ``between them'' like the one that is shaded in
Figure 5.4.  Call such a region a {\it buffer\/}.
  If the two octagons belong to the sector from Figure 5.4,
then the buffer is mapped between $N_{k}$ and $N_{k+1}$.
But then, by symmetry, all buffers are mapped between
$N_k$ and $N_{k+1}$.   But every point between $N_k$ and
$N_{k+1}$ is contained in a buffer.
\endproof

\subsection{Invariant Domains in the Strip}

Let $\Sigma'_0$ denote those points in $\Sigma_0$
that lie to the central octagon $P$.   A point
$p$ lies in $\Sigma_0'$ if every vector pointing
from $P$ to $p$ has positive $x$ coordinate.
Let $\Phi': \Sigma_0' \to \Sigma_0'$ be the
first return map of the $\phi^2$, the square
outer billiards map.

\begin{corollary}
\label{neck1}
For $k \geq 2$ the set $R_k$ is invariant under $\Phi'$.
\end{corollary}

\startproof
For $k \geq 2$, the region between $N_k$ and $N_{k+1}$
intersects $\Sigma_0'$ exactly in the set $R_k$.
\endproof

\noindent
{\bf Remark:\/}
The region between $N_1$ and $N_2$ intersects $\Sigma_0$ in a
set that is somewhat larger than $R_1$., as one can see by
looking carefully at Figure 5.3. It is for this reason
that Corollary \ref{neck1} fails for $R_1$.
\newline

\begin{corollary}
For $k \geq 4$ the set $R_k$ is invariant under $\Phi$,
the pinwheel map.
\end{corollary}

\startproof
For $k \geq 4$ the argument in Lemma \ref{far2} goes through,
and shows that $\Phi=\Phi'$.
\endproof

Below we will prove the following result.

\begin{lemma}
\label{ag1}
Let $\mu=1,2,3...$
Let $p \in R_1$ and 
$$q=p+2\mu (\sqrt 2+1,1) \in R_{1+\lambda}.$$
Then the arithmetic graphs $\Gamma(p)$ and $\Gamma(q)$
coincide.
\end{lemma}

\begin{corollary}
\label{neck2}
The set $R_k$ is $\Phi$-invariant for all $k \geq 1$.
\end{corollary}

\startproof
By Lemma \ref{ag1}, translation by the
vector $4(\sqrt 2+1,1)$ conjugates the action
of $\Phi$ on $R_1$ to the action of $\Phi$ on
$R_5$.  We already know that $R_5$ is
invariant under $\Phi$ and the conjugation
property implies the same result for $R_1$.
A similar argument takes care of $R_2$ and $R_3$.
We have already handled the remaining cases.
\endproof

\noindent
{\bf Remark:\/}
Note that $\Phi$ and $\Phi'$ do not coincide on
$R_1$, because Corollary \ref{neck1} is false for
$R_1$ whereas Corollary \ref{neck2} is true.
These maps in fact coincide for $k=2,3,4...$.
\newline

\noindent
{\bf Proof of the Invariance Lemma:\/}
We want to show that
$\Sigma_0^k$ is $\Phi$ invariant.
Let $N_k^+$ be the set of half-octagons
obtained by taking the right halves of all the
octagons in $N_k$.   Since $\phi^2$ is
a piecewise translation, the set $N_k^+$
is invariant under $\phi^2$.  But then
the set
$$S_k=N_k^+ \cap \Sigma_0$$ is $\Phi'$ invariant
for $k \geq 4$.   Hence, by the argument in
Lemma \ref{far2}, the set $S_k$
is $\Phi$-invariant for $k \geq 4$.
As in the proof of Corollary \ref{neck2}, we now conclude
that $S_k$ is $\Phi$ invariant for all $k=1,2,3...$
Finally
$$\Sigma_0^k=S_k \cup R_{k+1} \cup S_{k+1} \cup R_{k+2} \ldots.$$
We have decomposed $\Sigma_0^k$ into
$\Phi$ invariant sets.  Hence $\Sigma_0^k$ is also
$\Phi$ invariant.
\endproof

\noindent
{\bf Proof of Lemma \ref{main2}:\/}
Suppose we want to prove Theorem \ref{main} for
a sequence of periodic points $\{x_n\}$ in $\Sigma_0^1$.
Since $\Phi$ is the identity on each of the orbits
$N_k$, we have $x_n \in R_{k_n}$ for large $n$.
Chopping off the initial portion of the sequence,
we can assume that $x_n \in R_{k_n}$ for all $n$.

By Lemma \ref{ag1}, we can find a point
$y_n \in R_1 \cup R_2$ such that $x_n$ and $y_n$
have the same arithmetic graph.
Hence, it suffices to prove Theorem \ref{main}
for the sequence $\{y_n\} \in R_1 \cup R_2$.

The intersection $R_1 \cup R_2$ is a single point,
namely $\zeta=(2 + 2 \sqrt 2,1)$.    We call
two points {\it specially related\/} if
reflection through $\zeta$ interchanges the two
points.    We check by hand the following property,
which we call {\it Property P\/}:
If $y$ and $z$ are specially related then
$\Phi(y)$ and $\Phi(z)$ are specially related and
$\widehat \Phi(y)+\widehat \Phi(z)=0$.

From Property P we see that the arithmetic graphs of
$y$ and $z$ are isometric to each other.
Using this basic principle, we can replace each
$y_n \in R_2$ by the specially related point
$z_n$ without changing the isometry type 
of the rescaled limits of the arithmetic graphs.
\endproof

\noindent
{\bf Remark:\/}  In \S \ref{pinwheeldyn}
we will explain in detail how
we actually check Property P.

\subsection{Equivalent Points in Strips}
\label{stripmap}

The rest of the chapter is devoted to the proof of
Lemma \ref{ag1}.
In this section we make a general
observation about strip maps.
In the section following this one, we will
apply our observation to the case of the regular octagon.

Let $\Sigma_1$ and $\Sigma_2$ be non-parallel strips.
Let $V$ and $W_1,W_2$ be the three vectors shown in Figure 5.3.
We say that $p,q \in \Sigma_1$ are {\it related\/}
if $p-q=kW_1$ for some integer $k$.  Likewise, we say that
$p,q \in \Sigma_2$ are {\it related\/} if $p-q=kW_2$
for some integer $k$.

\begin{center}
\resizebox{!}{3in}{\includegraphics{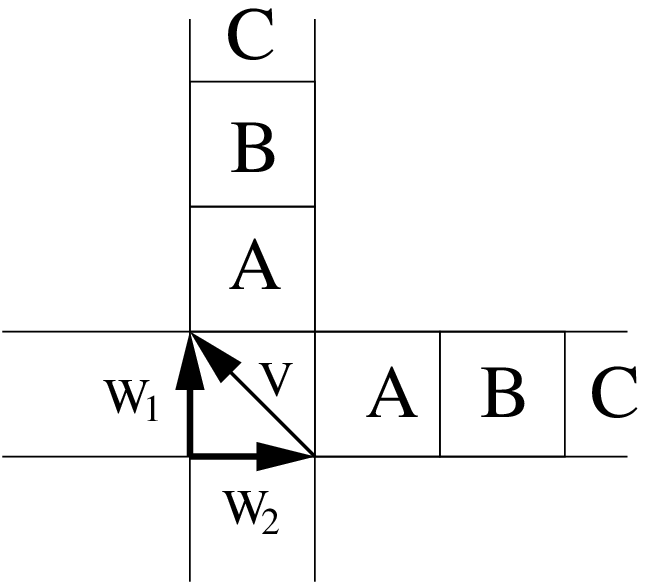}}
\newline
{\bf Figure 5.4:\/} A Strip map
\end{center}

\begin{lemma}
\label{propagate}
Let $f$ be the strip map associated to $(\Sigma_2,V)$.  If
$p,q \in \Sigma_1$ are related then $f(p),f(q) \in \Sigma_2$
are related.
\end{lemma}

\startproof
Our statement is an affinely invariant one.  So, we normalize
our strips so that they are perpendicular and each have width $1$.
In this case, $f$ translates the horizontal squares $A,B,C,...$ to
the vertical squares with the same labels.  If $p,q \in \Sigma_1$
are related then they are in the same positions relative to the
horizontal squares that contain them.  But then $f(p),f(q) \in \Sigma_2$
are in the same position relative to the vertical squares that
contain them.  Hence $f(p)$ and $f(q)$ are related.
\endproof

\subsection{Proof of Lemma \ref{ag1}}

We first prove Lemma \ref{ag1} in case that
$\lambda=2 \mu$ is even.

We define
\begin{equation}
W_{2n}=W_{2n+1}=2\omega^n (\sqrt 2 + 1,1).
\end{equation}
The vector $W_0=W_1$ is parallel to both
pointed strips $\Sigma_0$ and $\Sigma_1$.  The
vector $W_2=W_3$ is parallel to both
$\Sigma_2$ and $\Sigma_3$, and so on.

We say that two points $p,q \in \Sigma_{k-1}$ are {\it equivalent\/}
\begin{equation}
\label{equiv}
q-p=\mu W_{k-1}; \hskip 30 pt \mu \in \Z.
\end{equation}
We let $\mu_{k-1}(p,q)$ be the integer in Equation \ref{equiv}.

Let $f_k$ be the strip map associated to
$(\Sigma_k,V_k)$, as in \S \ref{ag0}.
Let $m_k(p)$ be as in Equation \ref{gr2}.  That is,
\begin{equation}
f_{k}(p)=p+m_{k}(p) V_k.
\end{equation}

\begin{lemma}
\label{odd}
Let $f=f_k$.
If $p,q \in \Sigma_{k-1}$ are equivalent and
$k$ is odd, then
$f(q),f(q)$ are equivalent in $\Sigma_k$ and
$$\mu_{k}(f(p),f(q))=\mu_{k-1}(p,q);
\hskip 30 pt m_k(p)=m_k(q).$$
\end{lemma}

\startproof
By rotational symmetry, it suffices to consider the
case when $k=1$.
Suppose that $k=1$.     Note that $\Sigma_0 \cap \Sigma_1$ is
a strip parallel to $W_0=W_1$. So, $p \in \Sigma_0 \cap \Sigma_1$
if and only if $q \in \Sigma_0 \cap \Sigma_1$.     
Suppose that $p,q \in \Sigma_0 \cap \Sigma_1$.
Then $f(p)=p$ and $f(q)=q$.  Also, $W_1=W_0$.
Hence $f(p)$ and $f(q)$ are equivalent in
$\Sigma_1$ and  
$$\mu_1(f(p),f(q))=\mu_0(p,q); \hskip 30 pt
f_1(p)=f_1(q)=0.$$
If $p,q \in \Sigma_0-\Sigma_1$ then
$p+V_1 \in \Sigma_1$ and likewise $q+V_1 \in \Sigma_1$.
Here $$f_1(q)-f_1(p)=q-p.$$  Hence
$$\mu_1(f(p),f(q))=\mu_0(p,q); \hskip 30 pt
f_1(p)=f_1(q)=1.$$
\endproof

\begin{lemma}
\label{even}
Let $f=f_k$.
If $p,q \in \Sigma_{k-1}$ are equivalent and
$k$ is even, then
$f(q),f(q)$ are equivalent in $\Sigma_k$ and
$$\mu_{k}(f(p),f(q))=\mu_{k-1}(p,q);
\hskip 30 pt m_k(q)=m_k(p)+\mu_{k-1}(p,q).$$
\end{lemma}

\startproof
By rotational symmetry, it suffices to consider the case $k=2$.
In this case, the strips $\Sigma_1$ and $\Sigma_2$, as well
as the vectors
$W_1$ and $W_2$ and $V=V_2$ are related precisely as
discussed in \S \ref{stripmap}.   Lemma
\ref{propagate} now implies that $p$ and $q$ are
equivalent in $\Sigma_2$ and
$\mu_2(p,q)=\mu_{1}(p,q)$.

Referring to Figure 5.4, which is an affine image
of the octagon picture, the number $\mu_1(p,q)$
counts the number of parallelograms one needs to
shift in order to move $p$ over to $q$.  For
instance, if $p$ were in parallelogram $A$ and
$q$ were in parallelogram $C$ in Figure 5.4, then
$\mu_0(p,q)=2$.    From Figure 5.4, we see
clearly that  $m_2(q)=m_k(p)+\mu_{1}(p,q).$
\endproof

Now we can finish the proof of Lemma \ref{ag1}.
Suppose that $p,q \in \Sigma_0$ are equivalent.
Let $\mu=\mu_0(p,q)$.  Then, referring
to Equation \ref{spectrum}, we have
\begin{equation}
\widehat \Phi(q) - \widehat \Phi(p)=
\sum_{i=1}^8 (m_k(q)-m_k(p)) \widetilde V_k=^1
\sum_{k\  {\rm even\/}} \mu \widetilde V_k=^20.
\end{equation}
\label{cancel}
Equality 1 comes from Lemmas \ref{odd} and \ref{even}, and
Equality 2 comes from the fact that
$\widetilde V_6-=\widetilde V_2$
and $\widetilde V_8=-\widetilde V_4$.

\newpage
\section{A Toy Model}

\subsection{The Map}

The work in the previous chapter reduces the proof of Theorem \ref{main}
to the study of the pinwheel map $\Phi: R_1 \to R_1$, where
 $R_1$ is  shown in Figure 5.2.
This first return map turns out to be a polygon exchange map.
In this chapter, we will briefly discuss a simpler but related polygon
exchange map.   

Though the simpler map is not the one we will
ultimately study in our proof of Theorem \ref{main}, we think
that it will prepare the reader for the more complicated system
we do study.    First,  the basic shapes that arise in the
simple system here arise in the more complicated system we
study later.  Second, both the system here and the one we
do study have a {\it renormalization scheme\/}.  It is the
renormalization scheme for the first return map
that leads to the self-similar nature
of the arithmetic graphs.

The dynamical system we study here is well known,
and indeed also arises in the study of outer billiards:
It is the first return map to a certain
 invariant domain for a dynamical system generated
by the outer billiards map and the rotational symmetry
group of the regular octagon.
Compare [{\bf T2\/}] and [{\bf BC\/}].
The system is a self map of the kite-shaped region $X$, shown as the
shaded region in Figure 6.1.  The region $X$ is best defined in terms
of a large regular octagon.    The smaller regular octagon, which
is a subset of $X$, is drawn in
for later reference.

\begin{center}
\resizebox{!}{1.6in}{\includegraphics{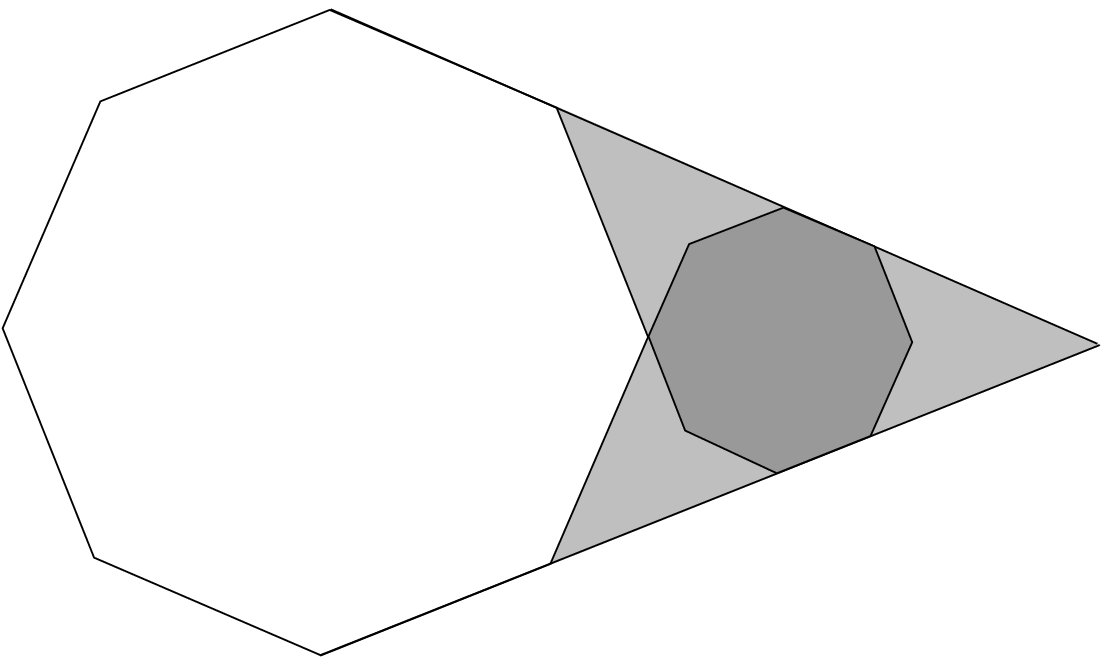}}
\newline
{\bf Figure 6.1:\/} The region $X$.
\end{center}

As shown in Figure 6.2, the
 region $X$ can be partitioned in two ways into right angled
isosceles triangles,
\begin{equation}
X=A_1 \cup B_1 = A_2 \cup B_2.
\end{equation}

There is a unique orientation-preserving isometry $\phi_A: A_1 \to A_2$.
The map $\phi_A$ rotates $3 \pi/4$ radians clockwise about the center
of the small octagon in Figure 6.1.  Likewise, there is a unique
orientation preserving isometry $\phi_B: B_1 \to B_2$.  The
map $\phi_B$ rotates $\pi/4$ radians counterclockwise about
the center of the big octagon in Figure 6.1.
  We have a map
$\phi: X \to X$ which restricts to $\phi_A$ on $A_1$ and $\phi_B$
on $B_1$.  The map $\phi$ is defined except on the segment
common to $A_1$ and $B_1$.

\begin{center}
\resizebox{!}{1.5in}{\includegraphics{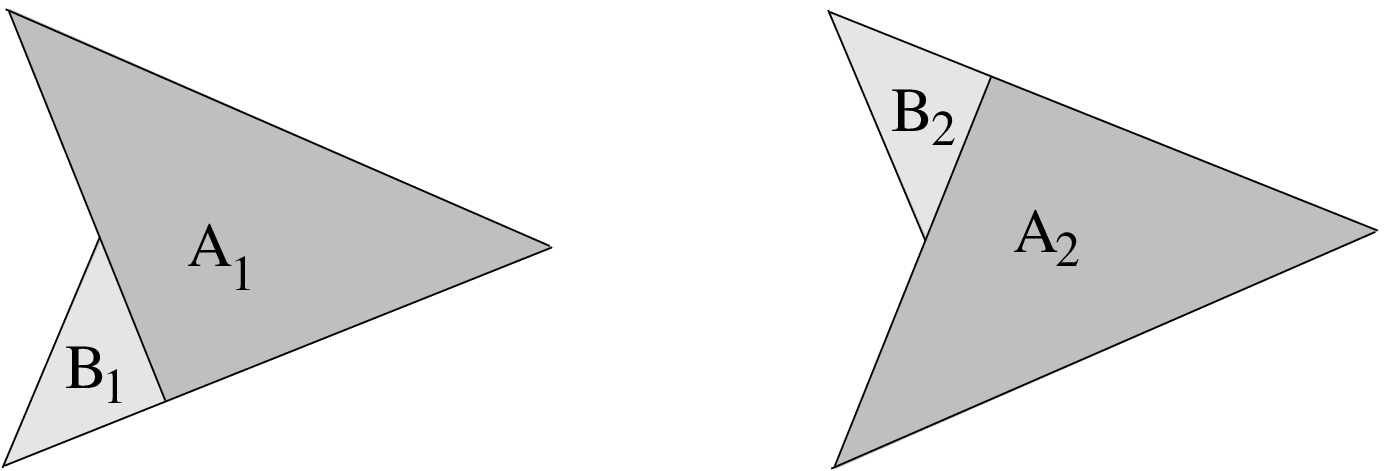}}
\newline
{\bf Figure 6.2:\/} Two partitions of $X$.
\end{center}

\subsection{The Renormalization Scheme}

Let $O \subset X$ denote the inner octagon in Figure 6.1.
By construction, the center of $O$ is a fixed point of
$\phi$, and every point of $O$ has period $8$ with
respect to $\phi$.

\begin{center}
\resizebox{!}{2.8in}{\includegraphics{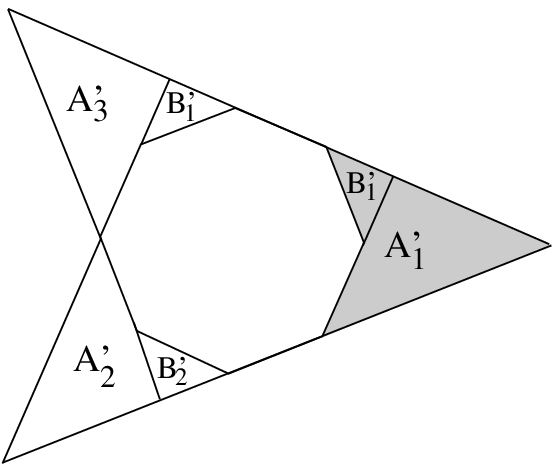}}
\newline
{\bf Figure 6.3:\/} Renormalization Scheme
\end{center}

Let $X'$ be the shaded component of $X-O$ shown in Figure 6.3.
We have partitioned $X'$ into two pieces, which we call
$A_1'$ and $B_1'$.    The third power $\phi^3$ maps
$X'$ to itself.  Figure 6.2 shows the pieces
$A_2'=\phi(A_1')$ and $A_3'=\phi(A_2')$ and
$B_2'=\phi(B_1')$ and $B_3'=\phi(B_2')$.

We have not drawn $A_4', B_4' \subset X'$, but
we observe that $A_4'$ and $B_4'$ relate to
$A_1'$ and $B_1'$ in the same say that
$A_2$ and $B_2$ relate to $A_1$ and $B_1$.
More precisely, there is an orientation reversing
dilation $\Theta: X' \to X$ such that
\begin{equation}
\label{renorm1}
\Theta^{-1} \circ \phi \circ \Theta=\phi^3
\end{equation}
on $X'$.
The map $\Theta$ is the {\it remormalization map\/}.
Having $\Theta$ by itself is nice, but we have a second
property that really pins things down.  Namely,
\begin{equation}
\label{renorm2}
X-O=X' \cup \phi(X') \cup \phi^2(X').
\end{equation}
In other words, every point of $X-O$ lies in the
orbit of a point in $X'$.

\subsection{The Periodic Points}

Our remormalization scheme allows us to get a complete
understanding of the periodic points of $\phi$.
The octagon $O$ consists of points having period $8$.
(The center is the one point of $O$ that has period $1$.)
For our purposes, it is nicer to think of $O$ as having
period $1$ {\it as a tile\/}.  The set $O$ is
preserved by $\phi$.  We set $O_0=O$.

By equation \ref{renorm1}, the region $X'$ contains a small
octagon 
\begin{equation}
O'_0=\Theta^{-1}(O_0)
\end{equation}
that has tile-period $3$.
We set
\begin{equation}
O'_{k+1}=\phi(O'_k).
\end{equation}
The tiles $O'_k$ for $k=0,1,2$ all have tile-period $3$.
They comprise a tile-orbit.

Now we iterate.  By Equation \ref{renorm1}, the three octagons
\begin{equation}
O''_{3k}=\Theta^{-1}(O'_k); \hskip 30 pt k=0,1,2
\end{equation}
all have tile-period $9$. 
Setting \begin{equation}
O''_{k+1}=\phi(O''_k),
\end{equation}
we produce $9$ small octagons
$O''_0,...,O''_8$ that all have tile-period $9$.
Next, we produce $27$ small octagons
$O^{(3)}_k$ for $k=0,...,26$, all having
tile-period $27$.   Figure 6.4 shows some of these
octagons.

\begin{center}
\resizebox{!}{3.5in}{\includegraphics{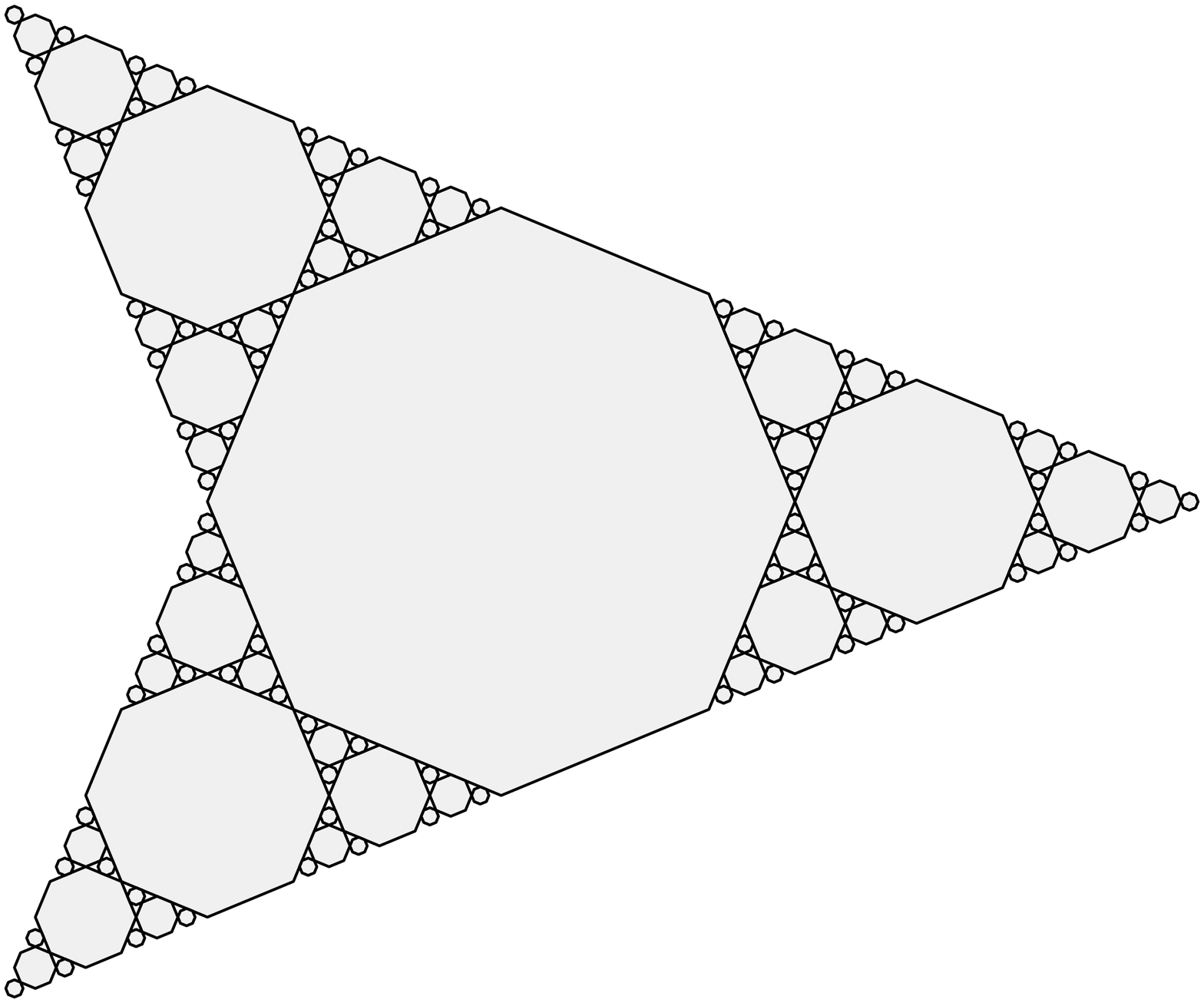}}
\newline
{\bf Figure 6.4:\/} Cascade of periodic tiles
\end{center}

Thus, we have produced, for each $n=0,1,2,3...$, a periodic 
octagonal tile having tile-period $3^n$.   We claim that
every periodic point is contained in one of these octagons.
To see this, we use Equation \ref{renorm2}.
Clearly, the only periodic points of order $1$ lie in $O_0$.
Suppose $N$ is the smallest positive integer for which
we have not proven our result.  Let $y'$ be a periodic
point of period $N$.

Since $N>1$, we have $y' \in X-O$.  By Equation \ref{renorm2}, the point
$y'$ lies in the orbit of some point $x' \in X'$.   Let
$x=\Theta(x')$.  Note that
$\phi(X')$ and $\phi^2(X')$ are disjoint from $X'$.
Since $\phi^N(x')=x'$, we must therefore have $N=3M$ for
some integer $M$.   But then, by Equation \ref{renorm1},
we see that $x$ has period $M$.     By induction, $x$ lies
in one of our octagonal tiles.  But then, by construction,
so does $x'$.   Now we know that $x'$ is one of the
periodic points we have already studied.  So, the
same goes for $y'$, the original point of interest to us.

Using the renormalization scheme, we have been able to
classify all the periodic points of the system.  Incidentally,
we note that the complement of the octagons has dimension
larger than $1$, and the set of points where some iterate of
$\phi$ is undefined has dimension $1$.  Therefore,
$\phi$ has some aperiodic points.  This is the same argument
Tabachnikov [{\bf T2\/}] gives for the regular pentagon.

\newpage

\section{The Pinwheel Dynamics}
\label{pinwheeldyn}

\subsection{The Partition}

We will work with the region $R_1$ from Lemma \ref{main2}.
Both sides of Figure 7.1 show $R_1$.   Note that $R_1$ is the union
of $4$ copies of the region $X$ from the previous chapter.

\begin{center}
\resizebox{!}{4.5in}{\includegraphics{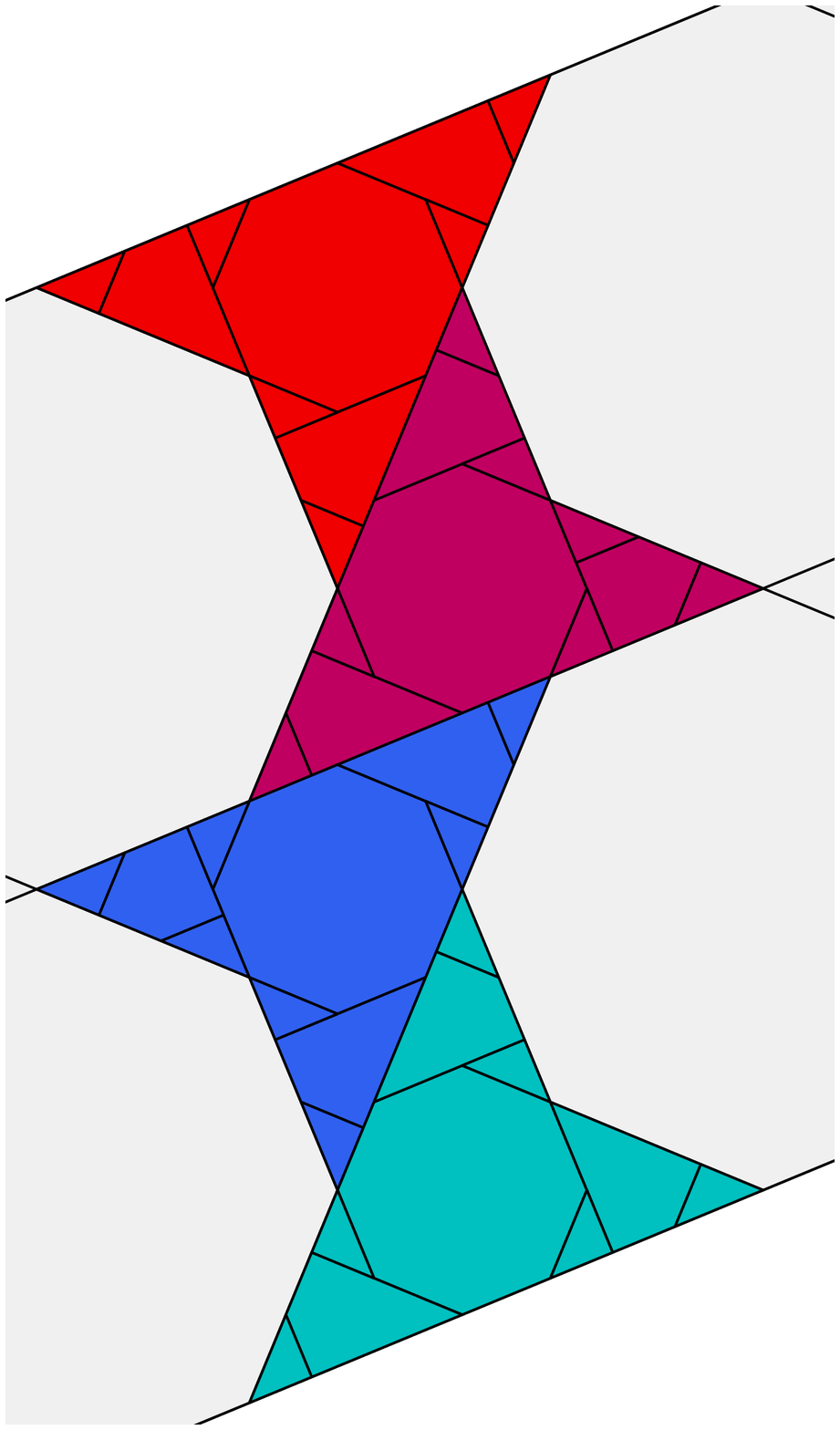}}
\resizebox{!}{4.5in}{\includegraphics{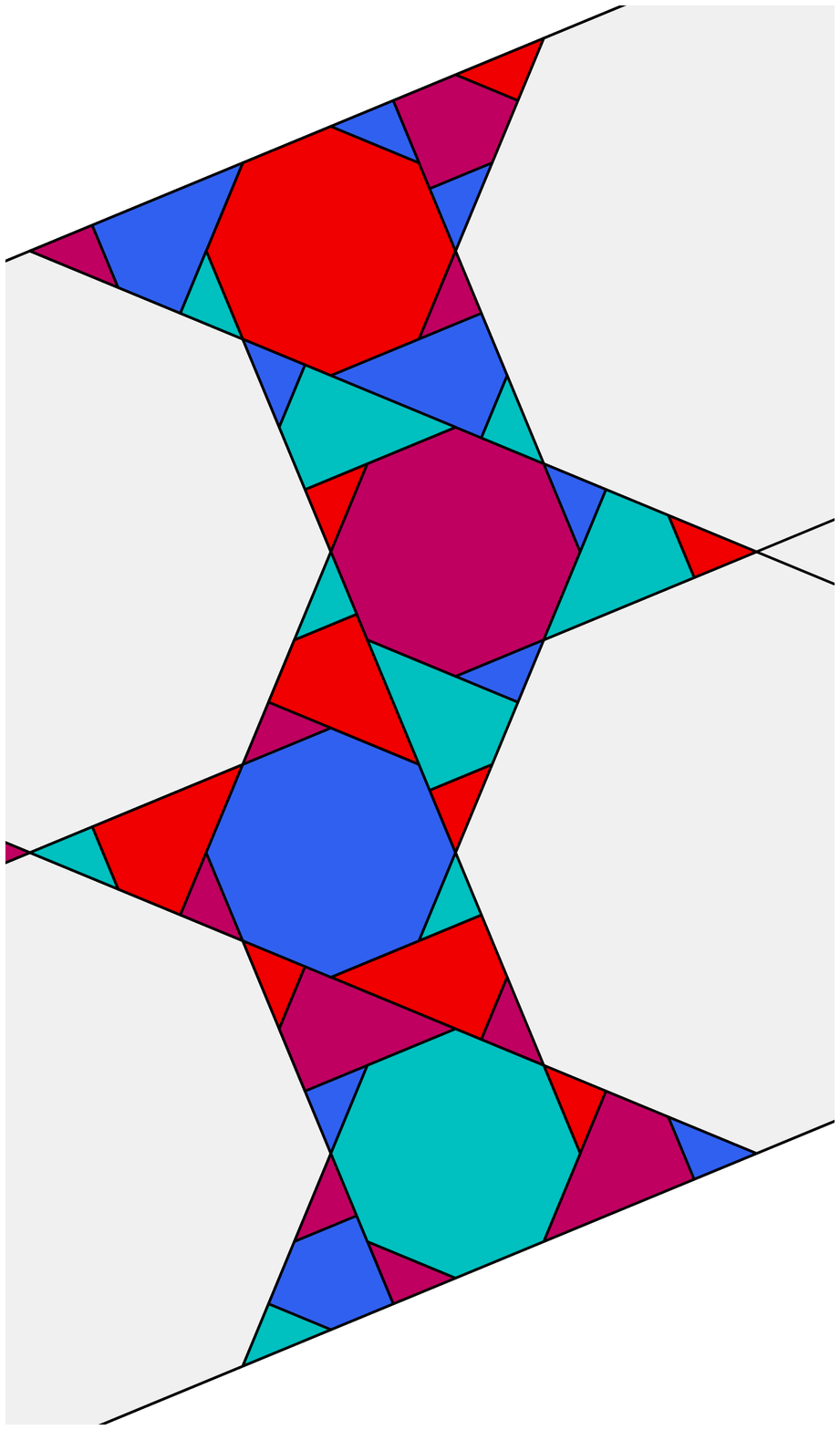}}
\newline
{\bf Figure 7.1:\/} The region $R_1$, inside the strip $\Sigma_0$.
\end{center}

Figure 7.1 shows the action of $\Phi$ on $R_1$.  The left
hand side of the figure shows the partition of $R_1$ into
the maximal pieces on which $\Phi$ and $\widehat \Phi$
are constant, and the right hand side shows the images of
these pieces under $\Phi$.  We will give coordinates for
these tiles below.

Just knowing that $\Phi$ is a piecewise translation, we can
nearly pin down the action of $\Phi$ from Figure 7.1.
The problem is that there are some pairs of triangles,
having identical shading, that are translates of each other.
Such pairs are specially situated:  One triangle in each such
pair lies on the left hand side of $R_1$ and the other triangle
lies in the right hand side.   To pin down the action
of $\Phi$ exactly, we mention 
that $\Phi$ maps any triangle on the left (respectively right) hand side of
$R_1$ to a triangle on the right (respectively left) hand side of $R_1$.
\newline
\newline
{\bf Remark:\/}
It seems worth also mentioning that
$\Phi$ maps any quadrilateral or pentagon on the
left (respectively right) hand side of $R_1$ to another quadrilateral or
pentagon on the left (respectively right) hand side of $R_1$.  
Finally, $\Phi$ fixes each of the central octagons.
\newline
\newline
$\Phi$ has one obvious symmetry and one hidden symmetry.
Let $\rho_1$ be the order-$2$ rotation about the center point
of $R_1$.   Then $\Phi$ commutes with $\rho$.

To describe the hidden symmetry, we observe that it is better 
to think of $R_1$ as a subset of an infinite cylinder.  
Define $\sigma(x,y)=(x,y+2)$.
The squared-map $\sigma^2$ maps the bottom boundary of
$\Sigma_0$ to the top one.  The quotient, obtained
by identifying the two boundary components, is a cylinder.
We think of $R_1$ as a subset of this cylinder, so that certain
pieces on the top of $R_1$ are contiguous with pieces on
the bottom.

This point of view makes the picture look more symmetric.
Consider the magenta pentagon $K'$ shown on the middle
left hand side of Figure 7.1.   Let $K''$ be the magenta triangle such that
$K=K' \cup K''$ is a kite, isometric to the rest of the kites in the picture.
$\Phi(K')$ lies at the very top of
$R_1$ and $\Phi(K'')$ lies at the very bottom.  Under
the identification we have been discussing,
$\Phi(K'')$ is just the continuation of $\Phi(K')$.  So, when we
think of $R_1$ as a subset of a cylinder, $K$ is a single tile
on which $\Phi$ is constant.
The other pentagon/triangle pair
on the left hand side of Figure 7.1 has the same analysis.
Considered this way, we see that the partition of $R_1$ consists
just of kites, triangles, and octagons.  

Finally, we mention the hidden symmetry.   Once
we interpret $R_1$ as being a subset of the cylinder,
$\Phi$ commutes with $\sigma$.
\newline
\newline
\noindent
{\bf Remark:\/} It turns out that
$\widehat \Phi$ does not share the same symmetries.
That is, it is not true that the maps
$\widehat \Phi$ and $\widehat \Phi \circ \rho$ and
$\widehat \Phi \circ \sigma$ coincide.  However,
as we will see, there is some relation between these maps.

\subsection{Notation and Coordinates}

This section does not contribute to the mathematics at all. We include
it for  readers who would like to see explicit coordinates to that they
can reproduce our experiments and
calculations on their own.
Let
\begin{equation}
s=\sqrt 2 +1.
\end{equation}
 Let $O_3$ be the red octagon in
Figure 7.1.   The octagon $O_3$ has radius $1/s$.  By
this we mean that the distance from the center of $O_3$
to any vertex is $1/s$.  We point out
another octagon, $O_4$, of radius $1/s^2$.
The octagon $O_4$ is homothetic to $O_3$, and
the rightmost vertex of $O_4$ coincides with the
leftmost vertex of $O_3$.   Both $O_3$ and $O_4$ are
periodic tiles for $\Phi$.  See Figure 7.2 below.
Now we distinguish $4$ special points.

\begin{itemize}
\item The point 
$c_1=(3s,1)/2$ is the center of $R_1$.
\item The point 
$c_2=(3s,3)/2$ 
is the center of the top half of $R_1$.
\item The point $c_3=(2+\sqrt 2,2)$ is the center of $O_3$.
\item The point $c_4=(2\sqrt 2,2)$
is the center of $O_4$.
\end{itemize}
Let $\rho_k$ denote counterclockwise rotation by
$2 \pi k/8$ about $c_k$.    

\begin{center}
\resizebox{!}{2.8in}{\includegraphics{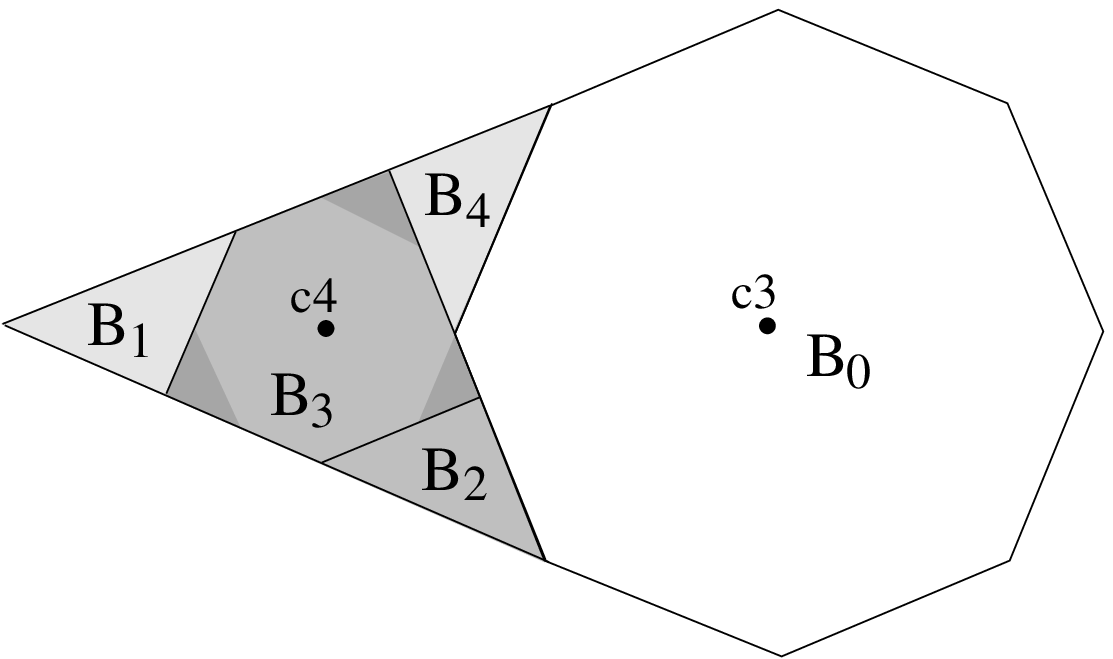}}
\newline
{\bf Figure 7.2:\/} Some polygons of interest.
\end{center}

Figure 7.2 shows $4$ shaded polygons.  
For these polygons, we will use a second labelling scheme
in which $O_4=B_0$.   Figure 7.2 shows the $B$ labels.
 The union $B_2 \cup B_3$ is
the kite that lies to the left of the red octagon in Figure 7.1.
We have tried to shade $B_3$ in such a way that we point out
$O_4 \subset B_3$.    

The vertices of $2B_3$ are
\begin{equation}
(5s-6,2s-1) \hskip 8 pt
(4s-4,4s-6) \hskip 8 pt
(3s-2,2s-1) \hskip 8 pt
(s+3,3s-3) \hskip 8 pt
(2s+1,s+2).
\end{equation}
We divide these coordinates in half to get the coordinates of $B_3$.

The vertices of $2B_1$ are given by
\begin{equation}
(2s,4); \hskip 30 pt
(3+s,3s-3); \hskip 30 pt
(3s-2,2s-1).
\end{equation}
We divide these coordinates in half to get the coordinates of $B_3$.

We can deduce the remaining coordinates just from what we have
already given.  For instance,
\begin{equation}
B_4=\rho_4^5(B_1); \hskip 30 pt
B_2=\rho_4^3(B_1).
\end{equation}
The remaining polygons in the partition of $R_1$ all have
the form
\begin{equation}
\rho_1^{a_1} \circ \rho_2^{a_2} \circ \rho_3^{a^3}(B);
\hskip 30 pt
a_1,a_2 \in \{0,1\}; \hskip 30 pt a_3 \in \{3,5\}.
\end{equation}
Here either $B=B_j$ or $B=B_2 \cup B_3$.

Now we explain briefly how we verify that the return
map $\Phi$ is as stated.   We explain this for $B_3$.
The check is the same for the remaining files.
Recall that $\Phi=f_{16} \circ ... \circ f_1$,
where $f_k$ is the strip map associated to the
pair $(\Sigma_k,V_k)$.    We verify that there are
integers $m_1,...,m_{16}$ such that
\begin{itemize}
\item $B_3^1:=B_3+m_1V_1 \in \Sigma_1$.
\item $B_3^2:=B_3^1 + m_2 V_2 \in \Sigma_2.$
\item $B_3^3:=B_3^2 + m_3 V_3 \in \Sigma_3.$
\end{itemize}
And so on.  By convexity we just have to check this
on the vertices of $B_3$.  This is a short calculation
we omit.  This calculation shows not only that
$\Phi(B_3)=B_3^{16}$, but also that
$\widehat \Phi$ is constant on the interior of $B_3$.

\subsection{The Compressed System}
\label{compress}

Let $R \subset R_1$ denote the top half of $R_1$.
The region $R$ is colored red and magenta on the
left half of Figure 7.1.    Note that $R$ is a
fundamental domain for the action of
$\sigma$ on $R_1$.    We define a new dynamical
system $\Psi: R \to R$ as follows:
$\Psi(p)=\Phi(p)$ if $\Phi(p) \in R$ and
$\Psi(p)=\sigma \circ \Phi(p)$ if
$\Phi(p) \not \in R$.  In other words,
we use the action of $\sigma$ to move everything
into $R$.

\begin{center}
\resizebox{!}{2.6in}{\includegraphics{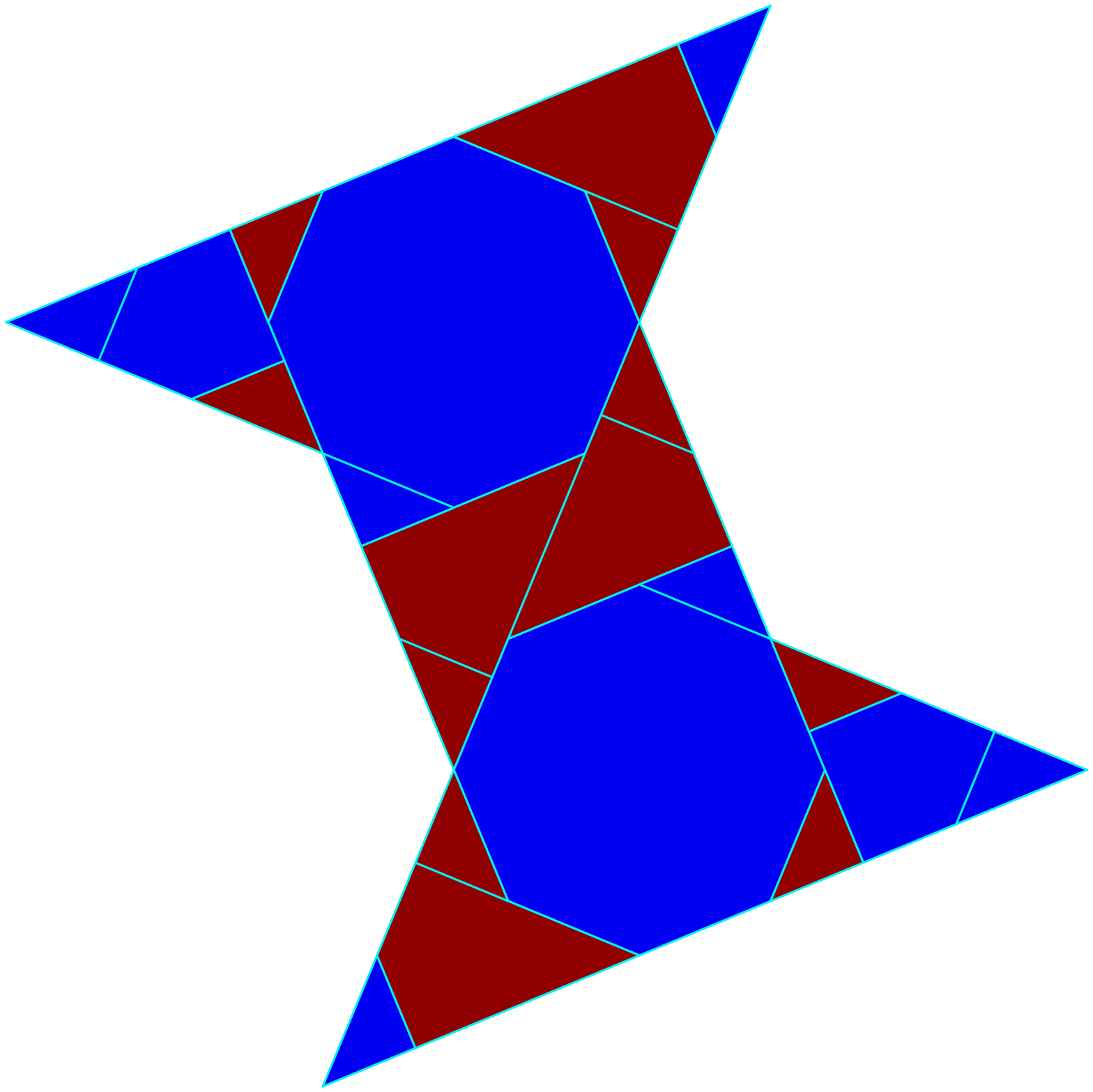}}
\resizebox{!}{2.6in}{\includegraphics{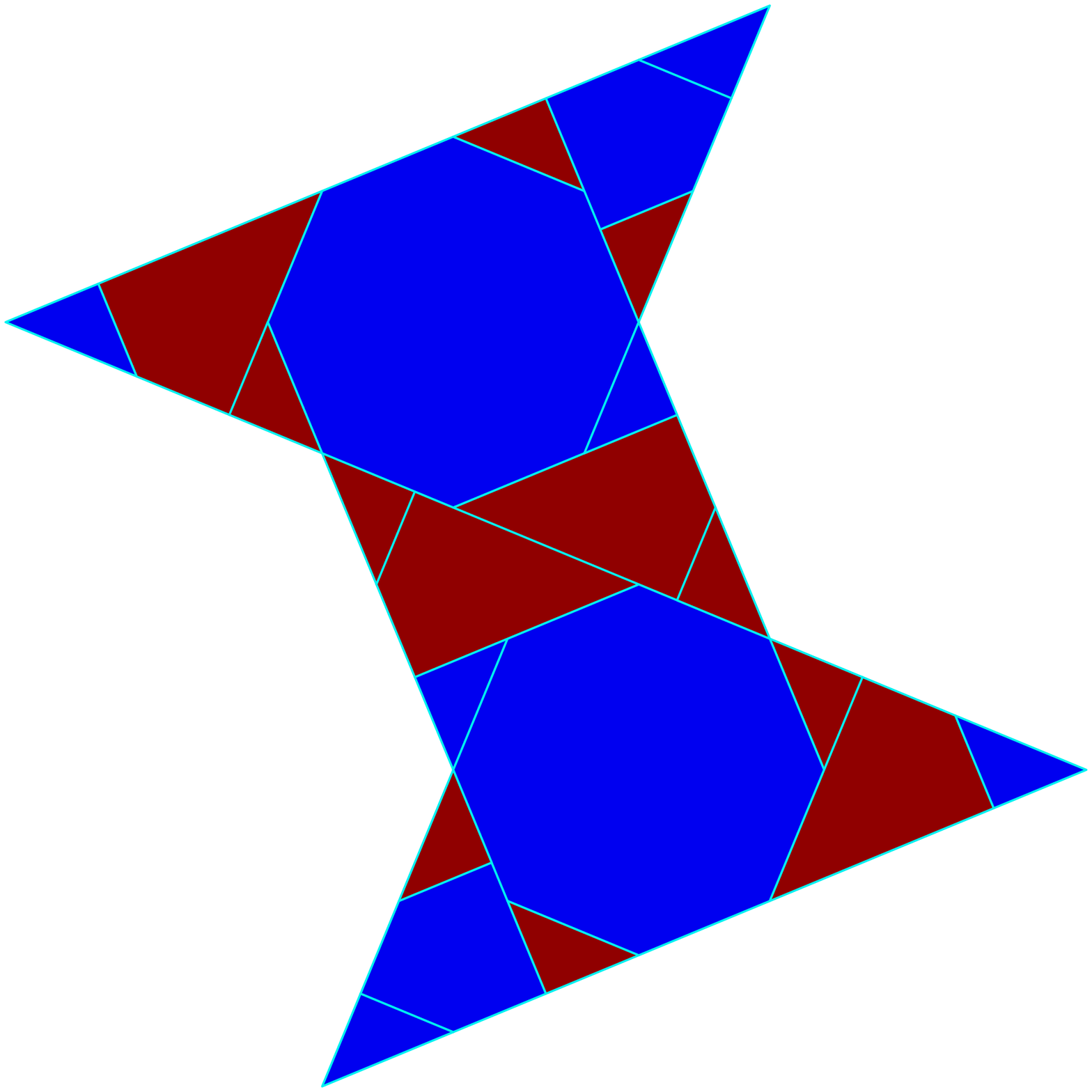}}
\newline
{\bf Figure 7.3:\/} The compressed system
\end{center}

The compressed system has $22$ regions.  The regions
$B_0,...,B_4$ are as in Figure 7.2.   For the record, here are the remaining regions.
\begin{itemize}
\item $B_5=\rho_3^5(B_1)$ and $B_6=\rho_3^5(B_2 \cup B_3)$ and $B_7=\rho_3^5(B_4)$.
\item $B_8=\rho_3^2(B_1)$ and $B_9=\rho_3^2(B_2 \cup B_3)$ and $B_{10}=\rho_3^2(B_4)$.
\item $B_{k+11}=\rho_2^4(B_k); \hskip 30 pt k=0,...,10$.
\end{itemize}
The regions $B_5,B_6,B_7$ all lie in the upper right corner of $R$.
The regions $B_8,B_9,B_{10}$ lie in the middle of $R$ on the left side.
The remaining $11$ regions are symmetrically placed with respect
to the first $11$.   Again for the record, the red colored regions
are
\begin{equation}
B_{i+11j}; \hskip 30 pt
i=2,4,6,7,8,9; \hskip 30 pt j=0,1.
\end{equation}
We will use this notation system in the next chapter, just
so that we can list things in a completely explicit way.

Figure 7.3 does for $\Psi$ what Figure 7.1 does
for $\Phi$.  The left hand side shows the
partition of $R$ into maximal regions on which
$\Psi$ is constant.  The blue polygons are the
ones on which $\Psi=\Phi$ and the
red polygons are the ones on which
$\Psi=\sigma \circ \Phi$.
The right hand side shows the images of
these pieces under $\Psi$.  Once we
mention that $\Psi$ is a piecewise translation
that maps all the triangles
on the left (respectively right) side of $R$ to triangles on the
right (respectively left) side of $R$, the partitions alone
determine $\Psi$.    We also remark that
$\Psi$ maps  the  kites and pentagons on
the left (respectively right)
side of $R$ to the left (respectively right) side of $R$.

The map $\Psi$ has a symmetry that it inherits from
$\Phi$, and also a new one.  First, $\Psi$ commutes
with the rotation $\rho_2$ about the center point of
$R$.   Second, the partition on the right hand side
of Figure 7.3 is a mirror image of the one on the left.
(The line of symmetry joins the centers of the
red octagons.)
Finally, $\Psi$ moves the one partition to
its mirror image by way of translations.

\subsection{Classification of Periodic Orbits}
\label{renorm}

\begin{lemma}
\label{period1}
$p \in R$ is a periodic point for $\Psi$ if and only
if $p$ is a periodic point for $\Phi$. 
\end{lemma}

\startproof
We write $\phi_0=\psi_0=p$.  Let
$\phi_n$ be the $n$th iterate of $\phi_0$ under $\Phi$ and
let $\psi_n$ be the $n$th iterate of $\psi_0$ under $\Psi$.
Since $\Phi$ commutes with $\sigma$, we have
\begin{equation}
\phi_n=\sigma^{d_n}(\psi_n)
\end{equation}
where $d_n$ counts the number of $k \in \{0,...,(n-1)\}$ for which
$\psi_k$ lies in a blue region in the partition of $R$ (on the left
side of Figure 7.3.)

Suppose that $p$ has $\Phi$-period $n$.  
Then $\phi_n=\phi_0$.  This means that
$\sigma^{d_n}(\psi_n)=\psi_0$.
Since both $\psi_n$ and $\psi_0$ belong to
$R$, we must have $d_n$ even.  Hence
$\psi_n=\psi_0$.   This shows that the
$p$ is periodic for $\Psi$.
On the other hand,
suppose that $p$ has $\Psi$-period $n$.
Then $\psi_n=\psi_0$.
If $d_n$ is even then $\phi_n=\phi_0$.
If $d_n$ is odd, then $d_{2n}=2d_n$ because
$\psi_{n+k}=\psi_k$ for each $k=0,...,n$.
Hence $\phi_{2n}=\phi_0$.  
Hence $p$ is $\Phi$-periodic.
\endproof

Figure 7.4 shows the beginning stages of a packing of octagons
in $R$.  As we will see momentarily, every periodic point
of $\Psi$ lies in one of the octagons in this packing.

\begin{center}
\resizebox{!}{5.2in}{\includegraphics{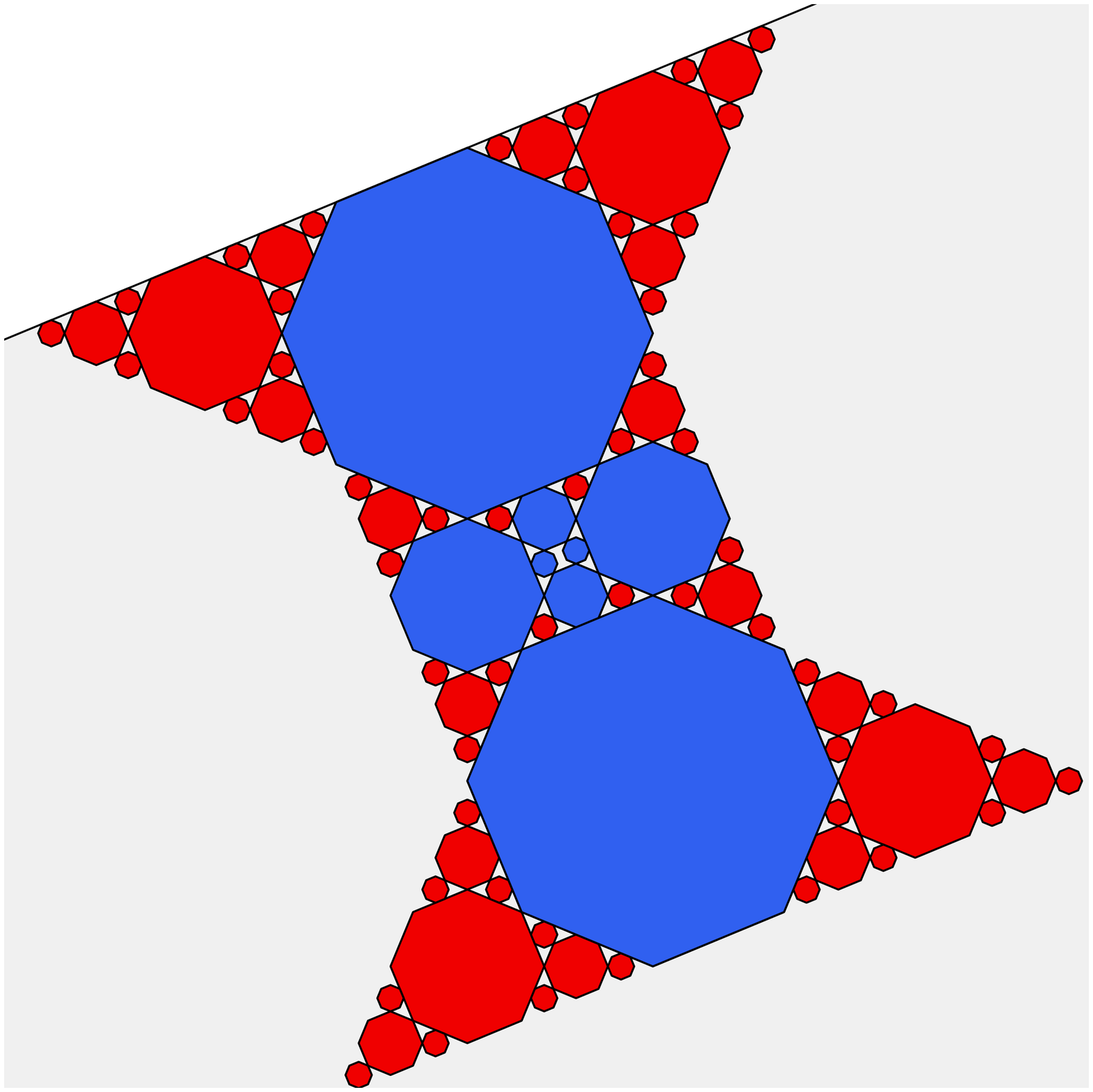}}
\newline
{\bf Figure 7.4:\/} Image of the central portion
\end{center}

We will establish the result about the packing using a
renormalization scheme that is very similar to the
one used in the previous chapter.
The remormalization scheme is based on the dilation
\begin{equation}
\Theta(z)=\overline{a z + b}; \hskip 30 pt
a=(s+1)+i(s+1); \hskip 30 pt
b=-3i(s+1)
\end{equation}
We have $\Theta(S)=R$, where $S \subset R$ is the portion of $R$
between the two octagons $O_1,O_2 \subset R$.  The
left hand side of Figure 7.4 below shows $R$ in blue.
We will use the renormalization scheme to classify the
periodic points of $\Psi$.

\begin{lemma}
\label{period2}
Let $O_1^k=\Theta^{-k}(O_1)$
and $O_2^k=\Theta^{-k}(O_2)$.  
Then
\begin{itemize}
\item There is a self-similar packing of regular octagons in $S$,
each of which is translation equivalent to $O_1^k$ (or $O_2^k$)
for some $k$.  (See Figure 7.4.)
\item Every periodic point of $\Psi$ belongs to one of the octagon
interiors, and its period is the same as the
tile-period of the octagon containing it.
\item Every octagon in the packing has tile period $3^k$ for some
$k=0,1,2,...$, and such an octagon either lies in the orbit of
$O_1^k$ or the orbit of $O_2^k$.  These are the blue octagons
in Figure 7.4.
\end{itemize}
\end{lemma}

\startproof
We see by inspection that $\phi^3(S') \subset S'$.  Indeed,
the left hand side of Figure 7.5 shows $S$ in blue,
the middle picture shows $\Psi(S)$ in blue, and the
right hand side shows $\Phi^2(S)$ in blue.

\begin{center}
\resizebox{!}{1.7in}{\includegraphics{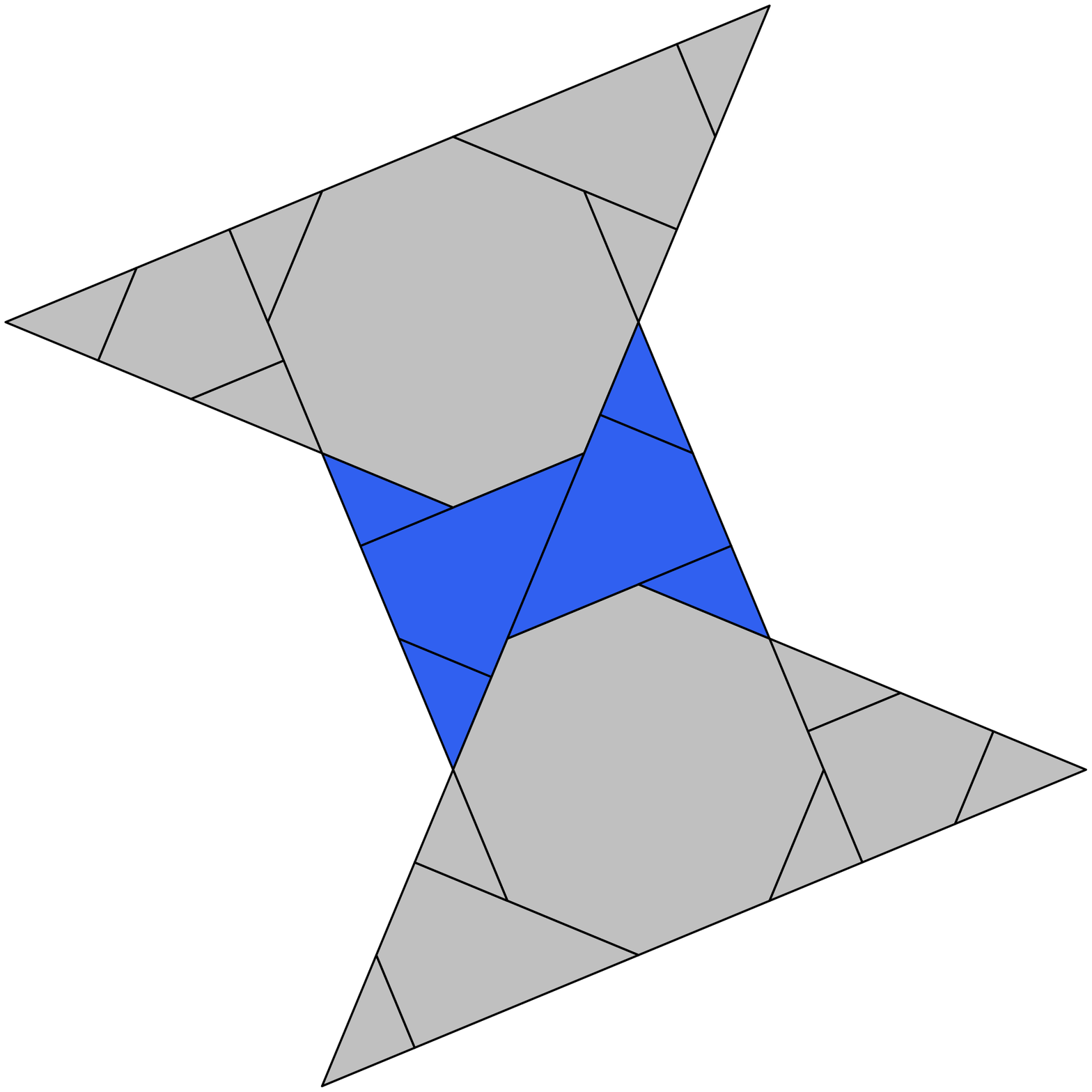}}
\resizebox{!}{1.7in}{\includegraphics{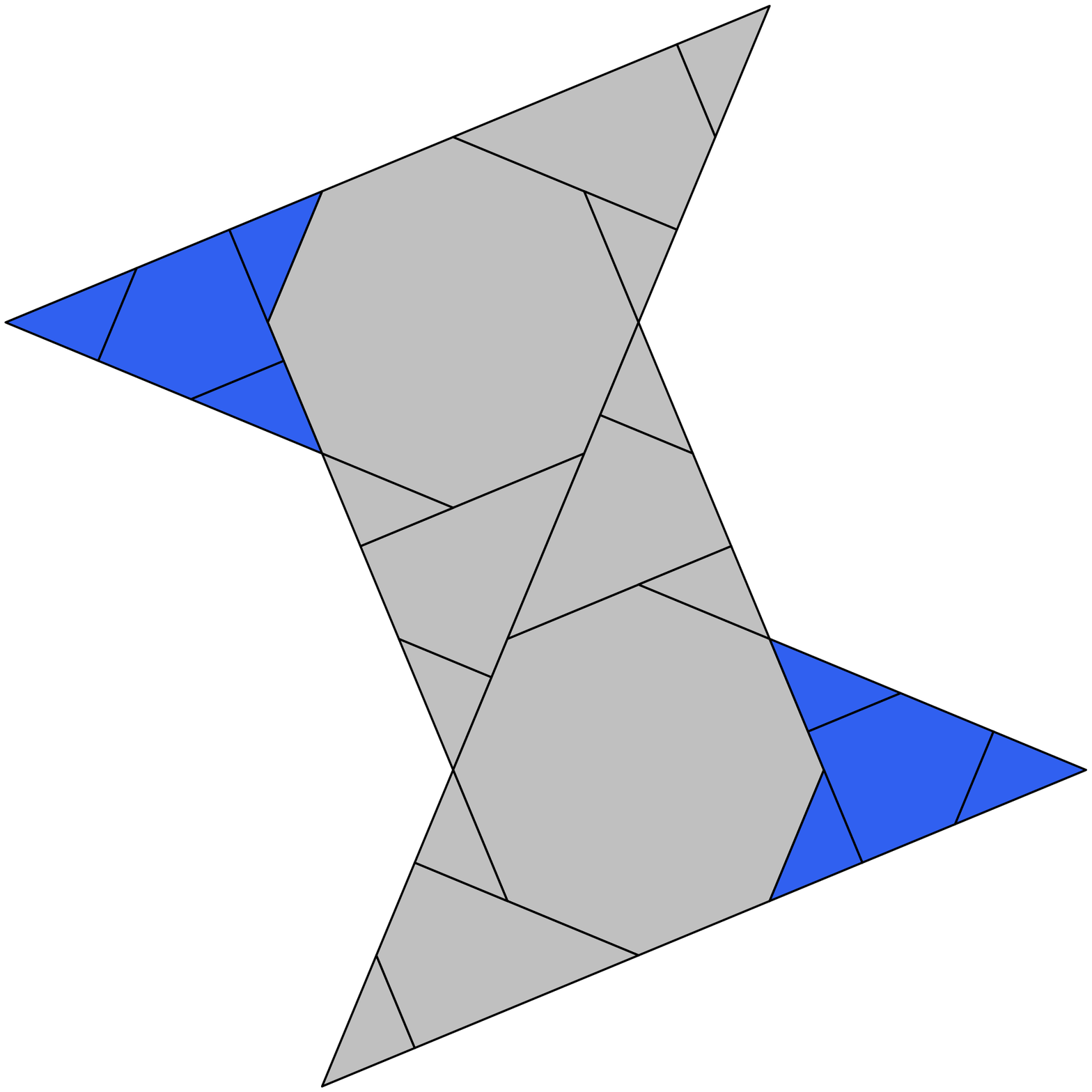}}
\resizebox{!}{1.7in}{\includegraphics{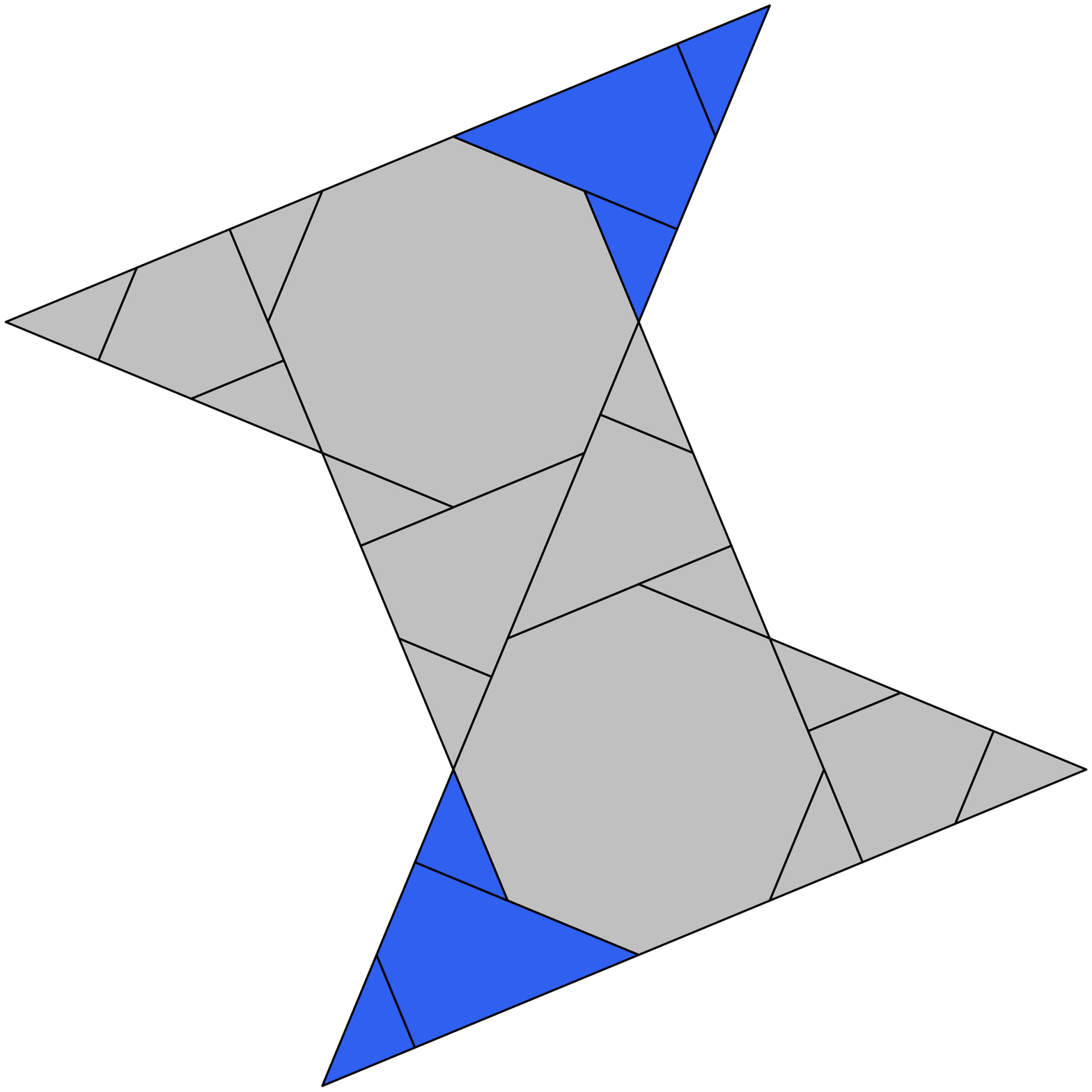}}
\newline
{\bf Figure 7.5:\/} Image of the central portion
\end{center}

\begin{lemma}
Every $\Psi$-periodic point in $R$ has period $3^k$
for some $k=0,1,2...$
\end{lemma}

\startproof
A short calculation establishes the following two facts.
\begin{equation}
\label{renorm3}
\Theta^{-1} \circ \Psi \circ \Theta|_S=\Psi^3; \hskip 30 pt
R-O_1-O_2=S \cup \Psi(S) \cup \Psi^2(S).
\end{equation}
The same analysis as in the previous chapter produces the
packing of octagons, together with the claims about
their tile periods.
The only new feature of the situation here is that the map
$\Psi$ is a piecewise translation, whereas the map from
the previous chapter had a rotational component.  Since
$\Psi$ is a piecewise translation, each point in a octagonal
tile has exactly the same period as the tile itself.
\endproof

\subsection{Proof of Lemma \ref{ocount}}
\label{countproof}

Thanks to the analysis done in \S \ref{small} it
suffices to prove that
every periodic point of $R_1$ has $\Phi$-period $3^k$ for
some $k=0,1,2,...$. 
We have already worked out the structure
of the periodic points of the $\Psi$-system,
and we will use this to get the desired
information about the periodic points of
the $\Phi$-system.

There is a $2$ to $1$ map $F: R_1 \to R$.
which, by construction, carries periodic points
of $\Phi$ to periodic points of $\Psi$.
Moreover, this map conjugates the 
action of $\Phi$ on $R^1$ to the
action of $\Psi$ on $R$.   From this,
we see that every $\Phi$-periodic
orbit has period either $3^k$ or
$2 \times 3^k$.

A parity check rules out the second possibility.
Given a $\Psi$-periodic orbit
$p_0,p_1,...,p_{N-1}$, we let
$\sigma_i=0$ if $p_0$ lies in a blue region of
$R$ and other wise $\sigma_i=1$.  We just
need to check that
$\sum_{i=0}^{N-1} \sigma_i$ is even.
We check that the parity of the sum
is preserved under the renormalization operation.
Since the sum obviously works out for the
large octagons -- i.e., when $N=1$, it works
out in general.

The parity check is a finite calculation, as we now explain.
Suppose we start with a point $z \in R$.  Let
$w_0=\Theta^{-1}(R)$.    Let
$w_j=\Psi^j(w_0)$.   We just have to prove that
\begin{equation}
\sigma \equiv \sigma_0+\sigma_1+\sigma_2 \hskip 15 pt {\rm mod\/}\ \ 2.
\end{equation}
Given the nature of $\Psi$, we only need to check this
equation for $1$ point in each of the $22$ regions.
Here are the results for the first $11$ regions
$$\matrix{
0: & 0 && 1&1&1 \cr
1: & 0 && 0&0&0 \cr
2: & 1 && 0&0&1 \cr
3: & 0 && 0&0&0 \cr
4: & 1 && 1&0&0 \cr
5: & 0 && 1&0&1 \cr
6: & 1 && 1&1&1 \cr
7: & 1 && 1&1&1 \cr
8: & 1 && 1&1&1 \cr
9: & 1 && 1&1&1 \cr
10: & 0 && 1&0&1}$$
The pattern repeats exactly for the second $11$ regions.

This completes the proof of Lemma \ref{ocount}.

\newpage

\section{Generating the Arithmetic Graph}

\subsection{The Basic Approach}

Let $R_1$ be our basic region, partitioned as on the left hand side of Figure 7.1.
To each region $r \in R_1$ we choose and $p \in r$ and define
\begin{equation}
g_k(r)=\pi_k(\widehat \Phi(p)).
\end{equation}
Here $\widehat \Phi(p) \in \Z^8$ is the vector  from Equation \ref{spectrum}. 
 Since $\widehat \Phi$ is constant on the region $r$,
we see that our definition is independent of the choice of $p \in r$.
If $p$ is part of some orbit, then $g_k(p) \in \R^2$ is the side of
the edge of $\pi_k(\Gamma)$ that corresponds to $(p,\Phi(p))$.

Given an initial point $p_0 \in R_1$ we consider the forward orbit $p_0,p_1,p_2,...$  We
have a sequence of regions $r_0,r_1,r_2,...$ and a corresponding sequence
$g_k(r_0), g_k(r_1), g_k(r_2),...$ 
These vectors are the successive sides of the projection $\pi_k \circ \Gamma(p)$.

\subsection{The Compressed Approach}

Given the nice renormalization scheme we discussed in the previous
chapter, we prefer to work with the compressed system
$\Psi: R \to R$, as defined in \S \ref{compress}.  The region $R$ is
colored as in Figure 7.3.   

Recall that $R$ is the top half of $R_1$.
Let $R'$ be the bottom half. 
Let $p_0,p_1,p_2,...$ be a $\Psi$ orbit in $R$.  Let
$q_0,q_1,q_2$ be a $\Phi$ orbit in $R_1$ such that
$p_0=q_0$.   Let $\sigma_i=0$ if $p_i$ lies in a
blue region of $R$ on the left hand side of Figure 7.3,
and $\sigma_i=1$ otherwise.

We define the {\it accumulated parity\/}
\begin{equation}
s_n=\bigg(\sum_{i=0}^{n-1} \sigma_i \bigg)\ {\rm mod\/}\ 2.
\end{equation}
That is, we count the number of times that the orbit lands in a red
region on the left hand side of Figure 7.3.
The sequence $\{s_n\}$ is a binary sequence that keeps track
of whether $p_n=q_n$ or not.

The accumulated parity is the additional information that lets us
``lift'' an $\Psi$ orbit on $R$ to a $\Phi $ orbit on $R_1$.

\begin{lemma}
$s_n=0$ if and only if $p_n=q_n$. 
\end{lemma}

\startproof
This follows from two facts.  First, $\sigma_i=0$ iff
$\Phi(q_i)$ and $q_i$ lie in the same half of $R_1$,
Second, $\Phi: R_1 \to R_1$ commutes
with the map interchanging $R$ and $R'$.
\endproof

To each region $r \in R$, we associate $4$ vectors according to the
following rule.
\begin{equation}
g_k(r,0)=g_k(r); \hskip 30 pt
g_k(r,1)=g_k(r'); \hskip 30 pt
r'=r-(0,1) \in R'.
\end{equation}
Given the orbit $p_0,p_1,p_2 \in R$ we have the corresponding regions
$r_0,r_1,r_2,...$ and the corresponding accumulated parities
$s_0,s_1,s_2,...$.   By construction, $g_k(r_n,s_n)$ is the $n$th
edge of $\pi_k \Gamma(p_0)$.

For reference, we list the vector assignments.
It is convenient to let $k$ stand for the (region,parity) pair
$(k,0)$ and let $k+22$ stand for the (region,parity) pair
$(k,1)$.   The indices $0,11,22,33$ correspond to the biggest octagons
and never arise, so we leave them off \footnote{For the record, $g_2(0)=(0,8)$ and
$g_3(0)=(0,0)$.} 
We consider the $g_2$ vectors first.  Letting
$T(x,y)=(-y,-x)$, we have the symmetries
\begin{equation}
\label{g2symm}
g_2(k+11)=T(g_2(k)); \hskip 30 pt
g_2(k+22)=-g_2(k); \hskip 30 pt k=0,...,10.
\end{equation}
So, the following assignment determines all $44$ vectors.
\begin{enumerate}
\item $g_2(1)=(6,2)$.
\item $g_2(2)=(2,-2)$.
\item $g_2(3)=(2,-2)$.
\item $g_2(4)=(-4,4)$.
\item $g_2(5)=(-6,2)$.
\item $g_2(6)=(0,-4)$.
\item $g_2(7)=(2,6)$.
\item $g_2(8)=(4,4)$.
\item $g_2(9)=(2,2)$.
\item $g_2(10)=(-2,6)$.
\end{enumerate}

In the case of $g_3$ we have the symmetries
\begin{equation}
\label{g3symm}
g_3(k+22)=g_3(k); \hskip 20 pt
g_3(k+11)=-g_3(k)+\sigma_k(0,4); \hskip 20 pt k=0,...,10.
\end{equation}
Here $\sigma_k \in \{0,1\}$ is the parity of the $k$th region.
We have
\begin{equation}
\sigma_k=1 \hskip 20 pt
\Longrightarrow \hskip 20 pt k \in \{2,4,6,7,8,9\}.
\end{equation}
In light of these symmetries, the following list
of $10$ vector assignments determines the whole thing.

\begin{enumerate}
\item $g_3(1)=(0,-1,-2,-1)$.
\item $g_3(2)=(2,1,0,-1)$.
\item $g_3(3)=(2,1,0,-1)$.
\item $g_3(4)=(-2,-2,2,0)$.
\item $g_3(5)=(-2,-1,0,-1)$.
\item $g_3(6)=(0,0,0,-2)$.
\item $g_3(7)=(-2,-1,4,1)$.
\item $g_3(8)=(-2,-2,2,0)$.
\item $g_3(9)=(-2,-1,0,-1)$.
\item $g_3(10)=(0,-1,2,1)$.
\end{enumerate} 
Our notation is as follows.
\begin{equation}
(a,b,c,d)=(a+b\sqrt 2,c+d\sqrt 2).
\end{equation}

In the next section we will explain how these vector
assignments combine with a combinatorial substitution
rule to generate the graphs.

\subsection{The Substitution Approach}
\label{sub}

Suppose that we have some orbit
$p_0,p_1,p_2,...$.  We let
$$p'_0=\Theta^{-1}(p_0)$$
and consider the renormalized orbit $p_0',p_1',p_2',...$.
Let
$s_j=\sigma_0+...+\sigma_{j-1}$
be the accumulated parity of the original orbit.
Likewise define $s'_j$.
Let $r_j$ be the region of $R$ containing $p_j$, and likewise
define $r_j'$.   

In \S \ref{countproof} we explained the sense in which
the renormalization scheme respects the accumulated parities.
This fact implies that
\begin{equation}
\label{paritykey}
s'_{3j}=s_j
\end{equation}
Given equation \ref{paritykey}, the pair
$(s_j,r_j)$ determines the triple
$$(r'_{3j},s'_{3j}) \hskip 15 pt 
(r'_{3j+1},s'_{3j+1}) \hskip 15 pt 
(r'_{3j+2},s'_{3j+2}).$$
Recall that we have chosen a particular listing for these
$44$ possible (region, parity)  pairs.   What we
are saying is that the numerical code for the original orbit
determines the numerical code for the renormalized orbit,
by way of a substitution scheme in which one number is
replaced by $3$ numbers.
Here is one quarter of the substitution rule.

$$\matrix{
0 & \to & 9&25&39 \cr
1 & \to & 10&14&5 \cr
2 & \to & 10&12&17 \cr
3 & \to & 10&12&16 \cr
4 & \to & 9&23&27 \cr
5 & \to & 8&36&28\cr
6 & \to & 8&35&17 \cr
7 & \to & 9&24&6 \cr
8 & \to & 9&26&6 \cr
9 & \to & 9&26&7 \cr
10 & \to & 9&25&40}$$

The rest of the pattern can be deduced from symmetry
and from what we have already written.  We point out
two involutions on the set $\{0,...,43\}$.  The first of these
has the action $n \to n \pm 22$.  The second one has
the action $n \to n \pm 11$, where the $+$ option is
taken if $n \in \{0,...,10\}$ or $n \in \{0,...,10\}+22$.
The substitution rule is invariant under these two
involutions.  For instance, 
$$\matrix{4 & \to & 9&23&27 \cr
15 & \to &20&34&38 \cr
26 & \to &31&1&5    \cr
37 & \to &42&12&16}.$$

Now we can explain a nice combinatorial way to generate the
arithmetic graph projections.
Referring to Lemma \ref{period2}, the two orbits
$O_1^k$ and $O_2^k$ both have period $3^k$.
Here $O_1^0$ and $O_2^0$ respective are the top and
bottom biggest octagons in Figure 7.4.  We will concentrate
on the orbits $O_1^k$.   The other half of the orbits have
a very similar treatment.

We start with $0$, and iterate the substitution rule $k$ times.
Then we replace each number by the relevant vector.  The
resulting list of vectors is the set of edges of the projection
of the arithmetic graph.  We call the resulting path
$G_2(k)$ or $G_3(k)$, depending on whether we
use the $G_2$ vectors or the $G_3$ vectors.

Here are two examples.  For $G_2(2)$ we make the substitution
$0 \to 9,25,39$, and then plug in the $G_2$ vectors:
\begin{itemize}
\item $9 \to (2,2)$.
\item $25 \to (-2,2)$.
\item $39 \to (0,-4)$.
\end{itemize}
We then get the vertices of $G_2(1)$ by accumulating the vectors:
$$(0,0), (2,2), (0,4), (0,0).$$
So, $G_2(1)$ is a closed path -- a triangle.

To  generate $G_2(2)$ we go one more step.
$$0 \to 9,25,39 \to 9,26,7,32,34,38,41,2,28.$$
Then we make all the substitutions and accumulate the vectors.
We discover in this case that
$G_2(2)$ is not closed.  
We check that $G_2(k)$ is closed if $k=1,3,5,7$ and open if $k=2,4,6,8$.
In \S \ref{closed} we will prove that this pattern persists for all $k$.
Our proof there, in particular, implies Statement 1 of the Main Theorem.

\newpage

\section{Taking a Better Limit}

\subsection{Subsequential Limits}
\label{gammadefine}

Note that $G_2(k)$ and $G_3(k)$ have natural parametrizations, coming
from the ordering on the edges.   Define
\begin{equation}
\Gamma_2(k)=(\sqrt 3)^{-k}G_2(k); \hskip 30 pt
\Gamma_3(k)=(1+\sqrt 2)^{-k}G_3(k).
\end{equation}
We notice that 
$\Gamma_2(k)$ and $\Gamma_2(k+\beta)$ closely follow each
other as parametrized curves when $\beta=4$. 
This property does not hold for smaller $\beta$.   When
$\beta=2$, the two curves $\Gamma_2(k)$ and $\Gamma_2(k+\beta)$
are close in the Hausdorff metric, but not close as parametrized
curves.
We also notice that $\Gamma_3(k)$ and $\Gamma_3(k+\beta)$ closely
follow each other closely when $\beta=2$ but not when
$\beta=1$.   

Note that $\beta=4$ is the lowest number that works well for both
curve families.  This accounts for the dependence on the mod $4$ congruence
that appears in the Main Theorem, which involves the subsequential limits
\begin{equation}
\lim_{j \to \infty} \Gamma_i(j+4k); \hskip 30 pt i=2,3; \hskip 30 pt j=1,3.
\end{equation}

The $4$th power of our substitution
rule is our basic object of study.  With this rule, each number is replaced
by $81$ numbers.  In geometric terms, each $G_2$ vector $g_2(n)$ is replaced
by a polygonal path $g'_2(n)$ of length $81$.  This path consists of $81$
other $G_2$ vectors, scaled down by $\sqrt 3$.  Likewise, each $G_3$ vector
$g_3(n)$ is replaced by a polygonal path $g'_3(n)$ of length $81$.   

Let $Rg_2(n)$ denote the vector that points
from the tail of $g'_2(n)$ to the head of $g'_2(n)$.
Likewise define $Rg_3(n)$.
We notice that 
$g_i(n)$ and $Rg_i(n)$ are nearly the same vector
for each choice $i=2,3$ and each $n=0,...,43$.
We think of $R$ as a {\it renormalization operator\/},
which replaces one vector by another one that
is very close to it.  It makes sense to define
$R^kg_i(n)$ as follows: 
\begin{itemize}
\item  We iterate the
substitution rule $4k$ times.
\item We replace
each number by the appropriate $g_i$ vector.
\item We rescale by $(s_i)^{-k}$, where
$s_2=\sqrt 3$ and $s_3=1+\sqrt 2$.
\end{itemize}
The key to analyzing our limits is to looked at the
fixed point of the renormalization operator.

\subsection{The Fixed Point of Renormalization}

For motivational purposes (only) we mention the Perron-Frobenius Theorem.

\begin{theorem} [Perron-Frobenius]  Let
$T: \R^n \to \R^n$ be a linear transformation
represented by a matrix with all positive entries.
Then $T$ has (up to scale) a unique unit positive
eigenvector $u$.  Moreover, if $v$ is any
positive vector $v$, then we have $v_n \to u$,
where $v_n=T^n(v)/\|T^n(v)\|$.
\end{theorem}

It turns out that $g_i(n)$ and $Rg_i(n)$ are
not in general the same vector. 
The Perron-Frobenius Theorem suggests to us that we
consider the alternate assignment
\begin{equation}
\label{limitvector}
n \to \lambda_i(n):=\lim_{s \to \infty} R^sg_i(n).
\end{equation}
Doing things this way is analogous to
finding the Perron-Frobenius eigenvector of a transformation $T$
by taking the limit of iterates of $T$ applied to some conveniently
chosen initial vector $v$.  In case the limit $u$ is expected
to have nice properties -- e.g. algebraic  coordinates -- one might hope
to guess $u$ by looking at the decimal expansions of a
fairly high iterate. Once we have made a guess, it is trivial
to verify that it works.  This is what we did.

Our new assignment for $G_2$ has the same symmeties as in
Equation \ref{g2symm},
so the $10$ vectors we list suffices
to describe the whole thing.
\begin{enumerate}
\item $\lambda_2(1)=^*(6,2)$.
\item $\lambda_2(2)=(3,-3)$.
\item $\lambda_2(3)=^*(2,-2)$.
\item $\lambda_2(4)=(-3,3)$.
\item $\lambda_2(5)=^*(-6,2)$.
\item $\lambda_2(6)=(-3,-1)$.
\item $\lambda_2(7)=(3,5)$.
\item $\lambda_2(8)=(5,3)$.
\item $\lambda_2(9)=(3,1)$.
\item $\lambda_2(10)=^*(-2,6)$.
\end{enumerate}
This assignment has the same symmetries as the
original one, and these symmetries let one deduce
the remaining vectors in the assignment.
The starred equalities indicate that no change has been
made.  Curiously, we have
\begin{equation}
\label{oldnew}
{\rm new\/} = {\rm old\/} + \epsilon (1,-1); \hskip 30 pt
\epsilon \in \{-1,0,1\}.
\end{equation}

In the case of $G_3$, the new assignment has nicer symmetries
than the original namely
\begin{equation}
\lambda_3(k)=-\lambda_3(k+11)=\lambda_3(k,1).
\end{equation}
That is, the factor of $\sigma_k$ is replaced by $0$ in
Equation \ref{g3symm}.  In light of these symmetries,
our $10$ vectors below determine the whole thing.

\begin{enumerate}
\item $\lambda_3(1)=^*(0,-1,-2,-1)$.
\item $\lambda_3(2)=(2,1,-2,-1)$.
\item $\lambda_3(3)=^*(2,1,0,-1)$.
\item $\lambda_3(4)=(-2,-2,0,0)$.
\item $\lambda_3(5)=^*(-2,-1,0,-1)$.
\item $\lambda_3(6)=(0,0,-2,-2)$.
\item $\lambda_3(7)=(-2,-1,2,1)$.
\item $\lambda_3(8)=(-2,-2,0,0)$.
\item $\lambda_3(9)=(-2,-1,-2,-1)$.
\item $\lambda_3(10)=^*(0,-1,2,1)$.
\end{enumerate} 
Again, the starred equalities indicate the ones where there is
no change.
This time we have
\begin{equation}
\label{oldnew2}
{\rm new\/} = {\rm old\/} + \epsilon (0,2); \hskip 30 pt
\epsilon \in \{-1,0,1\}.
\end{equation}

A direct calculation shows that
\begin{equation}
\label{fixed}
R\lambda_i(n)=\lambda_i(n); \hskip 30 pt i=2,3; \hskip 30 pt
n=1,...,43.
\end{equation}
As usual, we are omitting $n=0,11,22,33$.

Let $\Lambda_i(k)$ denote the version of $\Gamma_i(k)$ produced by the
new vector assignments.  Our new assignments are nicer, but we need to
prove that the $\Lambda$ paths have the same limits as the
$\Gamma$ paths.   Once we know this, we can throw out the
arithmetic graphs, and work with the much nicer $\Lambda$ paths.
The rest of the chapter is devoted to analyzing and equating the
relevant limits.

Figure 9.2 shows $\Gamma_3(3)$ and $\Gamma_3(7)$.
Notice that these curves follow each other pretty well, but that
a few of the vertices of $\Gamma_3(7)$ do not lie on
$\Gamma_3(3)$.   Figure 9.3 shows the improved
$\Lambda_3(3)$ and $\Lambda_3(7)$.   These curves
are better related to each other.    Note also that
$\Gamma_3(7)$ and $\Lambda_3(7)$ are really just
about the same curve.

\begin{center}
\resizebox{!}{5in}{\includegraphics{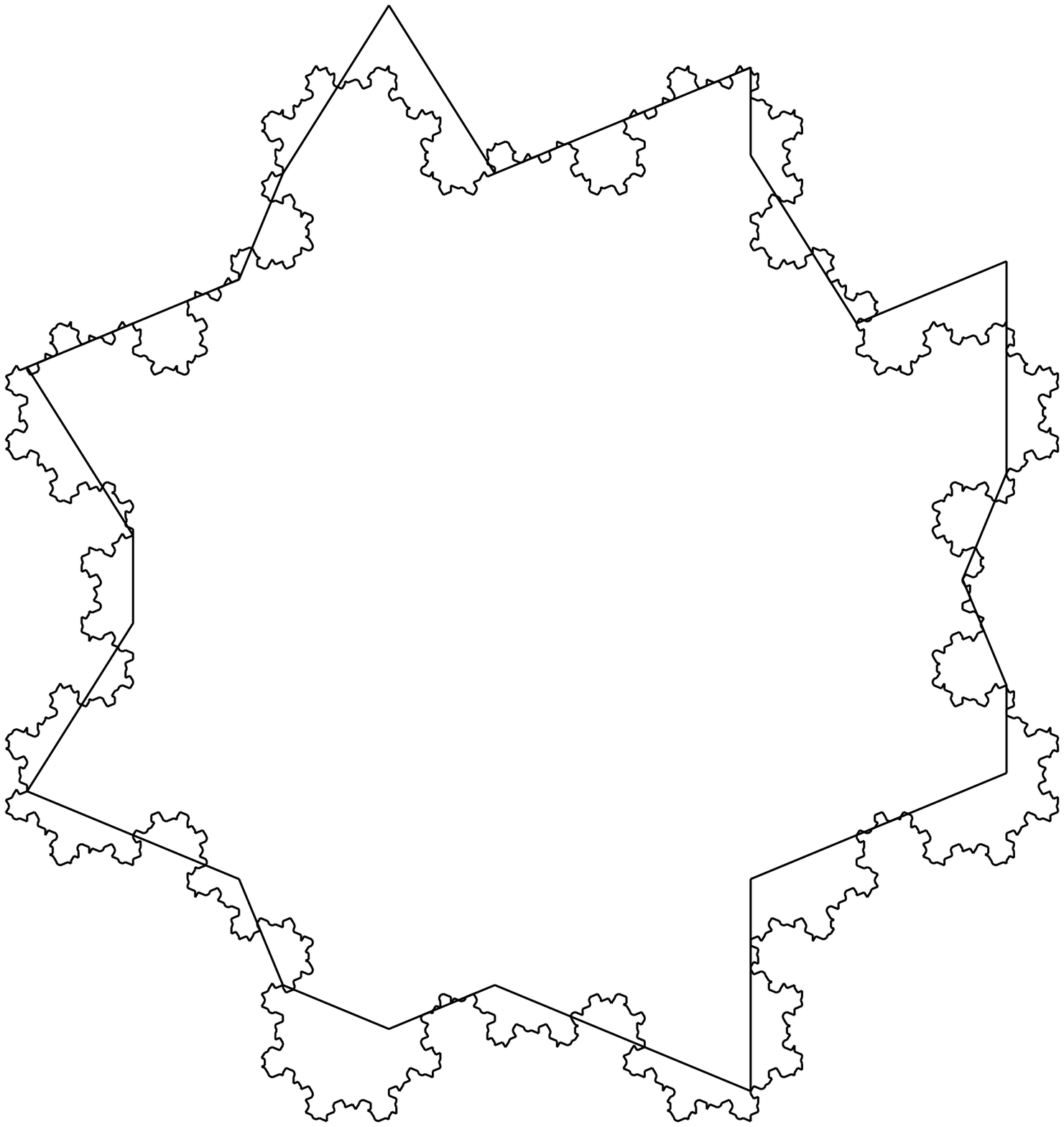}}
\newline
{\bf Figure 9.2:\/}  $\Gamma_3(3)$ and $\Gamma_3(7)$.
\end{center}  

\begin{center}
\resizebox{!}{5in}{\includegraphics{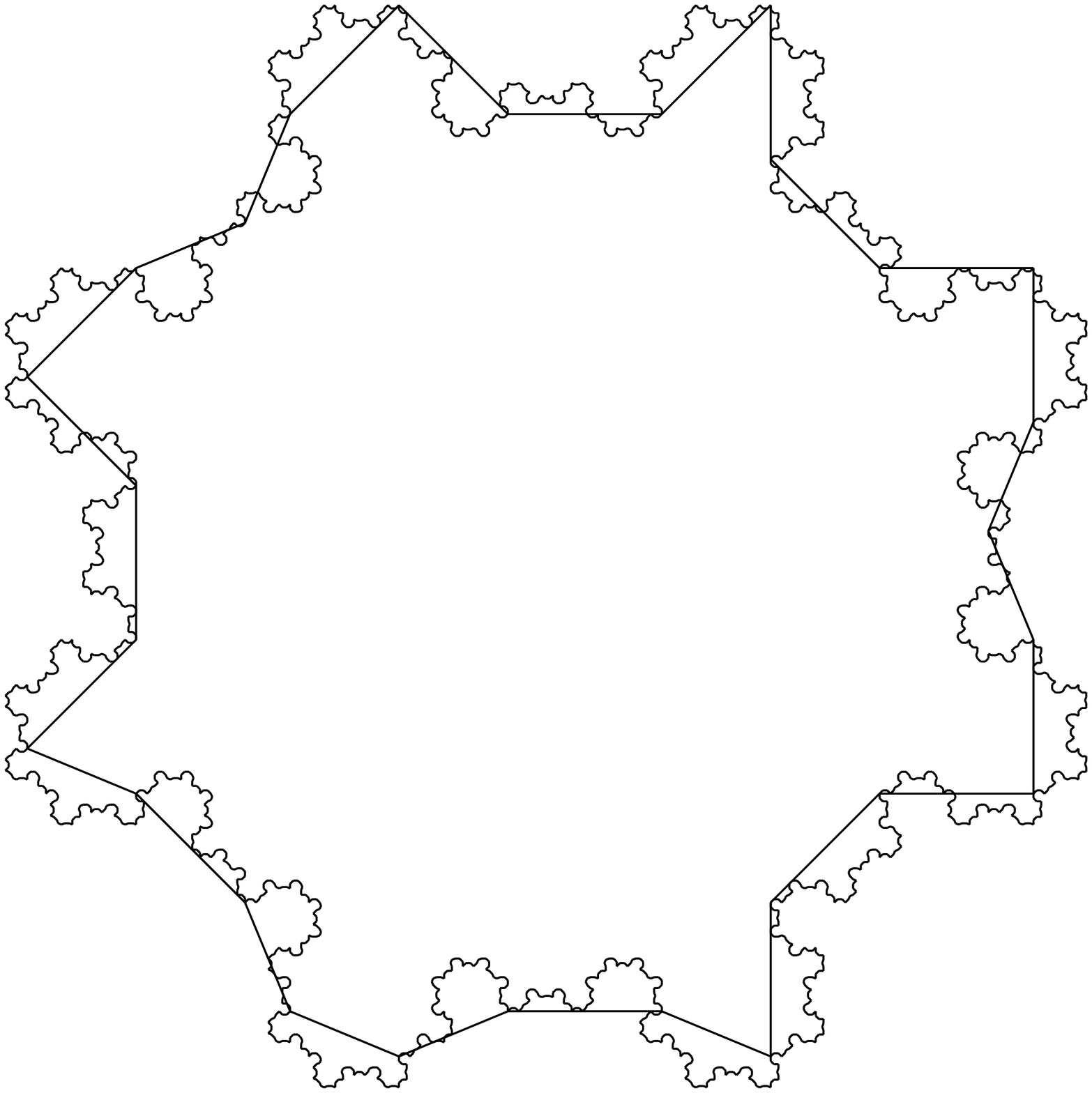}}
\newline
{\bf Figure 9.3:\/}  $\Lambda_3(3)$ and $\Lambda_3(7)$.
\end{center}  

\subsection{Proof of Statement 1}
\label{closed}

Statement 1 of the Main Theorem, which only deals with the
odd case, follows from the more complete result we now prove.

\begin{lemma}
\label{disc}
The following is true.
\begin{itemize}
\item $G_3(k)$ is a closed polygon for all $k=1,2,3,...$
\item $G_2(k)$ is a closed polygon if and only if $k$ is odd.
\end{itemize}
\end{lemma}

\startproof
Let $L_k$ be the unscaled version of $\Lambda_k$.
The main idea in our proof is showing that $G_k$ and $L_k$ have
the same endpoints.  Assuming this result, it suffices to prove
our result for the $L$ curves.  But the $\Lambda$ curves
are just scaled version of the $L$ curves.
So, it suffices to prove our result for the $\Lambda$ curves.
Given the endpoint-preserving property of our ideal substitution
rule -- see Equation \ref{fixed} -- we see that
$\Lambda_i(k)$ is closed if and only if $\Lambda_i(k+4)$ is closed.
We check the conclusion of this lemma for the first few cases, and
then apply the induction step to take care of the remaining cases.

Now we come to the interesting part of the proof, showing that each
$G$ curve has the same endpoints as the corresponding $L$ curve.
 We treat $G_2$ and $L_2$ first. 
All paths start at the origin, so we just have to see that the other
endpoints match.  That is, we need to prove that
\begin{equation}
\label{total}
\sum_{i=0}^{N-1} \delta(n_i) = 0; \hskip 30 pt
\delta(n_i)=g_2(n_i)-\lambda_2(n_i).
\end{equation}
Here $\{n_i\}$ is the numerical code for the relevant 
pair of paths, and $N$ is the total combinatorial length.

As we remarked above -- see Equation \ref{oldnew} -- we always
have $\delta(n)=\epsilon (1,-1)$ for $\epsilon \in \{-1,0,1\}$.
We check by hand the following $44$ equations.
\begin{equation}
\label{conserve}
\delta(n)=\delta(m_1)+\delta(m_2)+\delta(m_3).
\end{equation}
Here $(m_1,m_2,m_3)$ is the substitution for $n$.
Equation \ref{conserve} implies that the sum in
Equation \ref{total} is invariant under substitution.
We check that the sum is $0$ in the first $4$ cases,
and then we see by induction that it is always $0$.

The proof for $G_3$ and $L_3$ works in the same way.
This time, we have
$\delta(n)=\epsilon(0,2)$ for $\epsilon \in \{-1,0,1\}$.
\endproof

\subsection{Proof of Statement 2}

Let $G$ and $L$ be the unscaled versions
of $\Gamma$ and $\Lambda$ respectively, as in
the proof of Lemma \ref{disc}.  Say that a
{\it distinguished strand\/} of $G$ is one that
comes from repeated substitution applied
to a single edge.   Likewise define distinguished
strands of $L$.  Such strands consist of
$3^k$ segments.
Any distinguished strand $\sigma$ determines
a vector $V(\sigma)$ that points from the
starting point to the ending point.

\begin{lemma}
Let $g$ and $l$ be corresponding distinguished
strands of $G(n)$ and $L(n)$.  Then
$V(g)$ and $V(l)$ differ by at most $2$ units.
\end{lemma}

\startproof
Let $g'$ and $l'$ be the distinguished strands of
the paths $G(n-1)$ and $L(n-1)$ which give rise to
$g$ and $l$, respectively, upon substitution.
The same argument as in Lemma \ref{disc} shows
that $$V(l')-V(g')=V(l)-V(g).$$
By induction, we are reduced to the case of single edges.
Now we apply Equations \ref{oldnew} and \ref{oldnew2}.
\endproof

\begin{corollary}
\label{log}
The distance between any corresponding vertices of
$G_i(n)$ and $L_i(n)$ is at most $4 \log_3(n)$.
\end{corollary}

\startproof
Both $G_i(n)$ and $L_i(n)$ start at $(0,0)$.   
Suppose we want to consider the $k$th vertices.
Let $\gamma$ be the arc of $G_i$ connecting
$(0,0)$ to the $k$th vertex.  Likewise define
$\lambda$.  

Using the base $3$ expansion of $k$ as a guide, we
can decompose $\gamma$ and $\lambda$ into
at most $2 \log_3(n)$ distinguished arcs.
Our lemma now follows from the previous result
and the triangle inequality.
\endproof

\begin{lemma}
\label{converge0}
Let $i$ be $2$ or $3$.
The sequences $\{\Lambda_i(1+4k)\}$ and
$\{\Lambda_i(3+4k)\}$, considered as a sequence of
parametrized curves, converge in the $L_{\infty}$ norm as $k \to \infty$.
\end{lemma}

\startproof
We will consider the first of these sequences.  The second one has
the same proof.
The edges of $\Lambda_i(4k+5)$ are obtained by replacing each edge
of $\Lambda(4k+1)$ by a length $81$ polygonal path that has the
same endpoint.  Define
\begin{equation}
D_k=\sup_{s \in [0,1]} {\rm dist\/}\Big(\Lambda_i(4k+1;s),\Lambda_i(4k+5;s).\Big).
\end{equation}
  Each edge of $\Lambda(4k+5)$ is $$1/9=(\sqrt 3)^{-4}$$
times as long as an edge of the same type in $\Lambda(4k+1)$.   By similarity,
we have
\begin{equation}
\label{exponentially}
D_{k+1} = \frac{1}{\sqrt 3} D_k,
\end{equation}
once $k$ is large enough that $\Lambda_i(1+4k)$ contains all the edge types.
It suffices to take $k=1$, in fact.

Equation \ref{exponentially} implies that the distance between
two successive curves in our sequence decreases exponentially.
Hence, our curves form a Cauchy sequence in the space
${\rm Map\/}(S^1,\R^2)$, equipped with the $L_{\infty}$
topology.  Hence, our sequence converges.
\endproof

The following
result implies Statement 2 of the Main Theorem.

\begin{lemma}
\label{s1proof}
The sequences $\{\Gamma_i(1+4k)\}$ and
$\{\Gamma_i(3+4k)\}$, considered as a sequence of
parametrized curves, converge in the $L_{\infty}$ norm as $k \to \infty$.
The limits in this case coincide with the limits in Lemma \ref{converge0}.
\end{lemma}

\startproof
We prove this first for $i=2$. The $i=3$ case
is essentially the same.  By Corollary \ref{log},
the distance between corresponding
points of $G_2(1+4k)$ and $L_2(1+4k)$ is
at most $4\log(8k)$.
But then, by scaling, the distance between 
corresponding points of 
$\Gamma_2(1+4k)$ and $\Lambda_2(1+4k)$
is at most $4\log_3(8k)(\sqrt 3)^{-k}$, 
a quantity that tends to $0$ exponentially fast.
\endproof

\newpage

\section{The End of the Proof}

\subsection{Scaling Constants}

In this chapter, we prove Statements 3, 4, and 5 of the Main Theorem.

Let $I_2$ and $I_3$ be suitably scaled versions of the sets discussed in
\S \ref{fractal}. 
For $i=2,3$ and $j=1,3$ let
$\Lambda_{ij}$ be the limit of the sequence
$\{\Lambda_i(j+4k)\}$.    To finish the proof of the Main Theorem, we
need to show that
\begin{equation}
\label{biggie}
\Lambda_{21}=\Lambda_{23}=I_2; \hskip 30 pt
\Lambda_{31}=\Lambda_{33}=I_3.
\end{equation}
Our first order of business is to explain how to scale $I_2$ and $I_3$.

\begin{center}
\resizebox{!}{3.5in}{\includegraphics{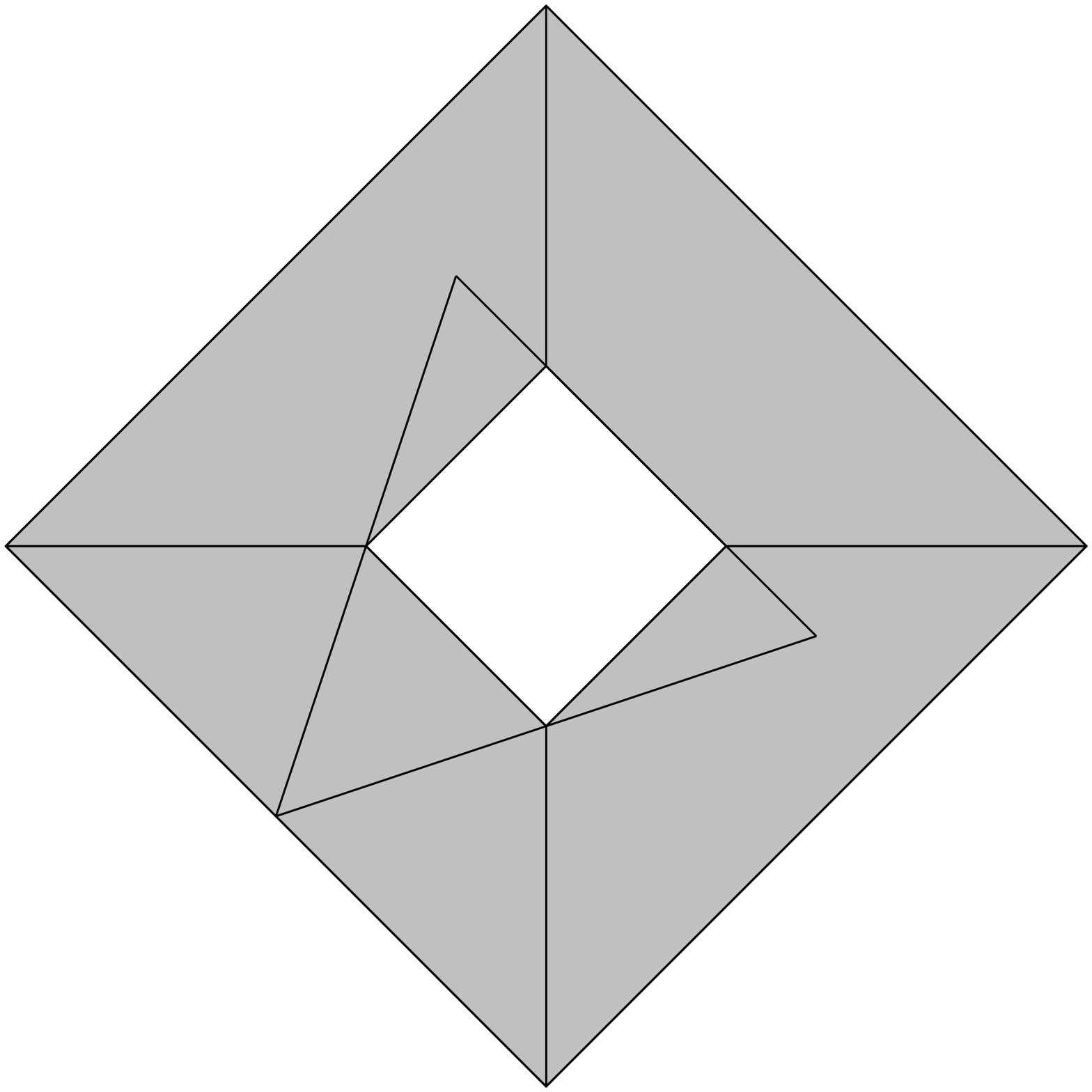}}
\newline
{\bf Figure 10.1:\/} $I_2$ seed and $\Lambda_2(1)$.
\end{center}

We will deal with $I_2$ first.  Figure 10.1 shows the initial seed for $I_2$,
scaled so that it is well situated with respect to $\Lambda_2(1)$.   
By construction, $\Lambda_2(1)$ is an isosceles triangle
with vertices $(0,0)$ and $(3,1)$ and $(1,3)$.   
Given the way that the $I_2$ seed sits with respect to this triangle, we
deduce that the center of $I_2$ is $(3/2,3/2)$ and the left corner is
$(-3/2,3/2)$.    This gives us the copy of $I_2$ mentioned in the
Main Theorem.

Now we deal with $I_3$.  Figure 10.2 shows $\Lambda_2(1)$ and
the union $I_3(1)$ of triangles obtained by applying the substitution
operator to the $I_3$ seed.

\begin{center}
\resizebox{!}{4.5in}{\includegraphics{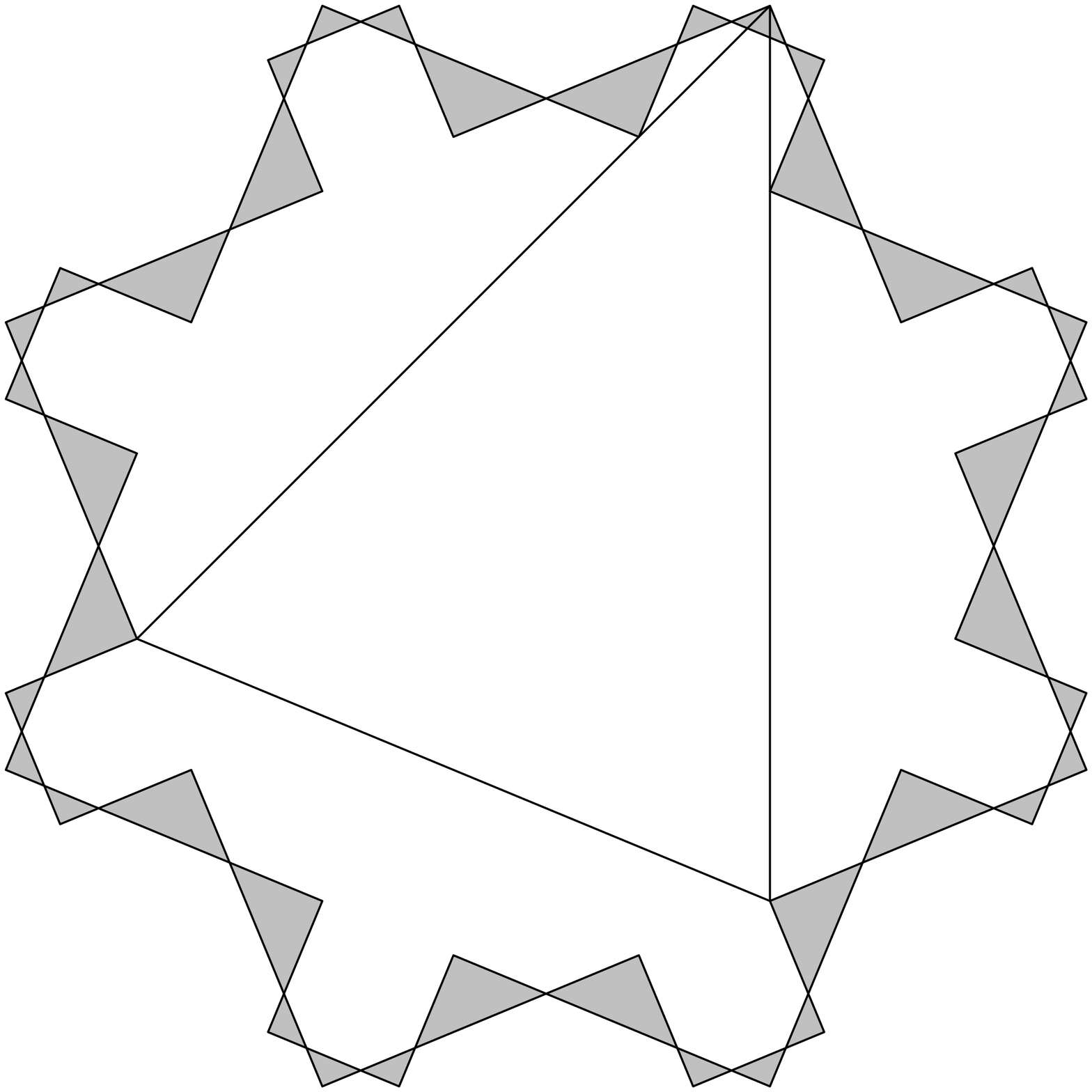}}
\newline
{\bf Figure 10.1:\/} $I_3(1)$ and $\Lambda_3(1)$.
\end{center}

$\Lambda_3(1)$ is an isosceles triangle, with vertices
\begin{equation}
(0,0); \hskip 30 pt
(-2-\sqrt 2,-2-\sqrt 2) \hskip 30 pt
(0,-2-2\sqrt 2).
\end{equation}
The point $(0,0)$ is the apex of the triangle, the top vertex.

This is enough to determine the placement of $I_3$.  
Here is the nicest way we can think of describing the scaling.
Let $s=1+\sqrt 2$.  This is the scaling factor that came up
in the previous chapter. The $I_3$ seed consists in $8$ isosceles
triangles.  One of the isosceles triangles has vertices
$(0,0)$, $(1/2,-s/2)$ and $(-s/2,-1/2)$.   This gives
us the version of $I_3$ mentioned in the Main Theorem.
Note that all the vertices of $I_3$, as we have scaled it, lie
in $\Q[\sqrt 2]$.

\subsection{Proof of Statement 3}

Let $I_2$ and $I_3$ be scaled as in the previous section.
By Equation \ref{fixed}, we see that
$\Lambda_{ij}$ contains the vertices of all the
polygons of which it is a limit.  Moreover, by
scaling,  the vertices of the approximating
polygons become dense in $\Lambda_{ij}$ as
the number of sides tends to $\infty$.
We passed to the limit 
of renormalization precisely to arrange these two
properties.  We want to understand the
sets in Equation \ref{biggie} just from the placement
various finite collections of vertices.

Given the properties enjoyed by the vertices, we can establish
Equation \ref{biggie} simply by showing, for $i=1,2$ and for all
odd $j$ that the vertices of the polygons
approximating $\Lambda_{ij}$ are contained in
$I_i$ and become dense in $I_i$ as $j \to \infty$.

We will deal with $i=3$ first.
Figure 10.3 shows a  piece of the $I_3(3)$, together
with a piece of $\Lambda_3(3)$.  Notice that
all the vertices of $\Lambda_3(3)$ in sight coincide
with right-angled vertices of $I_3(3)$.  Notice also
that there seems to be (most of) a smaller copy of
Figure 10.1 sitting near the top of Figure 10.3.

\begin{center}
\resizebox{!}{4in}{\includegraphics{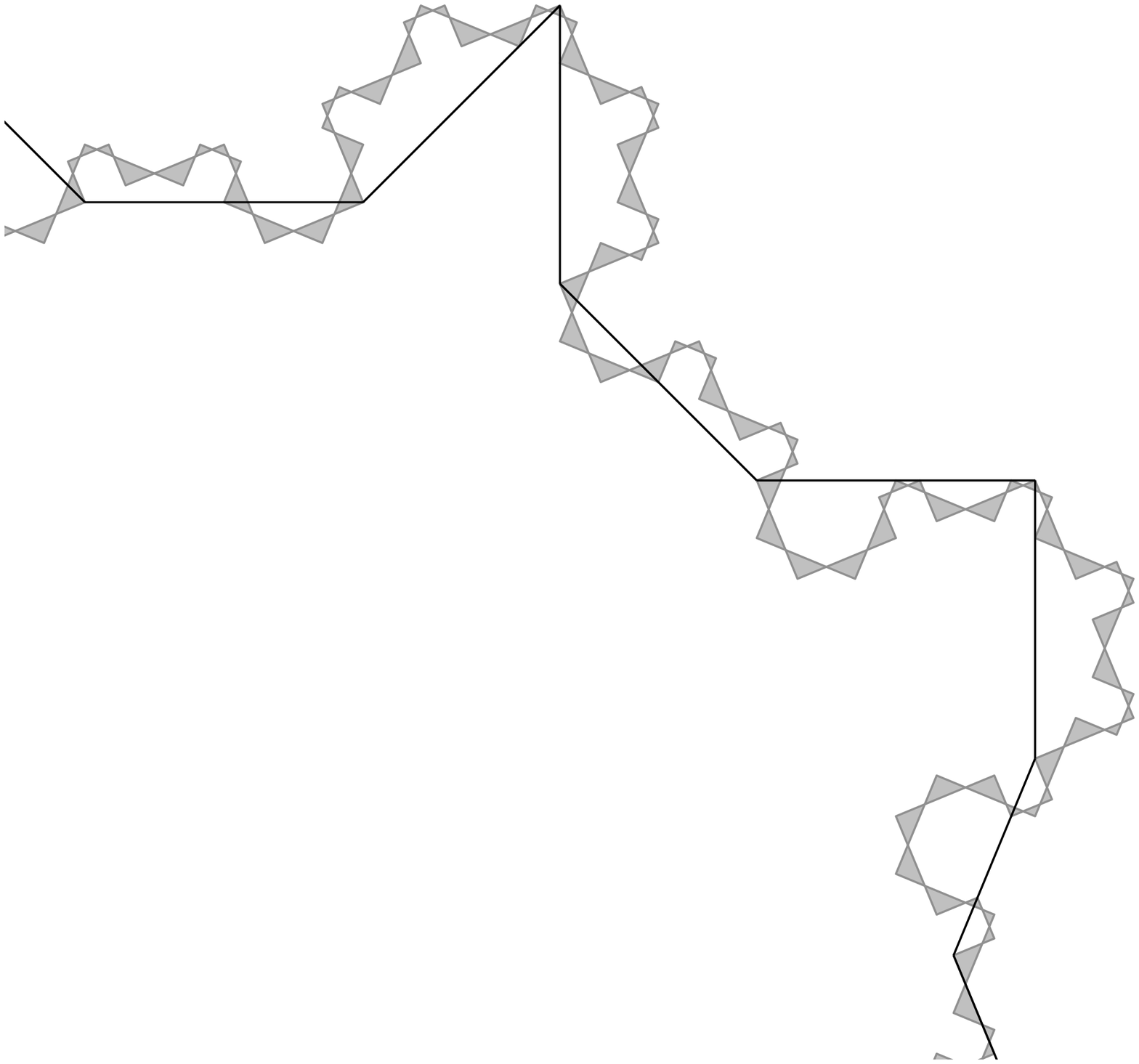}}
\newline
{\bf Figure 10.3:\/} $I_3(1)$ and $\Lambda_3(1)$.
\end{center}

Say that a {\it special vertex\/} of $I_3$ is a vertex of a right-angled
triangle of some $I_3(n)$.   We check by direct (computer)
calculation that
every vertex of $\Lambda_3(j)$ is a special vertex of $I_3$
for $j=1,3,5,7,9$.  This is an easy calculation that takes
place entirely inside $\Q[\sqrt 2]$.  

We can't make an infinite number of explicit calculations to deal
with each successive polygon, so we describe a method,
based on self-similarity, that proves everything we need
from a finite number 
of calculations.     For ease of exposition, we will work
entirely with the sequence $\Lambda_3(1+4k)$.  The
sequence $\Lambda_3(3+4k)$ has a similar treatment.
It is convenient to set $\Lambda_3=\Lambda_{31}$.
So, $\Lambda_3$ is the limit of the sequence we are
considering.

The polygon $\Lambda_3(1)$ is just a triangle, and
too small to be of any use to us (beyond what we have
already done with it above.)  So, we consider
the next polygon in the sequence.
Each edge of $\Lambda_3(5)$ is a scaled copy of some
$\lambda_3(n)$.  We call $n$ the {\it type\/} of the edge.
The type belongs to the
set $S=\{1,...,43\}-\{11,22,33\}$.  
Inspecting the numerical sequence associated
to $\Lambda_3(5)$ -- see \S \ref{sub} -- we check that
$\Lambda_3(5)$ has
every possible type edge.  We check the same
thing for $\Lambda_3(9)$.   

 Suppose that
$A \subset \Lambda_3(5)$ is an edge of
type $n$ and $B \subset \Lambda_3(9)$ is an
edge of type $n$.   The pair $(A,B)$ determines
a homothety 
\begin{equation}
H=H(A,B): A \to B.
\end{equation}
Given the way $\Lambda_3$ is the limit of our sequence, we have
\begin{equation}
H: \Lambda_3(A) \to \Lambda_3(B).
\end{equation}
Here $\Lambda_3(A)$ is the portion of $\Lambda_3$ between the endpoints
of $A$, and likewise for $\Lambda_3(B)$.   So, $H$ is a homothety that
carries a small portion of $\Lambda_3$ to an even smaller portion.
We call $H$ a {\it special homothety\/} associated to the pair $(A,B)$.
Every such pair gives rise to a special homothety.  We call $H$
{\it good\/} if
\begin{equation}
H: I_3(A) \to I_3(B).
\end{equation}
Here $I_3(A)$ is the portion of $I_3$ bounded by the endpoints of $A$.
Likewise we define $I_3(B)$.  This definition makes sense because
we have shown that the endpoints of $A$ and $B$ both lie in $I_3$.
Just to be clear, we mean that $H$ sets up a bijection between
$I_3(A)$ and $I_3(B)$ that maps endpoints to endpoints.

Below, we will prove the following result.
\begin{lemma}[Good Homothety]
Every special homothety is good.
\end{lemma}
The Good Homothety Lemma involves a finite number of special
homotheties.   We will show in the next section how one verifies
that a special homothety is good using a finite calculation.
Indeed, one can practically see this by inspection using our
program OctoMap 2.

Before we prove the Good Homothety Lemma, let us explain 
what it does for us.  Suppose we want to show that
the vertices of $\Lambda_3(13)$ lie in $I_3$.    Let $v$ be
a vertex of $\Lambda_3(13)$.   If $v$ is also a vertex
of $\Lambda_3(9)$, we are already finished.  Otherwise,
there is some edge $B$ of $\Lambda_3(9)$ such that
$v$ is one of the vertices associated to the replacement of
$B$.  That is, $v$ lies on the polygonal path of
$\Lambda_3(13)$ connecting the endpoints of $B$.
Choose and edge $A$ of $\Lambda_3(5)$ that has
the same type as $B$.   

The special homothety
$H=H(A,B)$ carries $A$ to $B$.  
Moreover, $H$ carries $\Lambda_3(9;A)$ to
$\Lambda_3(13;B)$.  Here
$\Lambda_3(9;A)$ denotes the portion of
$\Lambda_3(9)$ subtended by $A$.
At the same time as this, $H$ carries
$I_3(A)$ to $I_3(B)$.   But
$H^{-1}(v)$ is a vertex of
$\Lambda_3(9;A)$.   Hence $H^{-1}(v) \in I_3(A)$.
Hence $v \in I_3(B)$.  This proves that
all vertices of $\Lambda_3(13)$ lie in $I_3$.

Continuing in this way, we see that
all vertices of $\Lambda_3(1+4k)$ lie in
$I_3$ for all positive $k$.   To take
care of density, we let
$d_k$ denote the maximum distance
from a point of $I_3$ to a vertex of
$\Lambda_3(1+4k)$.    Since our
special homotheties scale distances by
a factor of $\sqrt 2 -1$, we see that
\begin{equation}
d_{k+1}=(\sqrt 2-1) d_k.
\end{equation}
In particular $d_k \to 0$ as $k \to \infty$.
This proves that the vertices of
$\Lambda_3(1+4k)$ becomes dense
in $I_3$ as $k \to \infty$.

In short, the Good Homothety Lemma implies Statement 3
of the Main Theorem.

\subsection{Proof of the Good Homothety Lemma}
\label{goodhomo}

In \S \ref{snowflake} we described the basic substitution rule
for the snowflake, in which each isosceles triangle is replaced
by $5$ smaller ones.   We started with a seed, consisting
$8$ isosceles triangles arranged in a cyclic pattern, and
then iteratively replaced each triangle by the smaller ones.

Here we describe a more conservative approach to the construction.
Let $L_1$ and $L_2$ be finite cyclically ordered lists of
isosceles triangles.   We write $L_1 \to L_2$ if one obtains
$L_2$ from $L_1$ by replacing one of the triangles of $L_2$
by the $5$ smaller ones from the substition rule and
adjusting the cyclic ordering in the obvious way.
We call a list $L$ of triangles {\it admissible\/} if we have
a finite chain $I_3(0)=L_0 \to ... \to L_n$.
The specific lists $I_3(n)$, produced in \S \ref{snowflake} are
examples of admissible lists.

We say that a {\it special vertex\/} of $I_3$ is a right-angled
vertex of some $I_3(n)$.  An equivalent definition is that a
special vertex of $I_3$ is a right-angled vertex of some
admissible list.   Let $\{v_k\}$ be a finite list of
special vertex.  We say that an admissible list
$L$ is {\it compatible\/} with $\{v_k\}$ if each point
$v_k$ is a vertex of some triangle in $L$.
We say that $L$ is {\it the minimal pattern\/} for
$\{v_k\}$ if $L$ is compatible with $\{v_k\}$ and
there is no $L'$ such that $L' \to L$ and
$L'$ is compatible with $\{v_k\}$.
\newline
\newline
{\bf Remark:\/}
For our proof of the Good Homotopy Lemma, we do not
need to know that the minimal pattern associated to
a finite list of vertices is unique, but we emphasize
that the minimial pattern is unique.  The uniqueness
comes from the  tree-like nature of the subdivision rule.
We simply start with the seed and replace
triangles only when necessary.  After a finite
number of steps we arrive at the minimial pattern.
\newline

Let $A$ be a line segment whose endpoints are special vertices
of $I_3$.   We let $J_3'(A)$ be the minimal pattern associated to the endpoints
of $A$.   We don't care so much about the whole pattern, but only
the portion related to $A$.  We let $J_3(A)$ denote the portion of $J_3'(A)$ that
lies between the two endpoints of $A$.  Figure 10.4 shows
the minimal pattern $J_3(A)$ when $A$ is a certain type-1
edge of $\Lambda_3(5)$.

\begin{center}
\resizebox{!}{2.3in}{\includegraphics{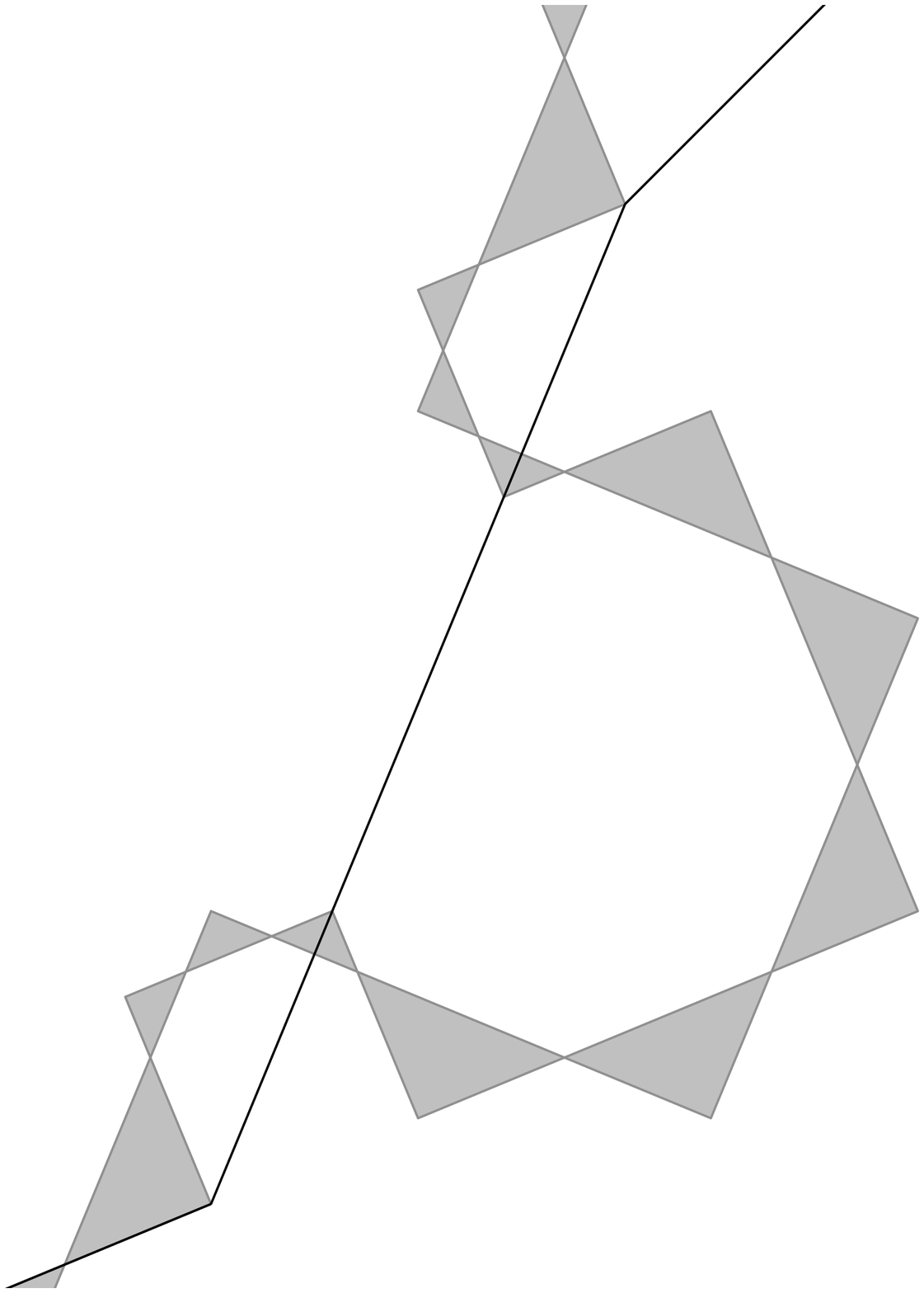}}
\newline
{\bf Figure 10.4:\/} $J_3(A)$ for $A$ a type $1$ edge of $\Lambda_3(5)$.
\end{center}

The figure we have shown is completely representative of the
local picture for the type 1 edges of $\Lambda_3(5)$.  The
associated minimal patterns are all the same.  This is not really
so surprising, because these vectors all are translates of each
other.    We also check that the local picture looks the
same for all type 1 edges of $\Lambda_3(9)$.  That is,
the picture looks the same up to a scale factor of
$\sqrt 2-1$.

So, suppose that $A$ is a type 1 edge of $\Lambda_3(5)$ and
$B$ is a type 1 edge of $\Lambda_3(9)$.   Then $H(A,B)$
carries $J_3(A)$ to $J_3(B)$ because the patterns are identical
up to scale.  But $J_3(A)$ and $J_3(B)$ respectively determine
$I_3(A)$ and $I_3(B)$.   Hence $H$ carries $I_3(A)$ to
$I_3(B)$, as desired.  This proves the Good Homotopy Lemma
for edges of type 1.  The key idea is that the minimal
pattern $J_3(A)$ only depends on the type of $A$, up
to scale.

\begin{center}
\resizebox{!}{1.7in}{\includegraphics{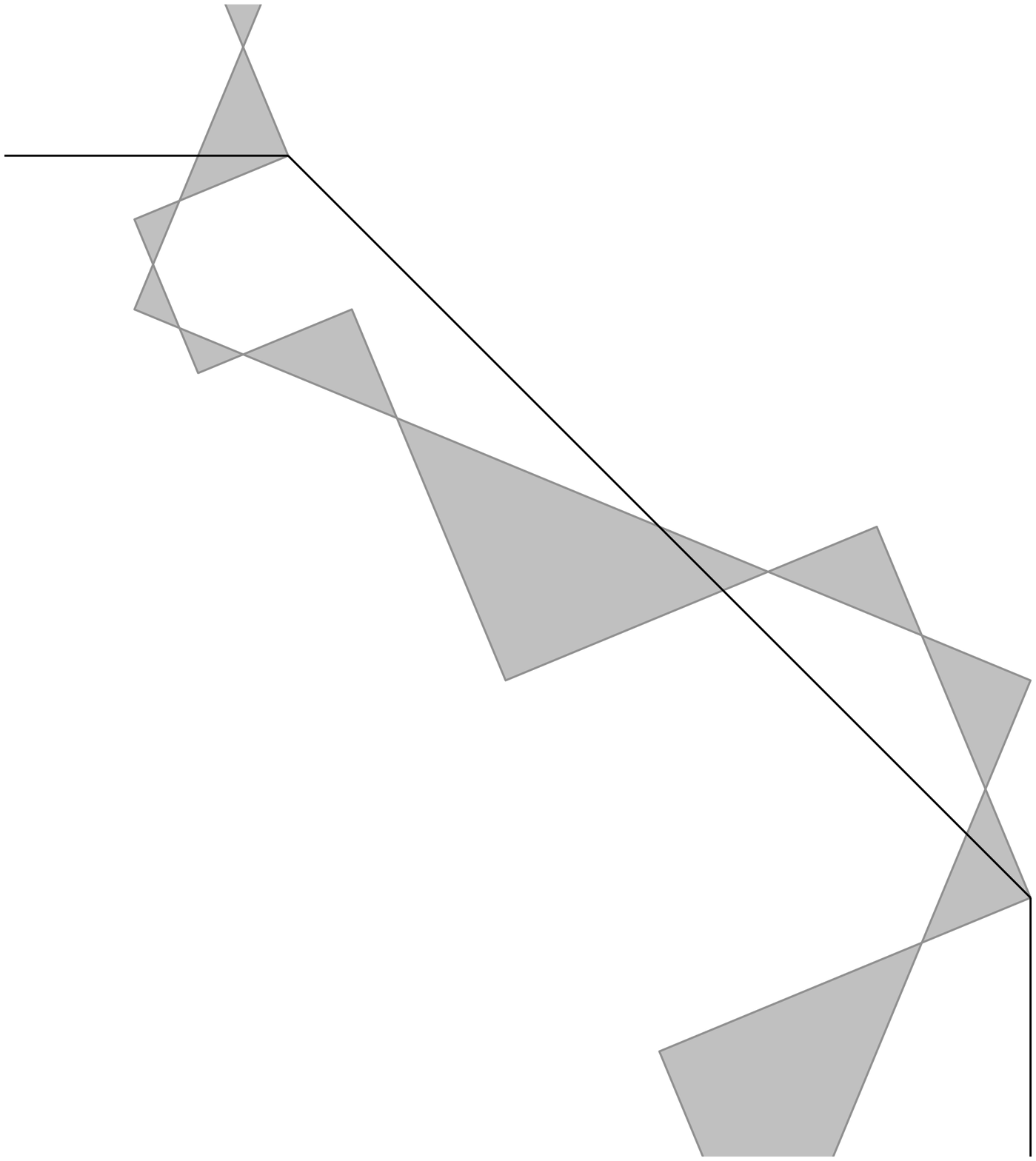}}
\resizebox{!}{1.7in}{\includegraphics{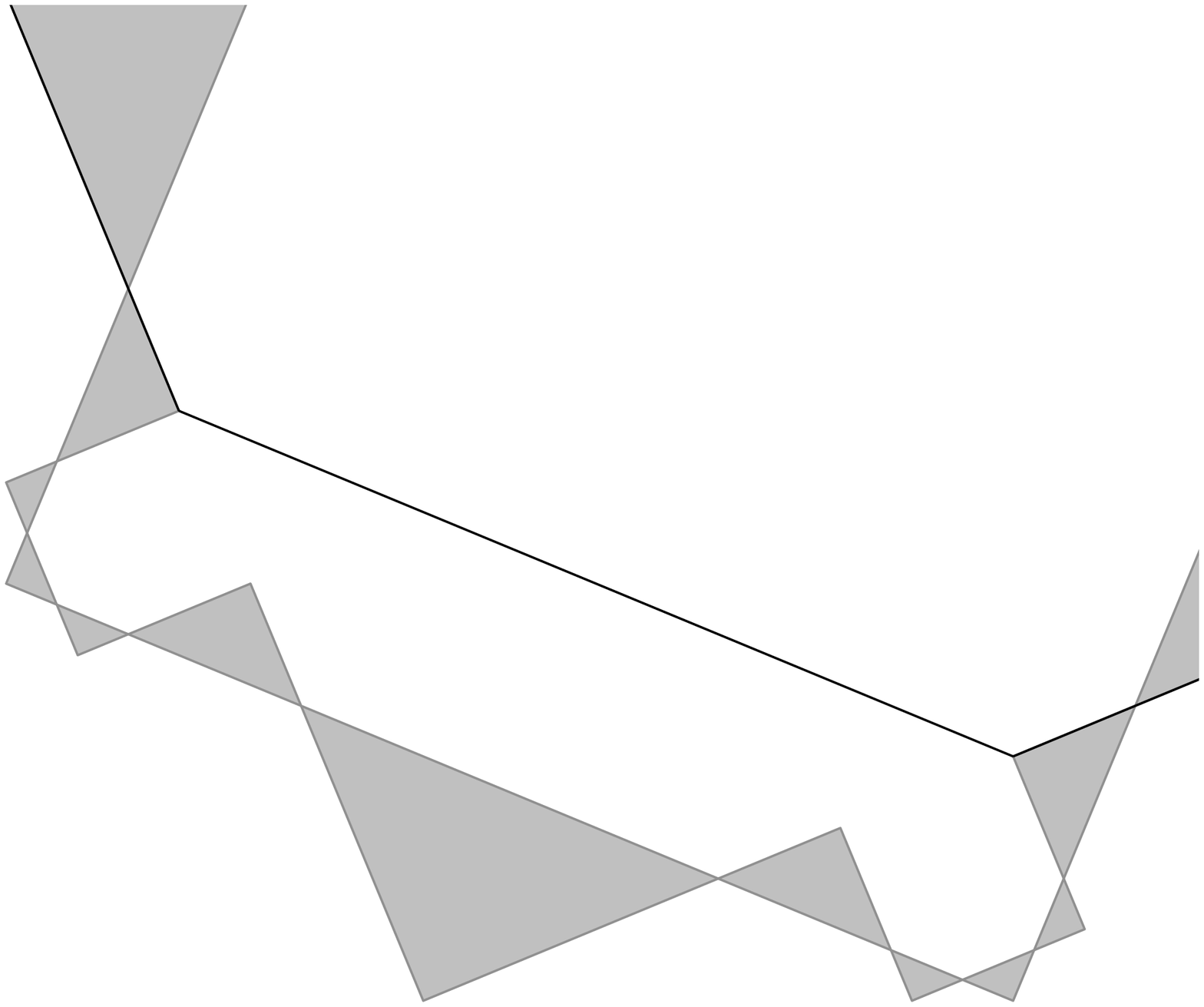}}
\resizebox{!}{1.7in}{\includegraphics{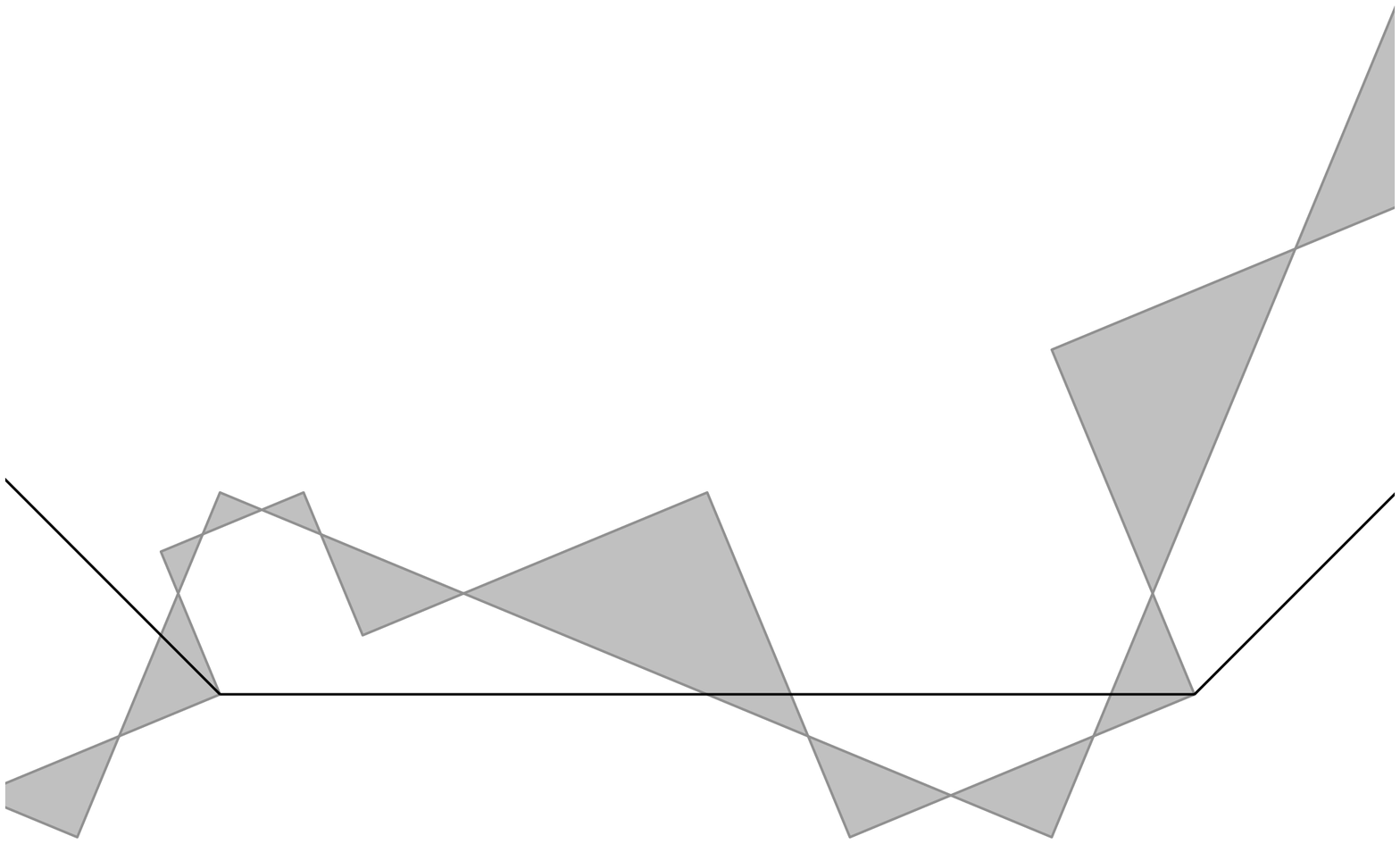}}
\newline
{\bf Figure 10.5:\/} The minimal patterns associated to types $2,3,4$.
\end{center}

Figure 10.5 shows the minimal patterns associated to the
types $2$, $3$, and $4$.   All the edge types except those
of the form $6+11k$ and $9+11k$ have the same
property we have already discussed.  The local picture is
the same, independent of the edge.   By inspection,
we see that the Good Homotopy Lemma holds for
for all edge types except the $8$ types we have mentioned.

For each of the remaining edge types, there are two local
pictures.  up to rotation, these two local pictures are the same for each
of the edge types.  Figure 10.6 shows the two pictures
associated to the edges of type $6+11k$.  The edges
are vertical.

\begin{center}
\resizebox{!}{2.5in}{\includegraphics{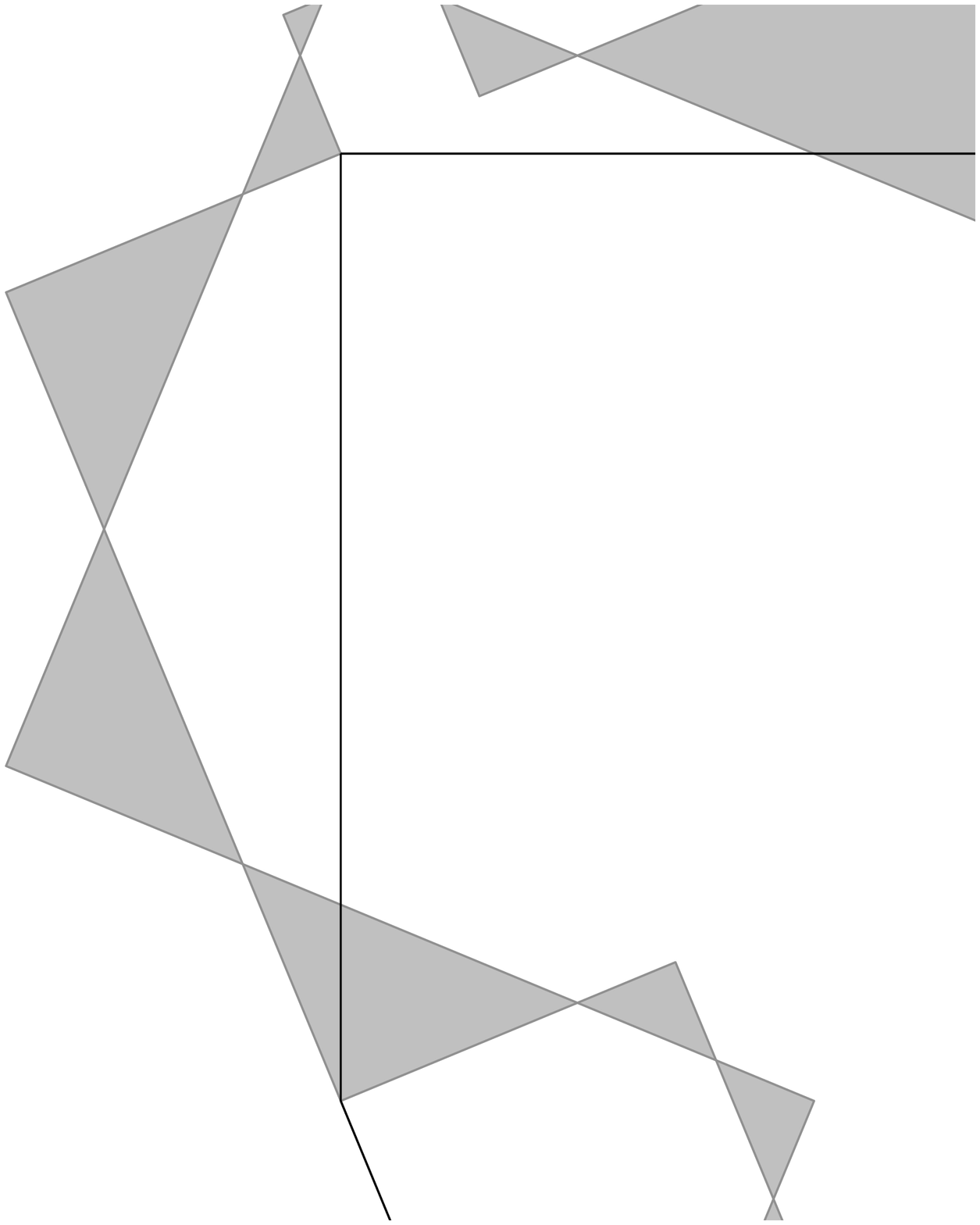}}
$\hskip 10 pt$
\resizebox{!}{2.5in}{\includegraphics{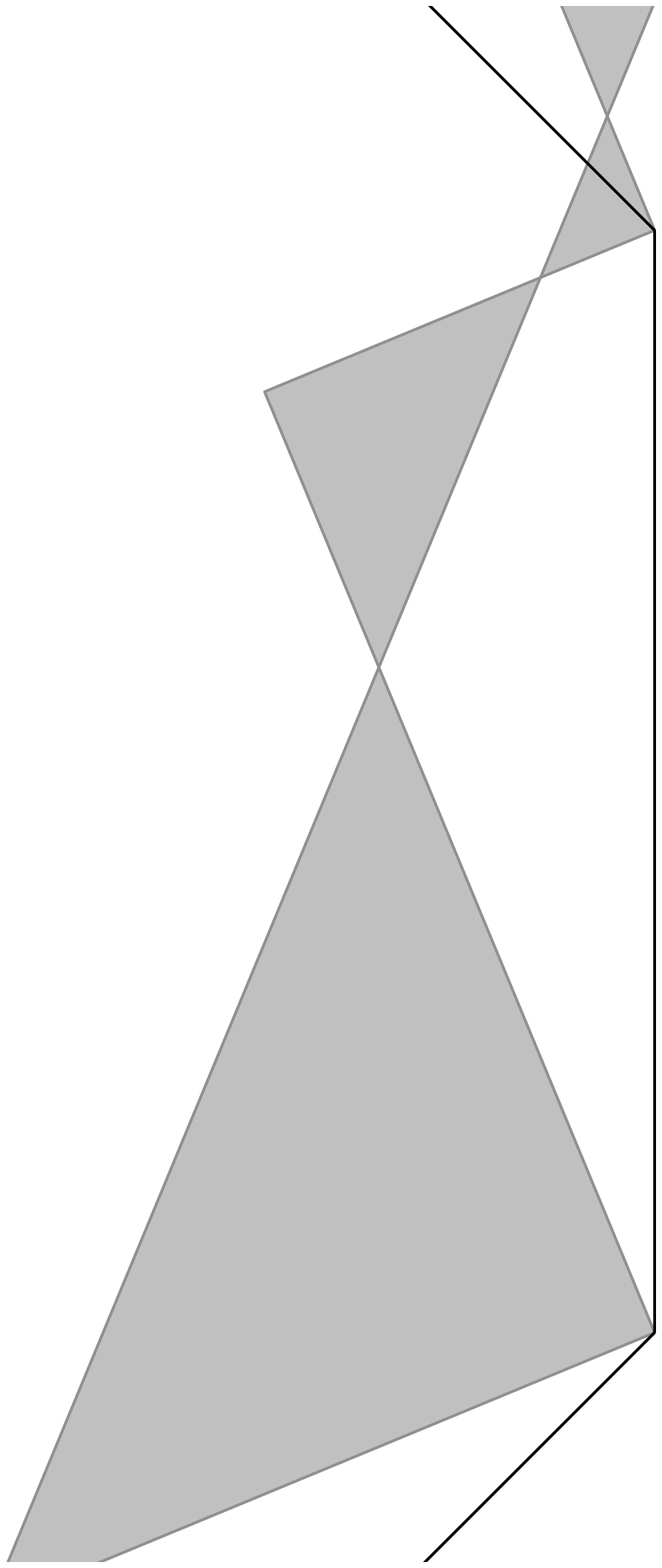}}
\newline
{\bf Figure 10.6:\/} The minimal patterns associated to the edges
of type $6+11k$.
\end{center}

Let $A$ and $A'$ be the two type 6 edges shown in Figure 10.6.
The pattern $J_3(A)$ consists of $4$ triangles $T_1 \cup T_2 \cup T_3 \cup T_4$
and $J_3(A')$ consists of $3$ triangles $T_1' \cup T_2' \cup T_3'$.
We list these triangles top to bottom.  There is a translation $\tau$ 
such that $\tau(A)=A'$. To finish our proof for
edges of type 
$6+11k$, we just have to see that $\tau$
carries $I_3(A)$ to $I_3(A')$.   Note that $\tau$ 
carries $T_j'$ to $T_j$ for $j=1,2$.  So, things
work out as desired for the ``top portions''
of the sets of interest to us.

Now we consider the bottom portions.
Let $S$ be the line segment on the left hand side
of Figure 10.6 that connects the bottom vertex of $A$ to
the bottom vertex of $T_2$.  Note that $S$ is made
from the hypotenuse of $T_3$ and a short side
of $T_4$.    Let $S'$ be the edge of $T_3$ that connects the
bottom of $A'$ to a vertex of $T_2'$.  We see
by direct calculation that $\tau(S)=S'$.  The
point is that $S$ and $S'$ have the same slope and length.
We just have to prove that $\tau$ carries $I_3(S)$ to $I_3(S')$.
But we have exactly the situation discussed in
\S \ref{hidden}, and the desired result here
follows from Lemma \ref{hiddenlemma}.

it remains to consider
The edges of type $9+11k$.  It turns out that the minimal
patterns associated to these edges look exactly like the ones
associated to the edges of type $6+11k$, except that the picture has
been rotated $45$ degrees.    The same argument as in the
previous case works here.

This completes the proof of the Good Homothety Theorem.
Again, the Good Homothety Theorem (and our previous work)
implies Statement 3 of the Main Theorem

\subsection{Proof of Statement 4}

In considering $I_2$ we make all the same definitions
that we made for $I_3$, except that we redefine what
we mean by a special point.  In the case of
$I_3$, the special points are the right-angled
vertices of the triangles.  In the case of $I_2$,
we call a point {\it special\/} if it has one of two properties:

\begin{itemize}
\item It is the center point of a parallelogram of some admissible list.
\item It is the center of the long edge of a trapezoid of some admissible list.
\end{itemize}
We check directly that every vertex of
$I_2(j)$ is a special point, for the initial values
$j=1,3,5,7,9$.  This check entirely involves calculations
with rational points.  We want to extend this result for
all odd $j$.    As in the $\Lambda_3$ case,
we ignore the sequence $\Lambda(4k+3)$ and
concentrate on the sequence $\Lambda(4k+1)$.
The other sequence has essentially the same treatment.

To finish the proof, we just have 
to prove the analogue of the Good Homothety Lemma
for $I_2$.   Given the combinatorially identical
subdivision rules that produce $I_2$ and $I_3$, it
makes sense to say that a minimal pattern
associated to $I_2$ is isomorphic to a minimal
pattern associated to $I_3$:  The canonical
bijection between the $I_3$-shapes and the
$I_3$-shapes carries the one pattern to the other.
For $j=1,3,5,7,9$ we make the following observation:
The minimial pattern associated
to some edge of $\Lambda_2(j)$ is isomorphic to
the minimal pattern associated to the corresponding
edge of $\Lambda_3(j)$.  There is a perfect
correspondence between the two kinds of
patterns.   

Given the combinatorial isomorphism between
the $I_3$-patterns and the $I_2$-patterns, the proof of
the Good Homothety Lemma goes through, word for word
for all the edges except those of type $6+11k$ and
$9+11k$.    In these cases, there are two combinatorial
types of minimial pattern, and we must make a
separarate analysis.  Figure 10.7 shows the relevant
portions of the minimal patterns associated to
edges of type $6+11k$.

\begin{center}
\resizebox{!}{2.4in}{\includegraphics{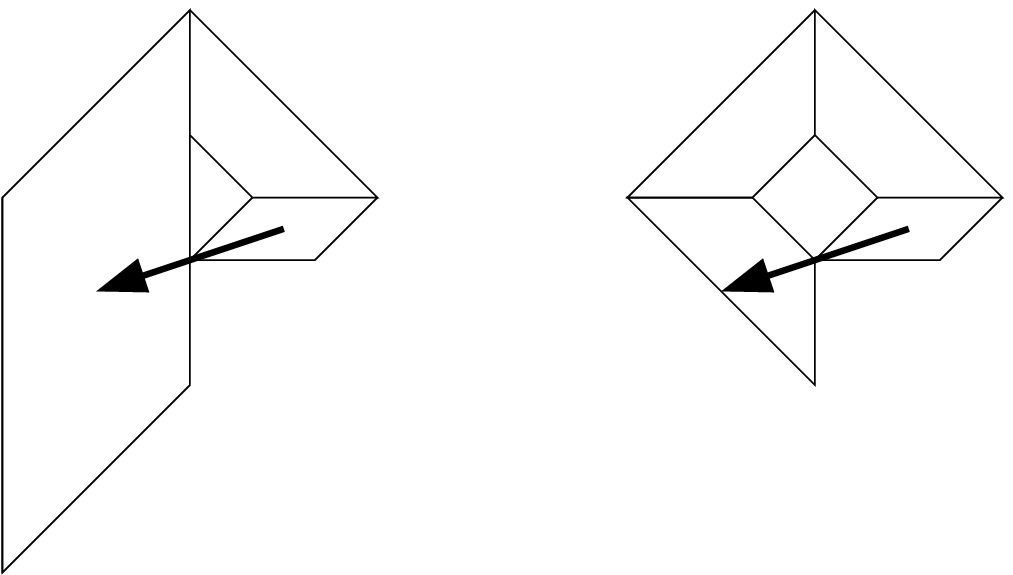}}
\newline
{\bf Figure 10.7:\/} The minimal patterns associated to
the edges of type $6+11k$.
\end{center}

On the left hand side we have a chain of $3$ quadrilaterals
and on the right hand side we have a chain of $4$ quadrilaterals.
The first two pieces agree on both sides.  When we chop
two pieces off, and superimpose what remains of the
left and right sides, we get exactly the left hand side
of Figure 4.7, the figure we used to illustrate the
hidden symmetry of the carpet.
The rest of the proof is just like what we did for the
$\Lambda_3$ case, except that we use
Lemma \ref{hiddenlemma2}in place of
Lemma \ref{hiddenlemma}.  

The edges of type $9+11k$ have the same treatment.
This completes the proof of the Good Homotopy Lemma
in the $I_2$ case, and thereby completes the proof
of Statement 4 of the Main Theorem.

\subsection{Proof of Statement 5}

Statement 5 takes a bit of unwrapping.   We start with
an arbitrary sequence $\{x_n\}$ of good points in the
half-strip $\Sigma_0^1$, and we let $A_n$ be the
arithmetic graph corresponding to $x_n$.   We want to
understand the map 
\begin{equation}
\phi_n:  S_{3,n}\pi_3(A_n) \to S_{2,n} \pi_2(A_n)
\end{equation}
Here $S_{2,n}$ and $S_{3,n}$ are the dilations
discussed in the Main Theorem.

In light of the analysis done in the previous chapter,
we can assume that
\begin{equation}
S_{k,n} \pi_k(A_n) = \Gamma_k(n).
\end{equation}
Here $n$ is odd, and $\Gamma_k(n)$ is as in 
\S\ref{gammadefine}.  For ease of exposition, we now make
the blanket assumption that $n \equiv 1$ mod $4$.
The case when $n \equiv 3$ mod $4$ has the same
treatment.

The two paths $\Gamma_2(n)$ and $\Gamma_3(n)$ have
the same number of sides, and so there is a canonical
map 
\begin{equation}
\phi_n: \Gamma_3(n) \to \Gamma_2(n).
\end{equation}

There is a canonical map from
$\Gamma_k(n)$ to $\Lambda_k(n)$, and this gives us a map
\begin{equation}
\psi_n: \Lambda_3(n) \to \Lambda_2(n).
\end{equation}
Our strategy is to first analyze the sequence $\{\psi_n\}$ of
maps and then use Lemma \ref{converge0} to get the
desired information about $\{\phi_n\}$.

There is a sense in which $\psi_n$ and $\phi$ are compatible,
we we now explain.
The map $\phi$ is the limit of a sequence of canonical
bijections between the admissible $I_3$-patterns
of triangles and the admissible $I_2$-patterns of
quadrilaterals.   We call this bijection
the {\it shape bijection\/}.   We denote it by
$\Phi$.  Note that $\Phi$ is canonical only
up to composition with isometries of $I_2$ and $I_3$.
This accounts for the isometry $F$ that appears in
Statement 5 of the Main Theorem.

Given an edge $e_3$ of $\Lambda_3(n)$,
there are two triangles $T_1$ and
$T_2$ in the minimal pattern
$\Lambda_3(e_3)$ whose right-angled
vertices are the endpoints of $e_3$.
Let 
$$e_2=\psi_n(e_3)$$ be the edge of $\Lambda_2$ corresponding
to $e_3$.    We say that the pair $(e_3,e_2)$
is {\it compatible\/} if the endpoints
of $e_2$ are special points of the two
shapes $\Phi(T_1)$ and $\Phi(T_2)$. 
In case $\Phi(T_1)$ is a parallelogram,
the relevant endpoint of $e_2$ should
be the center point of this parallelogram.
In case $\Phi(T_1)$ is a trapezoid,
the relevant endpoint of $e_2$
should be the midpoint of the long side
of this trapezoid.

We check that all relevant edge-pairs
are compatible for $n=1,3,5,7,9$.  But our
proofs of the Good Homotopy Lemma, in
the two cases, involved combinatorially
identical sets in each case.  That is,
whatever we did in the $I_3$ case
carries over to the $I_2$ case {\it via\/}
the shape bijection.  From this we see,
by induction and symmetry, that
all relevant edge-pairs are compatible
for all odd positive integers $n$.

Given this notation of compatibility, we have
\begin{equation}
\label{haus}
\lim_{n\to \infty} \psi_n =F \circ \phi,
\end{equation}
for some isometry $F$.
The convergence takes place
in the sense that the graph of $\psi_n$ converges,
to the graph of $\phi$ in the
Hausdorff topology.    In particular, the
left hand side of Equation \ref{haus} does have
a limit.  But Lemma \ref{s1proof} now says that
\begin{equation}
\lim_{n \to \infty} \phi_n=\lim_{n\to \infty} \psi_n.
\end{equation}
Putting the two equations together gives us Statement 5 of
the Main Theorem.

This completes our proof of the Main Theorem.

\newpage

\section{References}

[{\bf BC\/}] N. Bedaride, J. Cassaigne, {\it Outer Billiards outside Regular Polygons\/},
preprint (2009)
\newline
\newline
[{\bf DF\/}] D. Dolyopyat and B. Fayad, {\it Unbounded orbits for semicircular
outer billiards\/}, Annales Henri Poincar\'{e}, to appear.
\newline
\newline
[{\bf G\/}] A. Graham, {\it Nonnegative matrices and Applicable Topics in Linear Algebra\/},
John Wiley and Sons, New York (1987).
\newline
\newline
[{\bf GS\/}] E. Gutkin and N. Simanyi, {\it Dual polygonal
billiard and necklace dynamics\/}, Comm. Math. Phys.
{\bf 143\/}:431--450 (1991).
\newline
\newline
[{\bf Ko\/}] Kolodziej, {\it The antibilliard outside a polygon\/},
Bull. Pol. Acad Sci. Math.
{\bf 37\/}:163--168 (1994).
\newline
\newline
[{\bf M1\/}] J. Moser, {\it Is the solar system stable?\/},
Math. Intel. {\bf 1\/}:65--71 (1978).
\newline
\newline
[{\bf M2\/}] J. Moser, {\it Stable and random motions in dynamical systems, with
special emphasis on celestial mechanics\/},
Annals of Math Studies 77, Princeton University Press, Princeton, NJ (1973).
\newline
\newline
[{\bf N\/}] B. H. Neumann, {\it Sharing ham and eggs\/},
Summary of a Manchester Mathematics Colloquium, 25 Jan 1959,
published in Iota, the Manchester University Mathematics Students' Journal.
\newline
\newline
[{\bf S1\/}] R. E. Schwartz, {\it Unbounded Orbits for Outer Billiards\/},
J. Mod. Dyn. {\bf 3\/}:371--424 (2007). 
\newline
\newline
[{\bf S2\/}] R. E. Schwartz, {\it Outer Billiards on Kites\/},
Annals of Math Studies {\bf 171\/}, Princeton University Press (2009)
\newline
\newline
[{\bf S3\/}] R. E. Schwartz, {\it Outer Billiards and the Pinwheel Lemma\/},
preprint (2009)
\newline
\newline
[{\bf T1\/}] S. Tabachnikov, {\it Geometry and billiards\/},
Student Mathematical Library 30,
Amer. Math. Soc. (2005).
\newline
\newline
[{\bf T2\/}] S. Tabachnikov, {\it On the Dual Billiard Problem\/},
Advances in Mathematics {\bf 115\/} pp 221-249 (1995)
\newline
\newline
[{\bf VL\/}] F. Vivaldi and J. H. Lowenstein, {|it Arithmetical properties of a family
of irrational piecewise rotations\/}, {\it Nonlinearity\/} {\bf 19\/}:1069--1097 (2007).
\newline
\newline
[{\bf VS\/}] F. Vivaldi and A. Shaidenko, {\it Global stability of a class of discontinuous
dual billiards\/}, Comm. Math. Phys. {\bf 110\/}:625--640 (1987).

\end{document}